\documentclass[12pt, draft]{article}
\usepackage{amsmath}
\usepackage{amssymb}
\usepackage{latexsym}
\allowdisplaybreaks[4]
\date{}
\renewcommand{\theequation}{\arabic{section}.\arabic{equation}}
\newtheorem{Thm}{Theorem}[section]
\newtheorem{Prop}{Proposition}[section]
\newtheorem{Lemma}{Lemma}[section]
\newtheorem{Cor}{Corollary}[section]
\newtheorem{Def}{Definition}[section]
\newtheorem{Remark}{Remark}[section]
\newtheorem{Example}{Example}[section]
\renewcommand{\Bbb}{\bf}

\newcommand{\R}{{\mathbb R}}

\def\v#1{\vert #1 \vert}

\def\V#1{\Vert #1 \Vert}

\def\noytilde#1#2
{\left\| \{\tilde\rho ,\tilde v\} \right\|_{(s_0; #1, #2)}}
\title{An inverse problem for a three-dimensional heat equation in thermal imaging 
and the enclosure method}
\author{Masaru IKEHATA\footnote{
Laboratory of of Mathematics,
Institute of Engineering,
Hiroshima University, Higashi Hiroshima 739-8527, JAPAN.
e-mail address: ikehata@amath.hiroshima-u.ac.jp}
and
Mishio KAWASHITA\footnote
{Department of Mathematics,
Graduate School of Sciences,
Hiroshima University, Higashi Hiroshima 739-8526, JAPAN.
e-mail address: kawasita@math.sci.hiroshima-u.ac.jp}
}
\date{}
\begin{document}

\maketitle
\begin{abstract}
This paper studies a prototype of inverse initial boundary value
problems whose governing equation is the heat equation in three
dimensions.  An unknown discontinuity embedded in a three-dimensional
heat conductive body is considered.
A {\it single} set of the temperature and heat flux on the lateral boundary
for a fixed observation time is given as
an observation datum.  It is shown that this datum yields
the minimum length of broken paths that start at a given point 
outside the body, go to a point on the boundary of the unknown discontinuity
and return to a point on the boundary of the body
under some conditions on the input heat flux, the unknown discontinuity and the body.
This is new information obtained by using enclosure method.
\end{abstract}


\vskip 2truepc

\setcounter{section}{0}
\section{Introduction}\label{Introduction}

Let $\Omega$ be a bounded domain of ${\Bbb R}^3$ with 
$C^{2,\alpha_0}$ boundary and $0<\alpha_0\le 1$.
Let $D$ be an open subset of $\Omega$ with $C^{2,\alpha_0}$ boundary and satisfy that:
$\overline D\subset\Omega$; $\Omega\setminus\overline D$ is connected.
We denote by $\nu_x$, $\nu_y$ the unit outward normal vectors at $x\in\partial D$,
$y\in\partial\Omega$ on $\partial D$, $\partial\Omega$, respectively.
Let $T$ be an arbitrary {\it fixed} positive number and $\rho=\rho(x)\in C^{0,\alpha_0}(\partial D)$.
Given $f\in L^1(0,T;H^{-1/2}(\partial\Omega))$ let $u=u(x,t)$ be the weak solution of
\begin{equation}
\left\{
\begin{array}{l}
\displaystyle u_t-\triangle u=0\,\,\text{in}\,(\Omega\setminus\overline D)\times\,]0,\,T[,\\[2mm]
\displaystyle \frac{\partial u}{\partial\nu}+\rho u=0\,\,\text{on}\,\partial D\times]0,\,T[,\\[3mm]
\displaystyle
\frac{\partial u}{\partial\nu}=f\,\,\text{on}\,\partial\Omega\times\,]0,\,T[,
\\[2mm]
\displaystyle
u(x,0)=0\,\,\text{in}\,\Omega\setminus\overline D.
\end{array}
\right.
\label{1.1}
\end{equation}
For detailed information about the weak solution which follows \cite{DL},  see subsection 1.5 in this paper.

This paper is concerned with the following problem.
\vskip 0.5em
{\bf\noindent Inverse Problem.}  {\it Fix} $T>0$.  Assume that both $D$ and $\rho$ are {\it unknown}.
Extract information about the location and shape
of $D$ from the temperature $u$ on $\partial\Omega$ over finite time interval $]0,\,T[$
with a {\it fixed known} heat flux $f$.
\vskip 0.5em

\noindent This is a prototype of several inverse problems related
to thermal imaging, dynamical remote sensing and very important one.
$D$ is a mathematical model of unknown discontinuity embedded in a three-dimensional
heat conductive body.
There are extensive mathematical studies of uniqueness and stability issues 
of Inverse Problem. 
In particular, it is known that the observation data uniquely
determine general $D$ itself under a suitable condition on the
heat flux on $\partial\Omega$ in the case when $\rho\equiv 0$. See
Bryan-Caudill \cite{BC}, Canuto-Rosset-Vessella \cite{CRV}, Vessella \cite{Ve}
and his survey paper \cite{VS} together with references therein for more
information about these issues.

\subsection{An interpretation of previous one-space dimensional result}

In \cite{I4} Ikehata started a study that seeks an analytical and constructive
approach for the inverse problem.
He considered a one-space dimensional version of the problem and related ones.
The method used therein is called the {\it enclosure method} which was introduced
by himself in \cite{I1, I2}.
The enclosure method aims at extracting a domain that {\it encloses} an unknown
discontinuity, such as inclusion, cavity or crack in a known background medium
by observing a signal propagating inside the medium on the boundary of the surface surrounding the medium.
Then the Dirichle-to-Neumann map associated with the governing equation of the used signal
appears as an idealized mathematical model of the observed data.
The enclosure method constructs the so-called the {\it indicator function}
by using the Dirichlet-to-Neumann map or its partial knowledge combined with the {\it complex geometrical optics solution}
of the governing equation.  The indicator function has an independent variable 
which is contained in the complex geometrical optics solution as a large parameter.
The complex geometrical optics solution changes its {\it growing} and {\it decaying} property as the parameter goes to infinity
bordering on, for example, a plane in three dimensions.
The behaviour of the indicator function as the independent variable goes to infinity
depends on the relative position of the plane to unknown discontinuity
and enables us to obtain an enclosing domain.
In this sense this original enclosure method can be considered as a method of using the complex geometrical optics solutions.
However, note that the way of using this growing and decaying character positively differs from the well known method 
which goes back to Calder\'on \cite{C} and Sylvester-Uhlmann \cite{SU}
since their method is based on the oscillating character of the complex geometrical optics solutions
about the parameter.

Now let us describe one of the problems considered in \cite{I4}.
Let $u=u(x,t)$ with $u_x(0,t)\in L^2(0,\,T)$ be a solution of the problem 
\begin{equation}
\left\{
\begin{array}{l}
\displaystyle
u_t=u_{xx}\,\,\text{in}\,]0,\,a[\times]0,\,T[,
\\[2mm]
\displaystyle
u_x(a,t)+\rho u(a,t)=0\,\,\text{for}\,t\in\,]0,\,T[,
\\[2mm]
\displaystyle
u(x,0)=0\,\,\text{in}\,]0,\,a[.
\end{array}
\right.
\label{1.2}
\end{equation}
It is assumed that both constants $a > 0$ and $\rho \in {\Bbb R}$ in
(\ref{1.2}) are {\it unknown}.
He considered the problem: extract $a$ from $u(0,t)$ and
$u_x(0,t)$ for $0<t<T$. 
This inverse problem is the one dimensional version of our inverse problem
for (\ref{1.1}). In (\ref{1.2}), sets $]0,\,\infty[$, $]a,\,\infty[$ and $\{a\}$
correspond to $\Omega$, $D$ and $\partial D$ respectively.

In \cite{I4}, to extract $a$ from $u(0,t)$ and
$u_x(0,t)$ $(0 < t < T)$, 
he introduced an indicator function $I(\tau)$ of
independent variable $\tau>0$ given by the integral
$$
\displaystyle
I(\tau)
=\int_0^{T}\left(-v_x(0,t)\,u(0,t)
+u_x(0,t)\,v(0,t)\right)dt,
$$
where $v=v(x,t)$ is a solution of the one-dimensional backward heat equation
$v_t+v_{xx}=0$ of the following form:
$$\displaystyle
v(x,t)=e^{-\tau^2 t}e^{-\tau x}.
$$
For this indicator function $I(\tau)$, he showed that the formula
\begin{equation}
\displaystyle
\lim_{\tau\longrightarrow\infty}
\frac{1}{\tau}\log\vert I(\tau)\vert
=-2a
\label{1.3}
\end{equation}
is valid under the condition on $u_x(0,t)$: there exists a constant 
$\beta_0 \in {\Bbb R}$
such that 
\begin{equation}
\displaystyle
\liminf_{\tau\longrightarrow\infty}\tau^{\beta_0}\Big\vert\int_0^Tu_x(0,t)e^{-\tau^2 t}dt
\Big\vert>0.
\label{1.4}
\end{equation}

Formula (\ref{1.3}) means that the exact location of 
the unknown boundary $\{a\}$ of the inside cavity $]a,\,\infty[$ 
can be detected by a single set of $u(0,t)$ and $u_x(0,t)$ for $t\in]0,\,T[$
provided $u_x(0,t)$ satisfies (1.4).
Note that there are other choices of $v$ to
define $I(\tau)$ which is useful for detecting the unknown boundary $\partial{D}
= \{a\}$ (for the detail, see \cite{I4}).

Our aim is to seek formulae which enable us to extract information about the unknown
boundary $\partial{D}$ for the three-dimensional case. To be our problem clear,
we rewrite (\ref{1.3}) by using another solution of the backward heat equation.

Given $y\in {\Bbb R}^1$ the function
$$\displaystyle
\tilde
E_{\tau}(x,y)=\frac{1}{\tau}e^{-\tau\vert x-y\vert}
$$
satisfies the equation
$$\displaystyle
\tilde{E}''(x)-\tau^2 \tilde{E}(x)+2\delta(x-y)=0
$$
in the whole line.  Note that 
$-2^{-1}\tilde{E}_{\tau}(x,y)$ is a fundamental solution of the operator 
$ \partial_x^2 - \tau^2$.

Let $p$ be an arbitrary fixed point in $]-\infty, 0[$.
Then 
$v(x,t)=e^{-\tau^2 t}\tilde{E}_{\tau}(x,p)$ also satisfies 
the backward heat 
equation $v_t+v_{xx}=0$ for $(x,t)\in\,]0,\,a[\times ]0,\,T[$.
Using this function, we define another indicator function
$$
\displaystyle
\tilde{I}(\tau,p)
=\int_0^{T}\left(-v_x(0,t)\,u(0,t)
+u_x(0,t)\,v(0,t)\right)dt.
$$
Since $v(x,t)=e^{\tau\,p}e^{-\tau^2t}e^{-\tau x}/\tau$ on $[0,\,a]$, we have
$\tilde{I}(\tau,p)=e^{\tau\,p}I(\tau)/\tau$.  From this and (\ref{1.3}) we 
obtain another formula
\begin{equation}
\displaystyle
\lim_{\tau\longrightarrow\infty}
\frac{1}{\tau}\log\vert \tilde{I}(\tau,p)\vert
=p-2a.
\label{1.5}
\end{equation}
The point is the interpretation of this right-hand side of (\ref{1.5}).
Since $p<0$, one can write $p-2a=-(\vert p\vert+2a)$.  
Hence we can see that
the quantity $\vert p\vert+2a$ in formula (\ref{1.5}) coincides with
the length of the broken path that starts at $x=p$, goes to $a$ (the point of
the boundary $\{a\}$ of the cavity) and returns to $x=0$ (the point of
the boundary of medium).

In this paper we establish a three-dimensional analogue of formula
(\ref{1.5}) (which is equivalent to (\ref{1.3}) as mentioned above).

\subsection{Description of the main result}

First we introduce a three dimensional analogue of $\tilde{I}(\tau,p)$.

\begin{Def}\label{Def of the indicator function}
Let $p$ be an {\it arbitrary} point outside $\Omega$.
Define the {\it indicator function} 
for the solution 
$u_f(x,t)$ of (\ref{1.1}) with a fixed $f \in L^2(\partial\Omega{\times}]0,\,T[)$ 
by the formula
$$\displaystyle
I(\tau, p)
=\int_{\partial\Omega}
\int_0^T\left(\frac{\partial v}{\partial\nu}(y,t)u_f(y,t)
-f(y,t)v(y,t)\right)dtdS_y,
$$
where
$$\displaystyle
v(x,t)=e^{-\tau^2t}E_{\tau}(x,p)
$$
and
$$\displaystyle
E_{\tau}(x,y)=\frac{e^{-\tau\vert x-y\vert}}
{2\pi\vert x-y\vert},\,x\not=y,\,\,
\tau > 0.
$$
\end{Def}
\par
Note that 
$E(x)=E_{\tau}(x,y)$ satisfies the equation
$\displaystyle
(\triangle-\tau^2)E(x)+2\delta(x-y)=0$
in ${\Bbb R}^3$ in the sense of distribution.
Thus if $y\in {\Bbb R}^3\setminus\overline\Omega$, then 
$E(x)=E_{\tau}(x,y)$
satisfies the equation
\begin{equation}
\displaystyle
(\triangle-\tau^2)E(x)=0\,\,\,\text{in}\,\Omega.
\label{1.6}
\end{equation}
Hence, the indicator function $I(\tau, p)$ in definition
\ref{Def of the indicator function} is suited to treat
three dimensional analogue of formula (\ref{1.5}).

Throughout this paper, we always assume
that the heat flux $f(y, t)$ belongs to
the space $L^2(\partial\Omega{\times}]0,\,T[)$.
Since the weak solution $u_f$ of (\ref{1.1}) uniquely
exists, 
the indicator function $I(\tau, p)$ is well-defined.
Our purpose in this paper is to clarify what information
can be obtained from this indicator function. To describe them, 
we need to introduce the following notations:

\begin{Def}\label{Minimum of the broken path}
Given $p$ outside $\Omega$ define
$$\displaystyle
l(p, D)=\inf_{(x,y)\in\partial D\times \partial\Omega}l_p(x,y),
$$
where
$$\displaystyle
l_p(x,y)=\vert p-x\vert+\vert x-y\vert,\,\,(x,y)\in\, {\Bbb R}^3\times{\Bbb R}^3.
$$
\end{Def}
The quantity $l(p,D)$ can be interpreted as the minimum length of 
broken paths that start at $p$, go to a point on $\partial D$ and 
return to a point on $\partial\Omega$.

\vskip 0.5em
We also introduce some sets of pair of points on $\partial D$ and 
$\partial\Omega$ related to $l(p,D)$.

\begin{Def}\label{Sorts of points attaining the minimum}
Given $z$ outside $D$ define
$$\begin{array}{c}
\displaystyle
{\cal G}(z)=\{x\in\partial D\,\vert\,\nu_x\cdot(z-x)=0\},\\
\\
\displaystyle
{\cal G}^{\pm}(z)=\{x\in\partial D\,\vert\,\pm\nu_x\cdot(z-x)>0\}.
\end{array}
$$
Let $p$ be an arbitrary point outside $\Omega$.
Define
$$\begin{array}{c}
\displaystyle
{\cal M}(p)=\{(x,y)\in\partial D\times\partial\Omega\,\vert\,
l(p, D)=l_p(x,y)\},\\
\\
\displaystyle
{\cal M}_1(p)=\{(x,y)\in {\cal M}(p)\,\vert\,x\in {\cal G}^+(p)\cap {\cal G}^+(y)\},\\
\\
\displaystyle
{\cal M}_2^{\pm}(p)=\{(x,y)\in {\cal M}(p)\,\vert\,x\in 
{\cal G}^{\pm}(p)\cap {\cal G}^{\mp}(y)\},\,\\
\\
\displaystyle
{\cal M}_g(p)=\{(x,y)\in {\cal M}(p)\,\vert\,x\in {\cal G}(p)\}.
\end{array}
$$
\end{Def}

\vskip 0.5em
Now we state what the indicator function $I(\tau, p)$ gives.
We put
\begin{equation}
g(y, \tau) = \int_{0}^{T}e^{-\tau^2t}f(y, t)dt
\quad(y \in \partial\Omega, \tau > 0). 
\label{1.11}
\end{equation}

\begin{Thm}\label{Theorem 1.1}
Assume that $f \in L^2(\partial\Omega{\times}]0,\,T[)$ and
there exists a constant $\mu \in {\bf R}$ 
such that the function $g(y, \tau)$ defined by (\ref{1.11}) belongs to
$C^{0, \alpha_0}(\partial\Omega)$ for
all large $\tau > 0$ and satisfies
\begin{equation}
\begin{array}{c}
\displaystyle
0 < \inf_{y \in \partial\Omega}
\liminf_{\tau \longrightarrow \infty}
\tau^{\mu}{\rm Re}\, g(y, \tau)
\leq \limsup_{\tau \longrightarrow \infty}
\tau^{\mu}
\Vert g(\cdot, \tau) \Vert_{C^{0, \alpha_0}(\partial\Omega)} 
< \infty.
\end{array}
\label{1.12}
\end{equation}
Then, the formula 
\begin{equation}
\displaystyle
\lim_{\tau \longrightarrow \infty}
\frac{1}{\tau}\log\vert I(\tau, p) \vert = -l(p, D),
\label{1.13}
\end{equation}
holds if $\partial{D}$ and $\partial\Omega$ satisfy the following four 
conditions:
\par\noindent\hskip12pt
(I.1) $\partial{D}$ is strictly convex,
\hskip36pt
(I.2) ${\cal M}_g(p) = \emptyset$,
\hskip36pt
(I.3) ${\cal M}_2^{-}(p) = \emptyset$,
\par\noindent\hskip12pt
(I.4) every point $(x_0, y_0) \in \partial{D}\times\partial\Omega$ attaining
$l(p, D)$ is non-degenerate critical point 
of $l_p(x, y)$.
\end{Thm}

\begin{Remark}\label{there are many heat flux!!}
\par
There exist many $f \in L^2(\partial\Omega{\times}]0,\,T[)$
satisfying (\ref{1.12}). For example, (\ref{1.12}) with $\mu = 2$ holds for 
$f \in C^1([0, T]; C^{0, \alpha_0}(\partial\Omega))$ with
$
\inf_{y \in \partial\Omega}f(0, y) > 0.
$
Indeed, integration by parts implies that
$$
\V{\tau^2g(\cdot, \tau) - f(0, \cdot)}_{C^{0, \beta}(\partial\Omega)}
\leq \tau^{-2}\max_{0 \leq t \leq T}
\V{\partial_tf(t, \cdot)}_{C^{0, \beta}(\partial\Omega)}
\quad(0 \leq \beta \leq \alpha_0).
$$
\end{Remark}

%
%
%

Formula (\ref{1.13}) for the three-dimensional problem (\ref{1.1}) can be 
interpreted as the analogous formula of (\ref{1.5}) for the one-dimensional 
case (\ref{1.2}).  Note that in the one-dimensional case, 
$\Omega =\, ]0, \infty[$ and $\partial{D} = \{a\}$.
Hence the length $l(p, D)$ for a point $p \notin \, ]0, \infty[ (=\Omega)$
is just $2a-p$ as appeared in (\ref{1.5}). Thus, from formula
(\ref{1.5}) we can find the unknown boundary $\partial{D} = \{a\}$.

In section \ref{Proof of Theorem 1.1.}, we prove 
theorem \ref{Theorem 1.1}. 
We briefly introduce the decomposition of $I(\tau, p)$ into 
the main part
$I_0(\tau, p)$ and remainder term. This decomposition
enables us to reduce the problem to the study of the asymptotic behaviour of 
$I_0(\tau, p)$, which is stated as theorem \ref{Theorem 2.1}. 
Sections \ref{Preliminaries} 
to \ref{Asymptotic behaviour of F_j(x, p, lambda)} 
are devoted to the proof of theorem \ref{Theorem 2.1}. In the last
part of section \ref{Proof of Theorem 1.1.}, 
we explain the necessity of the succeeding sections
for the proof of theorem \ref{Theorem 2.1}.

\subsection{Other previous results using the enclosure method}

To obtain other information about $D$ one may think about replacing $v$ in $I(\tau,p)$ with other special solutions of
the backward heat equation $(\partial_t+\triangle)v=0$ in $\Omega$.

In three-space dimensional case, define the indicator function $J_v(\tau)$ by
\begin{equation}
J_v = 
\int_{\partial\Omega}
\int_0^T\left(\overline{
\frac{\partial v}{\partial\nu}(y,t)}u_f(y,t)
-f(y,t)\overline{v(y,t)}\right)dtdS,
\label{indicator1}
\end{equation}
where $u_f$ is the solution of (\ref{1.1}), $v(x, t)$ is a solution of the backward heat equation
$(\partial_t+\triangle)v = 0$ in $[0, T]\times\Omega$ having the form
$v=e^{-\tau^2 t}q(x,\tau)$ and thus $(\triangle-\tau^2)q=0$ in $\Omega$.


Note that there are several possibilities of the choice of $v$ and $f$ in (\ref{indicator1}).

\noindent
Case $(\infty)$: This is an ideal case. It is assumed that one can obtain
$u_f$ on $\partial\Omega\times\,]0, T[$ corresponding to infinitely many $f$. 
In this case, we can design input heat flux $f$ to obtain information of 
$D$.  
In what follows, for integer $k$, 
we denote by $H^k(\Omega)$ the $L^2-$Sobolev space
defined by $H^2(\Omega) = \{ u \in L^2(\Omega) \vert 
\partial_x^{\alpha}u \in L^2(\Omega) \text{ for } \vert \alpha \vert \leq 2 \}$,
where the derivative $\partial_x^{\alpha}u$ is in distribution sense.
For an appropriate $\varphi \in L^2(0, T)$ and a function
$q(x, \tau)$ satisfying $(\triangle - \tau^2)q = 0$ in $\Omega$ with
$\Vert q(\cdot, \tau)\Vert_{H^2(\Omega)} = O(e^{C\tau})$
$(\tau \longrightarrow \infty)$ for some fixed constant $C > 0$,
we input heat flux $f(x, t; \tau)$ depending on $\tau \geq 1$
as
$$
f(x, t; \tau) = \varphi(t)\frac{\partial q}{\partial\nu}(x, \tau)
\qquad\text{on}\quad\partial\Omega\times]0, T[.
$$
For each $\tau \geq 1$, we put 
$v(x, t; \tau) = e^{-\tau^2t}q(x, \tau)$. 
Since $f \in L^2(\partial\Omega\times]0, T[)$, from the definition
of the weak solutions for (\ref{1.1}) and
$v \in C^1([0, T]; H^2(\Omega))$, using (\ref{indicator1}),
we can define $I_q(\tau) = J_v$. As is in \cite{IF, IK00, IK0}, from
elliptic estimates, it follows that there exists a constant $C > 0$ such that
\begin{equation}
C^{-1}\Vert \nabla_xq(\cdot, \tau)\Vert_{L^2(D)}^2
\leq \vert I_q(\tau) \vert \leq C\{\Vert \nabla_xq(\cdot, \tau)\Vert_{L^2(D)}^2
+\tau^2\Vert q(\cdot, \tau)\Vert_{L^2(D)}^2\}
\quad(\tau \geq 1).
\label{estimate for elliptic type argument}
\end{equation}
From (\ref{estimate for elliptic type argument})
and the asymptotic behaviour of $q(x,\tau)$ on $D$ as $\tau\longrightarrow\infty$, 
one can extract several quantities such as
$h_D(\omega)=\sup_{x\in D}x\cdot\omega$, 
$d_{D}(p)=\inf_{x\in D}\vert x-p\vert$
and $R_D(y)=\sup_{x\in D}\vert x-y\vert$ when $q$ is chosen appropriately.
Note also that \cite{IF} covers the case where the background conductivity is {\it isotropic}, {\it inhomogeneous} and known.
It makes use of a complex geometrical optics solution constructed by using a Faddeev-type Green function
for the modified Helmholtz equation.

\noindent
Case (I): On the contrary to Case $(\infty)$, let us consider the case where
we can only use one
set of data $(f, u_f)$ on $\partial\Omega\times ]0, T[$
as the measurement. In this case, we can not design 
the indicator function like as Case $(\infty)$.  However, as is in \cite{IK00}, we can extract
${\rm dist}(\partial\Omega, D) =\inf_{y\in\partial\Omega,x\in D}\vert x-y\vert$
from $u_f$ on $\partial\Omega\times\,]0,\,T[$ for a fixed $f$.
More precisely, we introduce the function $g(y, \tau)$ defined by
(\ref{1.11}).
Taking a function $q(x, \tau)$ as the weak solution to
\begin{equation}
\left\{\hskip-12pt
\begin{array}{llll}
&(\triangle -\tau^2)q(x, \tau) = 0
\quad
&\mbox{in } \Omega, \\[2mm]
&\displaystyle\frac{\partial q}{\partial{\nu}}(x, \tau) =  g(x, \tau),  
 & \mbox{on } \partial\Omega, \\
\end{array}
\right.
\label{resolvent eq. for Omega}
\end{equation}
and putting $v(t, x; \tau) = e^{-\tau^2t}q(x, \tau)$ in 
(\ref{indicator1}), we define
$I_q(\tau) = J_v$ as the indicator function.
The point is: $v$ depends on $f$.
This idea comes from \cite{IWave} in which 
an inverse obstacle scattering problem in the time domain has been considered.
For this indicator function, estimate 
(\ref{estimate for elliptic type argument}) can be also 
shown similarly to Case $(\infty)$.
Hence, we can extract  
${\rm dist}(\partial\Omega, D)$ from the indicator function by studying the asymptotic behaviour
of $q(x, \tau)$ on $D$ as $\tau\longrightarrow\infty$.  Note that in the last step
we employ the potential theoretic construction of the solution of 
(\ref{resolvent eq. for Omega})(cf. \cite{IK00}).

In both cases, the limit
$$
\lim_{\tau\longrightarrow\infty}\frac{-1}{2\tau}\log\vert{I_q(\tau)}\vert (= d_0)
$$
gives various quantities related to $D$, as described above.

The results are listed as follows:

\par\noindent
\begin{tabular}{|c|c|c|c|}
\hline
 & & & \\[-4mm]
Case & $f(x, t)$ & $q(x, \tau)$ in $v(x, t; \tau)
= e^{-\tau^2}q(x, \tau)$& $d_0$ \\[1mm] 
 & & & \\[-4mm]
\hline
 & & & \\[-4mm]
 &  & $q = e^{\tau\omega\cdot{x}}$ with $\omega \in S^2$ 
& $h_D(\omega)$ \\[2mm] 
\cline{3-4}
 & & & \\[-4mm]
$(\infty)$  
&  $\displaystyle\varphi(t)\frac{\partial q}{\partial\nu_x}(x,\tau)$ &
$q = \displaystyle\frac{e^{-\tau\vert x-p\vert}}{\vert x-p\vert}$ 
with $p \in {\Bbb R}^3\setminus\overline{\Omega}$
& $d_{D}(p)$ \\[4mm] 
\cline{3-4}
 & & & \\[-4mm]
& & $q = 
\begin{cases}
\displaystyle
\frac{e^{\tau\vert x-y\vert}-e^{-\tau\vert x-y\vert}}{\vert x-y\vert}, 
&(x \neq y),
\\
2\tau, & (x = y).
\end{cases}
$
with $y \in {\Bbb R}^3$ & $R_{D}(y)$ \\[5mm]
& & & \\[-4mm]
\cline{3-4}
& & & \\[-4mm]
& & 
$q = e^{c\tau^2(\omega+i\lambda_\tau\omega^\bot)\cdot{x}}$ with $c\tau > 1$, & \\
& & 
$\omega\cdot\omega^\bot = 0$,
$\lambda_\tau = \sqrt{1-\frac{1}{c^2\tau^2}}$,
$\omega$, $\omega^{\perp}\in S^2$ & $h_D(\omega)$ \\
& & & \\[-4mm]
\hline
 & & & \\[-4mm]
(I) & A fixed $f$ & the solution to (\ref{resolvent eq. for Omega}) & 
${\rm dist}(\partial\Omega, D)$\\
&  & for $g(y, \tau)$ given by (\ref{1.11}) &\\
\hline
\end{tabular}
\par\noindent

Note that we can also apply the idea in Case (I) to one-space dimensional case (\ref{1.2})
and obtain ${\rm dist}(\partial\Omega, D) = a$.
However, this is different
from formula (\ref{1.3}) (and (\ref{1.5}))
since in this formula, $v(x, t)$ does not have any relation with
the heat flux $ f(0, t)$ !
Hence, for treating three-space dimensional analogue of formula
(\ref{1.5}) (or (\ref{1.3})), we need to choose
$v(x, t)$ in (\ref{indicator1}) being {\it independent} of $f(x, t)$.
%

In the following table, our result in this paper is described using $l(p,D)$.
However, there are places with question marks.
Those indicate that we do not know what kind of information about $D$ can
be extracted from the corresponding indicator function.
To fill the places with suitable quantities we need further investigation in future.

\par\noindent
\begin{tabular}{|c|c|c|c|}
\hline
 & & & \\[-4mm]
Case & $f(x, t)$ &$q(x, \tau)$ in $v(x, t; \tau)
= e^{-\tau^2}q(x, \tau)$& $d_0/2$ \\[1mm] 
 & & & \\[-4mm]
\hline
 & & & \\[-4mm]
 & & $q = e^{\tau\omega\cdot{x}}$ with $\omega \in S^2$  & ? \\[2mm] 
\cline{3-4}
 & & & \\[-4mm]
(I)  &  A fixed $f$
 
& $q = \displaystyle\frac{e^{-\tau\vert x-p\vert}}{\vert x-p\vert}$ 
with $p \in {\Bbb R}^3\setminus\overline{\Omega}$
& $l(p,D)$\\[4mm] 
\cline{3-4}
 & & & \\[-4mm]
 & & $q = 
\begin{cases}
\displaystyle
\frac{e^{\tau\vert x-y\vert}-e^{-\tau\vert x-y\vert}}{\vert x-y\vert}, 
&(x \neq y),
\\
2\tau, & (x = y)
\end{cases}
$ with $y\in{\Bbb R}^3$ & ? \\[5mm]
& & & \\[-4mm]
\cline{3-4}
& & & \\[-4mm]
 & &
$q = e^{c\tau^2(\omega+i\lambda_\tau\omega^\bot)\cdot{x}}$
& \\
& & 
with $\omega\cdot\omega^\bot = 0$,$\lambda_\tau = \sqrt{1-\frac{1}{c^2\tau^2}}$,
$c\tau>1$, $\omega, \omega^\bot \in S^2$
& ? \\
\hline
\end{tabular}
\par\noindent

%

Anyway, it seems that the result and proof of this paper
suggest us the difficulty of the reconstruction problem using a single set of data. 
It will be interesting to find a simpler proof of the result.

\subsection{What is a difference from one-space dimensional case?}

It may be suspicious that too many
assumption on $f$, 
$\partial\Omega$ and $\partial{D}$ 
appears in theorem \ref{Theorem 1.1}.
In this subsection, we will explain why those assumption is required for the proof of (\ref{1.13}).

\par
In one-space dimensional case, we have formula (\ref{1.5}) provided
the input heat flux at $t = 0$ on the boundary $\{0\}$ satisfies (\ref{1.4})
for some $\beta_0$.  This condition on the heat flux ensures the strength
of the input heat flux at $t=0$ from below implicitly.
In three-space dimensional case, assumption (\ref{1.12}) in theorem \ref{Theorem 1.1} 
corresponds to this condition. 
Moreover, theorem \ref{Theorem 2.1} in section 2 tells us that
we do not need to input the heat flux at $t = 0$ on the whole boundary 
$\partial\Omega$. If we know, in advance, the set of all points $y\in\partial\Omega$ 
such that there exists a point 
$x\in\partial D$ with $(x,y)\in {\cal M}_1(p)\cup {\cal M}_2^-(p)$,
then it suffices to input heat flux at $t = 0$ 
supplied only on such special points $y\in\partial\Omega$.
Thus (\ref{1.12}) can be replaced with weaker one if this is the case, however,
it is not practical to assume such a priori information.

In three-space dimensional case, there are several type of the points 
$(x_0, y_0) \in \partial D\times\partial{\Omega}$ that attain the minimum
length $l(p, D)$ (i.e. $(x_0, y_0) \in {\mathcal M}(p)$). 
One type consists of broken rays of geometrical optics passing through $y_0$, $x_0$ 
and $p$ in this order. The pairs of such points $(x_0, y_0)$ consist of the set 
${\mathcal M}_1(p)$.  Note that in a special case, there may exist a point 
$(x_0, y_0) \in {\mathcal M}_1(p)$ such that $y_0$ is contained in the line segment $px_0$. 
This case just corresponds to one-space dimensional case.

In three-space dimensional case, there may also exist points 
$(x_0, y_0) \in {\mathcal M}(p)$ such that $x_0$ is on the line segment $py_0$.
These points belong to one of the three types of disjoint sets 
${\mathcal M}_2^+(p)$, ${\mathcal M}_2^-(p)$ and ${\mathcal M}_g(p)$.
As it can be seen in the proof of theorem \ref{Theorem 1.1}, it is not easy to measure
the contribution of points in ${\mathcal M}_g(p)$ to the asymptotic behavior of 
$I(\lambda, p)$. 
We can also see that the contribution of points in ${\mathcal M}_2^-(p)$ 
to the asymptotic behavior of $I(\tau, p)$ 
may cancel the one of the points belonging to ${\mathcal M}_1(p)$
(cf. theorem \ref{Theorem 2.1}). 
In theorem \ref{Theorem 1.1}, to avoid these cancelations, we assume 
${\mathcal M}_g(p){\cup}{\mathcal M}_2^-(p) = \emptyset$
(i.e. (I.2) and (I.3)).

Thus, in three-space dimensional case, the structure of ${\mathcal M}(p)$
becomes complicated.
This is one of the different points from one-space dimensional case and makes 
the problem for three-space dimensional case harder. 
However, we can give a condition on $\partial\Omega$ that ensures
${\mathcal M}_g(p){\cup}{\mathcal M}_2^-(p) = \emptyset$
(cf. proposition \ref{Proposition 2.2.}). 
And also, in propositions \ref{Proposition 2.1.} and \ref{Proposition 4.2}
a condition to ensure that a point 
$(x_0,y_0)\in {\cal M}(p)\setminus {\cal M}_g(p)$
is a non-degenerate critical point of $l_p$ on $\partial D\times\partial\Omega$
(cf. propositions \ref{Proposition 2.1.} and \ref{Proposition 4.2}), is given.
Using these sufficient conditions,
we can give examples covered by theorem \ref{Theorem 1.1}. 

As the next step it would be interesting
to know what kind of information can be extracted
from $l(p, D)$ given at all or some $p \in {\bf R}^3\setminus\overline{\Omega}$.
To our best knowledge, the complete answer to the question is unknown.
However, in section \ref{Upper bound of the location of D} we show that $l(p,D)$ yields some information
about an upper bound of the location of $D$.

In theorem \ref{Theorem 1.1}, we also assume that $\partial{D}$ is strictly convex. 
It seems that this assumption is too strong for the applications to practical inverse problems.
However, at the present time, technically, to treat the case of \lq\lq one measurement", 
we need such kind of a priori information on the unknown object $\partial{D}$.
We can also show a similar result to the case that
$D$ consists of several disjoint strictly convex domains. 
However, to treat this
case, we need to repeat the argument
which was used in the proof of theorem \ref{Theorem 1.1}.
Hence to keep this paper in an appropriate length, we restrict ourselves within introducing
theorem \ref{Theorem 1.1}. 

%
%

\subsection{A remark on the solution class}

Before closing this section, following \cite{DL}, we describe the class of solutions of the
initial boundary value problem for the
heat equation (\ref{1.1}).

For $T > 0$ and a Hilbert space $H$, 
$L^2(0, T; H)$ denotes the space of $H$-valued $L^2$
functions in $t \in [0, T]$. 
For two Hilbert spaces $H$ and $V$
with $V \subset H \subset V'$, we also introduce the space
$W(0, T; V, V') = \{\, u \,\vert\, u \in L^2(0, T; V),
u' \in L^2(0, T; V')\,\}$, where $V'$ is the dual space of the 
Hilbert space $V$, and $u'$ means
the (weak) derivative in $t \in [0, T]$.

As is in \cite{DL}, for any $f \in L^2(0, T; H^{-1/2}(\partial\Omega))$,
we call $u \in W(0, T; H^1(\Omega\setminus\overline{D}), 
(H^1(\Omega\setminus\overline{D}))')$ satisfies
\begin{equation}
\left\{
\begin{array}{l}
\displaystyle u_t-\triangle u=0
\,\,\text{in}\,(\Omega\setminus\overline D)\times\,]0,\,T[,\\[2mm]
\displaystyle \frac{\partial u}{\partial\nu}+\rho u=0
\,\,\text{on}\,\partial D\times]0,\,T[,\\[3mm]
\displaystyle
\frac{\partial u}{\partial\nu}=f\,\,\text{on}\,\partial\Omega\times\,]0,\,T[
\end{array}
\right.
\label{1.1-1}
\end{equation}
in the weak sense if $u$ satisfies 
\begin{align*}
\displaystyle
\langle u'(t),\varphi\rangle_{H^1(\Omega\setminus\overline{D})}
&+\int_{\Omega\setminus\overline{D}}
\nabla u(x,t)\cdot\nabla\varphi(x)dx
\\&
- \langle{\rho}u(t),\,\varphi\vert_{\partial{D}}
\rangle_{H^{1/2}(\partial D)}
=\langle f(t),\,\varphi\vert_{\partial\Omega}
\rangle_{H^{1/2}(\partial\Omega)}
\,\,\text{in}\,]0,\,T[
\end{align*}
in the sense of distribution on $]0,\,T[$ for all 
$\varphi\in H^1(\Omega\setminus\overline{D})$ 
and a.e. $t\in]0,\,T[$. 
In the above, the bracket 
$ \langle \cdot\, , \cdot\,\rangle_{V}$ is the pairing between 
a Hilbert space $V$ and $V'$.

We see that every $u\in W(0,\,T;H^1(\Omega\setminus\overline{D}), 
(H^1(\Omega\setminus\overline{D}))')$
is almost everywhere equal to
a continuous function of $[0,\,T]$ in $L^2(\Omega\setminus\overline{D})$ 
(Theorem 1 on p.473 in \cite{DL}).  Further, we have the following inclusion:
$$\displaystyle
W(0,\,T;H^1(\Omega\setminus\overline{D}), 
(H^1(\Omega\setminus\overline{D}))')\hookrightarrow
C^0([0,\,T];L^2(\Omega\setminus\overline{D})),
$$
where the space $C^0([0,\,T];L^2(\Omega\setminus\overline{D}))$ 
is equipped with the norm of uniform convergence. 
Thus one can consider $u(t)$ ($0 \leq t \leq T$)
as elements of $L^2(\Omega\setminus\overline{D})$.
Then we see that for any given $f \in L^2(0, T; 
H^{-1/2}(\partial\Omega))$
and $u_0\in L^2(\Omega\setminus\overline{D})$, there
exists a unique 
$u \in W(0, T; H^1(\Omega\setminus\overline{D}), 
(H^1(\Omega\setminus\overline{D}))')$ 
satisfying (\ref{1.1-1}) in the weak
sense and the initial condition $u(0)=u_0$ 
(Theorems 1 and 2 on p.512 in \cite{DL}).
We denote by $u_f$ the weak solution
of (\ref{1.1-1}) with $u(0)=0$ and this is the meaning of the weak
solution of (\ref{1.1}).

%
\setcounter{equation}{0}
\section{Proof of theorem 1.1.}\label{Proof of Theorem 1.1.}
%

We begin with choosing the main term $I_0(\tau, p)$ of 
$I(\tau, p)$.
Define 
$$\displaystyle
w(x,\tau)=\int_0^T e^{-\tau^2 t}u_f(x,t)dt,\,\,x\in
\Omega\setminus\overline D,
\,\,\tau > 0.
$$
Since $f \in L^2(\partial\Omega{\times}]0,\,T[)$, 
$u \in W(0, T; H^1(\Omega\setminus\overline{D}), 
(H^1(\Omega\setminus\overline{D}))')$ is the weak solution 
of (\ref{1.1-1}). From these facts, we can see that 
$w(\cdot, \tau) \in H^1(\Omega\setminus\overline{D})$
is the unique solution of the following elliptic 
boundary value problem in the weak sense:
\begin{equation}
\left\{
\begin{array}{l}
\displaystyle
(\triangle-\tau^2)w = u(x,T)e^{-\tau^2T}\,\,\text{in}\,\Omega\setminus\overline D,
\\[2mm]
\displaystyle
\frac{\partial w}{\partial\nu}+\rho(x)w=0\,\,\text{on}\,\partial D,
\quad
\displaystyle
\frac{\partial w}{\partial\nu}=g(y,\tau)\,\,\text{on}\,\partial\Omega,
\end{array}
\right.
\label{2.1-1}
\end{equation}
where $g(y, \tau)$ is the function defined by (\ref{1.11}). 
Using $w(x, \tau)$, we obtain the expression
\begin{equation*}
I(\tau, p)
=\int_{\partial\Omega}\left(\frac{\partial E_{\tau}(y,p)}
{\partial\nu}w(y,\tau)
-\frac{\partial w(y,\tau)}{\partial\nu}E_{\tau}(y,p)\right)dS_y.
\end{equation*}

Let us consider the solution
$w_0(x; \tau)$ of the following elliptic boundary value
problem:
\begin{equation}
\left\{
\begin{array}{l}
\displaystyle
(\triangle-\tau^2)w_0=0\,\,\text{in}\,\Omega\setminus\overline D,
\\[2mm]
\displaystyle
\frac{\partial w_0}{\partial\nu}+\rho(x)w_0=0\,\,\text{on}\,\partial D,
\quad
\displaystyle
\frac{\partial w_0}{\partial\nu}=g(y,\tau)\,\,\text{on}\,\partial\Omega.
\end{array}
\right.
\label{2.2}
\end{equation}
Note that $g(\cdot, \tau) \in L^2(\partial\Omega)$
for $f \in L^2(\partial\Omega{\times}]0,\,T[)$.
Hence usual elliptic 
theory implies that
for any $\tau > 0$, 
there exists the unique solution 
$w_0(\cdot, \tau) \in H^1(\Omega\setminus\overline{D})$ 
of (\ref{2.2}) in the weak sense. Thus, for $\tau > 0$,
we can introduce 
$$\displaystyle
I_0(\tau, p)
=\int_{\partial\Omega}\left(\frac{\partial E_{\tau}(y,p)}
{\partial\nu}w_0(y,\tau)
-\frac{\partial w_0(y,\tau)}{\partial\nu}E_{\tau}(y,p)
\right)dS_y.
$$

We can show that there exist constants $C > 0$ and
$\mu_0 > 0$ depending on $\partial{D}$, 
$f$ and $\rho$ such that 
$$
\vert I(\tau, p) - I_0(\tau, p) \vert 
\leq C\tau^{-1/2}e^{-\tau^2T}
\quad(\tau \geq \mu_0).
$$
In what follows, when the above estimate holds, we merely write
\begin{equation}
I(\tau, p) = I_0(\tau, p)+O(\tau^{-1/2}e^{-\tau^2T})
\qquad \tau \to \infty.
\label{2.3}
\end{equation}

This reduction is well known
(cf. section 2 in \cite{IK0}), however, 
for this paper to be self-contained, we
show it in Appendix C.  


Now we state the asymptotic behavior of $I_0(\lambda, p)$ being the essential
part of this paper. 
\vskip0.5em

\begin{Thm}\label{Theorem 2.1}
Assume that $f \in L^2(\partial\Omega{\times}]0,\,T[)$, and
$\partial D$ and $\partial\Omega$ satisfy (I.1), (I.2) and (I.4) in
theorem \ref{Theorem 1.1}.
Then the set ${\cal M}(p)$ is finite.  Moreover, we have
\begin{equation}
\displaystyle
I_0(\tau,p)=\frac{1}{\tau}e^{-\tau l(p,D)}\left\{A(\tau,p)g
+\Vert g(\,\cdot\,,\tau)\Vert_{C^{0,\alpha_0}(\partial\Omega)}
O(\tau^{-\alpha_0/2})\right\}
\label{2.5}
\end{equation}
as $\tau\longrightarrow\infty$,
where
\begin{equation}
\begin{array}{c}
\displaystyle
A(\tau,p)g=\sum_{(x_0,y_0)\in {\cal M}_1(p)}
C(x_0,y_0)H^+(x_0,y_0,p)g(y_0,\tau)\\
\\
\displaystyle
+\sum_{(x_0,y_0)\in {\cal M}^{-}_2(p)}
C(x_0,y_0)H^{-}(x_0,y_0,p)g(y_0,\tau).
\end{array}
\label{2.6}
\end{equation}
In (\ref{2.6}), the coefficients 
$C(x_0,y_0)$ for each $(x_0,y_0)\in 
{\cal M}_1(p)\cup {\cal M}_2^{-}(p)$ is a positive constant
independent of $g$ (i.e $f$) 
and
$$
\displaystyle
H^{\pm}(x,y,p)
=\frac{1}{\vert x-p\vert\vert x-y\vert}
\nu_x\cdot\left\{\frac{p-x}{\vert p-x\vert}
\pm\frac{y-x}{\vert y-x\vert}\right\},\,(x,y)\in\partial D\times\partial\Omega.
$$

\end{Thm}

\begin{Remark}\label{Remark 2.1}
We have $H^+(x_0,y_0,p)>0$ for $(x_0,y_0)\in {\cal M}_1(p)$ and 
$H^-(x_0,y_0,p)<0$ for $(x_0,y_0)\in {\cal M}_2^-(p)$
(cf. (2) of proposition \ref{Proposition 4.1} and (\ref{4.8}) 
in subsection \ref{The structure of {cal M}(p)}).
\end{Remark}

\vskip 1em


Note that all the points of the set ${\cal M}(p)$ are critical points
of $l_p$ on $\partial D\times\partial\Omega$ and for each $(x,y)\in {\cal M}(p)$
the Hessian at $(x,y)$ of any local representation of $l_p$
in a neighbourhood of $(x,y)$ has no negative eigenvalues.  Thus
a point $(x,y)\in {\cal M}(p)$ is a non-degenerate critical point of $l_p$ on 
$\partial D\times\partial\Omega$ if and only if the Hessian at $(x,y)$ of a 
local representation of $l_p$ in a neighbourhood of $(x,y)$ is positive 
definite.  
Thus the conclusion of the finiteness of ${\cal M}(p)$ in 
theorem \ref{Theorem 2.1} is trivial.

\vskip 1em

Using theorem \ref{Theorem 2.1}, we can obtain theorem \ref{Theorem 1.1}. 
Here, we continue the proof of theorem \ref{Theorem 1.1} 
assuming theorem \ref{Theorem 2.1} holds. 

\vskip1em
{\it\noindent Proof of theorem \ref{Theorem 1.1}.}
Since we consider the case ${\cal M}_g(p)\cup {\cal M}_2^{-}(p) = \emptyset$,  
from (\ref{2.6}), (\ref{1.12}) and remark 2.1, it follows that 
there exist constant $C > 0$ and
$\mu_0 > 0$ such that
$$
C^{-1} \leq \tau^{\mu}{\rm Re}\, A(\tau, p)g \leq
\tau^{\mu}\vert A(\tau, p)g \vert \leq C
\quad(\tau \geq \mu_0).
$$
Combining this estimate with (\ref{2.5}) and 
(\ref{1.12}), we obtain 
\begin{equation}
C_1^{-1} \leq \tau^{\mu+1}\vert 
e^{\tau l(p, D)}I_0(\tau,p)\vert \leq C_1
\qquad(\tau \geq \mu_1)
\label{2.6.1}
\end{equation}
for some constants $C_1 > 0$ and $\mu_1 > 0$ 
independent of $\tau$. 
From the above estimate and (\ref{2.3}), it follows that 
$$
C_2^{-1} \leq \tau^{\mu+1}
e^{\tau l(p,D)}\vert I(\tau,p)\vert \leq C_2
\qquad(\tau \geq \mu_2)
$$
for some constants $C_2 > 0$ and $\mu_2 > 0$ independent 
of $\tau$. 
This estimate shows theorem \ref{Theorem 1.1} holds.

\noindent
$\Box$

\vskip1em

From the above proof of theorem \ref{Theorem 1.1}, we can 
see that formula (\ref{1.13}) in theorem \ref{Theorem 1.1}
is given by (\ref{2.6.1}).
Using $A(\tau, p)g$ in (\ref{2.6}), we can give sufficient
conditions for getting (\ref{2.6.1}) (i.e. (\ref{1.13})).

\begin{Cor}\label{Corollary 2.1}
Assume that there exists a positive number $\mu$ such that
the function $g(y, \tau)$ defined by (1.7) belongs to
$g(\cdot, \tau) \in C^{0, \alpha_0}(\partial\Omega)$ for
all $\tau  > 0$ large enough and satisfies
\begin{equation}\displaystyle
\liminf_{\tau\longrightarrow\infty}\tau^{\mu}\vert A(\tau,p)g\vert>0
\label{2.7}
\end{equation}
and
\begin{equation}\displaystyle
\lim_{\tau\longrightarrow\infty}\frac{\tau^{\mu}\Vert g(\,\cdot\,,\tau)\Vert_{C^{0,\alpha_0}(\partial\Omega)}}
{\tau^{\alpha_0/2}}=0.
\label{2.8}
\end{equation}
Then formula (\ref{1.13}), that is,
$$\displaystyle
\lim_{\tau\longrightarrow\infty}\frac{1}{\tau}\log\vert I(\tau,p)\vert=-l(p, D)
$$
is valid.
\end{Cor}

Note that in theorem \ref{Theorem 1.1}, we assume (\ref{1.12}) 
to ensure (\ref{2.7}) holds. 
However, (\ref{1.12}) is too strong.
We do not need to input the heat flux $f$ 
at $t = 0$ on the whole boundary $\partial\Omega$.
From the form (\ref{2.6}) of $A(\tau, p)g$, we can see that
it is enough to supply $f$ at $t = 0$ only on
the set of all points $y\in\partial\Omega$ such that there exists a point 
$x\in\partial D$ with $(x,y)\in {\cal M}_1(p)\cup {\cal M}_2^-(p)$. 
Note that condition (\ref{2.7}) also gives a lower bound estimate for the strength
of the input heat flux $f$ at $t=0$.
If both ${\cal M}_1(p)$ and ${\cal M}_2^-(p)$ are not empty,
a cancelation in $A(\tau,p)g$ may occur (see remark \ref{Remark 2.1})
and thus it is delicate whether (\ref{2.7}) holds or not.
Another condition (\ref{2.8}) is not a serious one. For example, 
if $f = 1$ on $\partial\Omega\times]0,\,T[$, 
then (\ref{2.8}) is satisfied with $\mu=2$. Note that in this case, (\ref{2.7}) 
also holds with $\mu=2$ if $A(\tau, p)g$ does not vanish.

It is crucial to represent 
the main term $I_0(\tau, p)$ 
by using Laplace type integrals (cf. 
proposition \ref{Proposition 3.1}) for the proof of theorem \ref{Theorem 2.1}.
This is done in subsection 
\ref{the decompostion of I_0(lambda,p)}. We construct
the solution $ w_0(x, \tau)$ of (\ref{2.2}) by single layer 
potentials on $\partial{D}$
and $\partial\Omega$ in potential theory. Using this expression, 
we decompose the main term into some parts. Each term can be 
reduced to a Laplace type integral
over $\partial\Omega\times\partial D$ with a large parameter 
$\tau$.

In each integral, the exponential terms are just given by
$e^{-{\tau}l_p(x, y)}$. Thus the points $(x_0, y_0) \in
\partial{D}\times\partial\Omega$ attaining the minimum $l(p, D)$
of $l_p(x, y)$, (i.e., $(x_0, y_0) \in {\cal M}(p)$) determine the
asymptotic behaviour of $I_0(\tau, p)$. 
In subsection \ref{The structure of {cal M}(p)} of
section \ref{the decompostion of I_0(lambda,p)}, we study 
the structure of the set ${\cal M}(p)$. 
\par

In section \ref{Proof of the main result}, we give a proof of
theorem \ref{Theorem 2.1} using the Laplace method.
Here, we need to have asymptotic behaviour of the amplitude
functions in the Laplace integrals. These key facts are described in lemma \ref{Lemma 5.1}.
In section \ref{Asymptotic behaviour of F_j(x, p, lambda)} the proof of lemma \ref{Lemma 5.1}
is given.

Since the amplitude
functions contain terms defined by using the inverse of an
integral operator on $\partial{D}$, the problem is
eventually reduced to obtaining some estimates of the kernel
$K_\tau^\infty(x, y)$ of an operator of the form 
$K_\tau(I - K_\tau)^{-1}$, where $K_\tau$ is an integral operator on
$\partial{D}$ with the kernel $K_\tau(x, y)$ estimated by
\begin{equation}
\vert {K_\tau(x, y)} \vert \leq C\tau\,{e}^{-\tau\vert
x - y \vert}\,\,(x, y \in \partial{D}, \tau > 0).
\label{2.10}
\end{equation}
We need to show that kernel $K_\tau^\infty(x, y)$ can be estimated by the same 
exponential term
${e}^{-\tau\vert x - y \vert}$ as in the estimate (\ref{2.10}).
Therefore we need more precise argument than that of usual 
classical potential theory although we study the kernels of 
the repeated integral operators
$K_\tau^n$ ($n = 1, 2, \ldots$) according to the classical 
approach. The needed estimates of the integral kernels
are given in \cite{IK2}.
Here, only the result used in this paper is summarized in 
subsection \ref{Basic estimates of integral kernels} 
(cf. theorem \ref{Theorem 3.1}). 
%
\par

The Laplace method requires the non-degenerateness of 
$l_p(x, y)$ at $(x, y) \in {\cal M}(p)$. In section \ref{Sufficient conditions},
sufficient conditions of non-degenerateness of $l_p(x, y)$ are given.
Using these conditions, we can give examples covered by 
theorems \ref{Theorem 1.1}, \ref{Theorem 2.1} and 
corollary \ref{Corollary 2.1}.
%

To make this paper self-contained
we add two appendixes A and B.
In Appendix A, we give a proof of one version of the Laplace 
method used to show the main result.
Appendix B is devoted to a computation of Weingarten map for ellipsoids, 
which is used to treat the examples in section \ref{Sufficient conditions}.

%
\setcounter{equation}{0}
\section{Preliminaries}\label{Preliminaries}
%

\subsection{the decompostion of $I_0(\tau,p)$}
\label{the decompostion of I_0(lambda,p)}

We employ the layer potential approach for 
the construction of $w_0$.

Given $g\in C(\partial\Omega)$ and $h\in C(\partial D)$ define
$$\begin{array}{c}
\displaystyle
V_{\Omega}(\tau)g(x)=\int_{\partial\Omega}E_{\tau}(x,y)g(y)dS_y,\,\,
x\in{\Bbb R}^3\setminus\partial\Omega,\\
\\
\displaystyle
V_D(\tau)h(x)=\int_{\partial D}E_{\tau}(x, z)h(z)dS_z,\,\,
x\in{\Bbb R}^3\setminus\partial D.
\end{array}
$$

We construct $w_0$ in the form
\begin{equation}
\displaystyle
w_0(x,\tau)=V_{\Omega}(\tau)\varphi(x,\tau)+V_{D}(\tau)\psi(x,\tau), 
\label{3.1}
\end{equation}
where $\varphi(\,\cdot\,,\tau)\in C(\partial\Omega)$ and 
$\psi(\,\cdot\,,\tau)\in C(\partial D)$ are
unknown functions to be determined.

Here we cite some well known facts for $V_{\Omega}(\tau)$ 
and $V_{D}(\tau)$ from potential theory
(cf. \cite{Ms}).

\noindent
$\bullet$  $V_{\Omega}(\tau)g$ satisfies 
$(\triangle-\tau^2)V_{\Omega}(\tau)g=0$ in 
${\Bbb R}^3\setminus\partial\Omega$.

\noindent
$\bullet$  $V_{D}(\tau)h$ satisfies 
$(\triangle-\tau^2)V_{D}(\tau)h=0$ in 
${\Bbb R}^3\setminus\partial D$.

\noindent
These yield that $w_0$ having the form (\ref{3.1}) satisfies 
the equation $(\triangle-\tau^2)w_0=0$ in 
$\Omega\setminus\overline D$.

In what follows, we denote by $B(X, Y)$ the space consisting of 
continuous linear operators from 
a normed space $X$ to a Fr\'echet space $Y$. Note that $B(X, Y)$ is the space 
consisting of all bounded linear operators when $X$ and $Y$ are Banach spaces. 
We also put $B(X) = B(X, X)$.

\noindent
$\bullet$  $V_{\Omega}(\tau)\in B(C(\partial\Omega), 
C^{\infty}({\Bbb R}^3\setminus\partial\Omega)\cap C({\Bbb R}^3))$
and the Neumann derivative for $V_{\Omega}(\tau)g$ at 
$x\in\partial\Omega$
$$\displaystyle
\frac{\partial}{\partial\nu_x}V_{\Omega}(\tau)g\vert_{\partial\Omega}(x)
=\lim_{\epsilon\downarrow 0}\sum_{j=1}^3(\nu_x)_j
\left(\frac{\partial}{\partial x_j}V_{\Omega}(\tau)g\right)(x-\epsilon\nu(x))
$$
exists and is given by the formula
$$\displaystyle
\frac{\partial}{\partial\nu_x}V_{\Omega}(\tau)g\vert_{\partial\Omega}(x)
=g(x)+S_{\partial\Omega}(\tau)g(x),
$$
where
$$
\displaystyle
S_{\partial\Omega}(\tau)g(x)=
\int_{\partial\Omega}
\frac{\partial}{\partial\nu_x}E_{\tau}(x,y)g(y)dS_y,
\,\,x\in\partial\Omega.
$$

\noindent
$\bullet$  $V_{D}(\tau) \in 
B(C(\partial D), C^{\infty}({\Bbb R}^3\setminus\partial D)
\cap C({\Bbb R}^3))$
and the Neumann derivative for $V_{D}(\tau)h$ at 
$x\in\partial D$
$$\displaystyle
\frac{\partial}{\partial\nu_x}V_{D}(\tau)h\vert_{\partial D}(x)
=\lim_{\epsilon\downarrow 0}\sum_{j=1}^3(\nu_x)_j
\left(\frac{\partial}{\partial x_j}V_{D}(\tau)h\right)
(x+\epsilon\nu(x))
$$
exists and is given by the formula
$$\displaystyle
\frac{\partial}{\partial\nu_x}V_{D}(\tau)h\vert_{\partial D}(x)
=-h(x)+S_{\partial D}(\tau)h(x),
$$
where
$$
\displaystyle
S_{\partial D}(\tau)h(x)=\int_{\partial D}
\frac{\partial}{\partial\nu_x}E_{\tau}(x,z)h(z)dS_z,
\,\,x\in\partial D.
$$

\noindent
$\bullet$  For $\tau > 0$, 
$S_{\partial\Omega}(\tau)\in B(C(\partial\Omega))$
and $S_{\partial D}(\tau)\in B(C(\partial D))$.  
Moreover there exists a positive constant
$C$ such that these operator norms are bounded by 
$C\tau^{-1}$.

Using these properties, we can show that 
$w_0$ having the form (\ref{3.1}) 
satisfies the boundary conditions in (\ref{2.2})
if and only if $\varphi$ and $\psi$ satisfies the system of integral 
equations on $\partial\Omega\cup\partial D$:
\begin{equation}\begin{array}{c}
\displaystyle
\varphi(x,\tau)+S_{\partial\Omega}(\tau)\varphi(x,\tau)
+X_{\partial\Omega}(\tau)\psi(x,\tau)=g(x,\tau)
\,\,\text{on}\,\partial\Omega, 
\\
\\
\displaystyle
\psi(x,\tau)-(X_{\partial D}(\tau)+\rho(x)V_{\Omega}(\tau))
\varphi(x,\tau)
-(S_{\partial D}(\tau)+\rho(x)V_D(\tau))\psi(x,\tau)=0
\,\,\text{on}\,\partial D,
\end{array}
\label{3.2}
\end{equation}
where
$$\begin{array}{c}
\displaystyle
X_{\partial\Omega}(\tau)\psi(x,\tau)
=\int_{\partial D}\frac{\partial}{\partial\nu_x}E_{\tau}(x,z)
\psi(z,\tau)dS_z\,\,\text{on}\,\partial\Omega,\\
\\
\displaystyle
X_{\partial D}(\tau)\varphi(x,\tau)
=\int_{\partial\Omega}\frac{\partial}{\partial\nu_x}E_{\tau}(x,y)
\varphi(y,\tau)dS_y\,\,\text{on}\,\partial D.
\end{array}
$$
For the concise expression of $\varphi$ and $\psi$ we 
introduce the $2\times 2$ matrix operator acting on 
$C(\partial\Omega)\times C(\partial D)$
\begin{equation*}
\displaystyle
Y(\tau)=(Y_{ij}(\tau))
=\left(\begin{array}{cc}
\displaystyle -S_{\partial\Omega}(\tau) & 
\displaystyle -X_{\partial\Omega}(\tau)\\
\\
\displaystyle X_{\partial D}(\tau)+\rho(x)V_{\Omega}(\tau)
&
\displaystyle S_{\partial D}(\tau)+\rho(x)V_D(\tau)
\end{array}
\right).
\end{equation*}
Using $Y(\tau)$, we can write the equations (\ref{3.2}) as 
$$\displaystyle
(I-Y(\tau))\left(\begin{array}{c} \varphi\\
\\
\psi\end{array}\right)
=\left(\begin{array}{c} g\\
\\
0\end{array}
\right).
$$
Using a similar argument for the proof of the boundedness for 
$S_{\partial\Omega}(\tau)$ and $S_{\partial D}(\tau)$,
we know that: if $\tau>0$, then 
$X_{\partial\Omega}(\tau)\in B(C(\partial D), C(\partial\Omega))$, 
$X_{\partial D}(\tau)\in B(C(\partial\Omega), C(\partial D))$, and 
there exists a positive constant $C$ such that
these operator norms are bounded by $C/\tau$.
For $V_{\Omega}(\tau)$ and $V_{D}(\tau)$, we can show that
$V_{\Omega}(\tau)\in B(C(\partial\Omega), C(\partial D))$, 
$V_D(\tau)\in B(C(\partial D))$
and they have similar estimates.

Therefore we conclude that there exists a positive constant $C$
such that, for all $\tau > 0$
$$\displaystyle
\Vert Y(\tau)\Vert_{B(C(\partial\Omega)\times C(\partial D))}
\le C\tau^{-1}.
$$
This ensures that if $\tau$ is large enough, then 
the Neumann series
$\sum_{n=0}^{\infty}Y(\tau)^n$ is absolutely convergent with the operator norm and
coincides with $(I-Y(\tau))^{-1}$.
$\varphi$ and $\psi$ are given by
\begin{equation}
\displaystyle
\left(\begin{array}{c}
\displaystyle
\varphi
\\
\\
\displaystyle
\psi
\end{array}
\right)
=(I-Y(\tau))^{-1}\left(\begin{array}{c} g\\
\\
0
\end{array}
\right).
\label{3.4}
\end{equation}
This completes the construction of $w_0$.

Next, we write $I_0(\tau,p)$ in terms of only 
$\varphi$ given by (\ref{3.4}).
For the definition of $I_0(\tau,p)$, it follows that
$$\displaystyle
I_0(\tau,p)
=\int_{\partial D}\left(\frac{\partial E_{\tau}}{\partial\nu}
+\rho E_{\tau}\right)(y, p)w_0(y,\tau)dS_y.
$$
Indeed, integration by parts implies that
$$
\begin{array}{l}
\displaystyle
I_0(\tau,p)= 
\int_{\partial\Omega}\left(\frac{\partial E_{\tau}}
{\partial\nu}w_0-\frac{\partial w_0}{\partial\nu}E_{\tau}
\right)dS_y
\\[3mm]
\displaystyle
= 
\int_{\Omega\setminus\overline{D}}\left(\big((\triangle-\tau^2) E_{\tau}\big)w_0
- \big(({\triangle}-\tau^2)w_0\big)E_{\tau}\right)dx
\displaystyle
+
\int_{\partial{D}}\left(\frac{\partial E_{\tau}}
{\partial\nu}w_0-\frac{\partial w_0}{\partial\nu}E_{\tau}
\right)dS_y
\\[3mm]
\displaystyle
= \int_{\partial{D}}\left(\frac{\partial E_{\tau}}
{\partial\nu}w_0+{\rho}w_0E_{\tau}
\right)dS_y
= \int_{\partial{D}}\left(\frac{\partial E_{\tau}}
{\partial\nu}+{\rho}E_{\tau}\right)w_0dS_y.
\end{array}
$$
Using the above equality and (\ref{3.1}), one has the decomposition
\begin{align}
I_0(\tau,p)&=J_1(\tau,p)+J_2(\tau,p)
\nonumber
\\
&\equiv
\int_{\partial D}\left(\frac{\partial}{\partial\nu}+\rho\right)
E_{\tau}(x,p)V_{\Omega}(\tau)\varphi(x,\tau)dS_x
\label{3.5}
\\&\hskip10mm
+
\int_{\partial D}\left(\frac{\partial}{\partial\nu}+\rho\right)
E_{\tau}(x,p)V_{D}(\tau)\psi(x,\tau)dS_x.
\nonumber
\end{align}
A direct computation gives
\begin{equation}\displaystyle
\left(\frac{\partial}{\partial\nu_x}+\rho(x)\right)E_{\tau}(x,y)
=\frac{1}{2\pi}e^{-\tau\vert x-y\vert}H(x,y,\tau)
\quad(x \in \partial{D}\cup\partial\Omega, y \in \R^3, x \neq y),
\label{3.6}
\end{equation}
where
$$\displaystyle
H(x,y,\tau)
=\frac{\nu_x\cdot(y-x)}{\vert x-y\vert}
\left(\frac{\tau}{\vert x-y\vert}+\frac{1}{\vert x-y\vert^2}\right)
+\frac{\rho(x)}{\vert x-y\vert}.
$$
This yields
\begin{equation}\displaystyle
J_1(\tau,p)
=\left(\frac{1}{2\pi}\right)^2
\int_{\partial\Omega}dS_y\varphi(y,\tau)
\int_{\partial D}\frac{H(x,p,\tau)}{\vert x-y\vert}
e^{-\tau\,l_p(x,y)}dS_x.
\label{3.7}
\end{equation}

Set $w_2(x,\tau)=V_D(\tau)\psi(x,\tau)$ and write
$$\displaystyle
J_2(\tau,p)=\int_{\partial D}\frac{\partial}{\partial\nu}E_{\tau}(x,p)w_2(x,\tau)dS_x
+\int_{\partial D}\rho E_{\tau}(x,p)w_2(x,\tau)dS_x.
$$
Note that $w_2 \in C^\infty({\Bbb R}^3\setminus\overline{D})$ 
satisfies the equation $(\triangle -\tau^2)w_2=0$ 
in ${\Bbb R}^3\setminus\overline D$. 
For sufficiently large $R > 0$, this function 
belongs to $H^2$ for $\v{x} > R$ and
$\displaystyle\lim_{h \to 0}\frac{{\partial}w_2}{\partial\nu}
(x+h\nu_x)$ exists in $C(\partial\Omega)$.
Since $E_\tau(x, p)$ satisfies $\displaystyle
(\triangle_x-\tau^2)E_\tau(x, p)+2\delta(x-p)=0$, 
integration by parts and the property of $w_2$ mentioned above yield
$$
\displaystyle
\int_{\partial D}\frac{\partial}{\partial\nu}E_{\tau}(x,p)w_2(x,\tau)dS_x
=2w_2(p,\tau)+\int_{\partial D}E_{\tau}(x,p)
\frac{\partial w_2}{\partial\nu}(x,\tau)dS_x.
$$
From the property of $V_D(\tau)$ and 
the second equation in (\ref{3.2})
we obtain
$$
\displaystyle
\left(\frac{\partial}{\partial\nu}+\rho\right)w_2(x,\tau)=-Y_{21}(\tau)\varphi(x,\tau)\,\,\text{on}\,\partial D.
$$
Therefore we have
\begin{equation}\displaystyle
J_2(\tau,p)
=2w_2(p,\tau)-\int_{\partial D}E_{\tau}(x,p)Y_{21}(\tau)\varphi(x,\tau)dS_x.
\label{3.8}
\end{equation}

From (\ref{3.6}) we know that
\begin{equation}
Y_{21}(\tau)\varphi(x,\tau)
=\frac{1}{2\pi}\int_{\partial\Omega}e^{-\tau\vert x-y\vert}H(x,y,\tau)\varphi(y,\tau)dS_y.
\label{3.9}
\end{equation}
This yields
\begin{equation}
\begin{array}{c}
\displaystyle
\int_{\partial D}E_{\tau}(x,p)Y_{21}(\tau)\varphi(x,\tau)dS_x\\
\\
\displaystyle
=\left(\frac{1}{2\pi}\right)^2
\int_{\partial\Omega}dS_y\varphi(y,\tau)
\int_{\partial D}\frac{H(x,y,\tau)}{\vert x-p\vert}
e^{-\tau\,l_p(x,y)}dS_x.
\end{array}
\label{3.10}
\end{equation}

Note also that
$$\displaystyle
\psi(x,\tau)=(I-Y_{22}(\tau))^{-1}Y_{21}(\tau)\varphi(x,\tau),
\,\,\tau>>1.
$$
In what follows we denote by ${}^tY_{22}(\tau)$ the formal adjoint operator
defined by
$$
\int_{\partial{D}}({}^tY_{22}(\tau)f)(x)h(x)dx 
= \int_{\partial{D}}f(x)(Y_{22}(\tau)h)(x)dx
\qquad(f, h \in C(\partial{D})).
$$
From the definition of ${}^tY_{22}(\tau)$, it follows that
$ {}^t((I-Y_{22}(\tau))^{-1}) = (I-{}^tY_{22}(\tau))^{-1}$.
From these facts, it holds that
\begin{equation}
\begin{array}{c}
\displaystyle
w_2(p,\tau)=V_D(\tau)\psi(p,\tau)\\
\\
\displaystyle
=\int_{\partial D}E_{\tau}(p,x)(I-Y_{22}(\tau))^{-1}
Y_{21}(\tau)\varphi(x,\tau)dS_x\\
\\
\displaystyle
=\frac{1}{2\pi}
\int_{\partial D}\frac{e^{-\tau\vert x-p\vert}}{\vert x-p\vert}
(I-Y_{22}(\tau))^{-1}Y_{21}(\tau)\varphi(x,\tau)dS_x\\
\\
\displaystyle
=\frac{1}{2\pi}
\int_{\partial D}Y_{21}(\tau)\varphi(x,\tau)
\left((I-{}^tY_{22}(\tau))^{-1}
\frac{e^{-\tau\vert\,\cdot\,-p\vert}}{\vert\,\cdot\,-p\vert}\right)dS_x\\
\\
\displaystyle
=\frac{1}{2\pi}
\int_{\partial D}e^{-\tau\vert x-p\vert}Y_{21}(\tau)\varphi(x,\tau)
\cdot e^{\tau\vert x-p\vert}
\left((I-{}^tY_{22}(\tau))^{-1}
\frac{e^{-\tau\vert\,\cdot\,-p\vert}}{\vert\,\cdot\,-p\vert}\right)dS_x.
\end{array}
\label{3.11}
\end{equation}
Define
\begin{equation}
F(x,p,\tau)
=e^{\tau\vert x-p\vert}
\left((I-{}^tY_{22}(\tau))^{-1}\frac{e^{-\tau\vert\,\cdot\,-p\vert}}{\vert\,\cdot\,-p\vert}\right)(x).
\label{3.12}
\end{equation}
A combination of (\ref{3.9}) and (\ref{3.11}) gives
\begin{equation}
\displaystyle
w_2(p,\tau)
=\left(\frac{1}{2\pi}\right)^2
\int_{\partial\Omega}dS_y\varphi(y,\tau)\int_{\partial D}
e^{-\tau\,l_p(x,y)}H(x,y,\tau)F(x,p,\tau)dS_x.
\label{3.13}
\end{equation}

Finally from (\ref{3.5}), (\ref{3.7}), (\ref{3.8}), (\ref{3.10}), 
(\ref{3.13}), we obtain the representation formula of $I_0(\tau,p)$:
\begin{equation}
\begin{array}{c}
\displaystyle
(2\pi)^2I_0(\tau,p)=\int_{\partial\Omega}dS_y\varphi(y,\tau)\\
\\
\displaystyle
\times
\int_{\partial D}
e^{-\tau\,l_p(x,y)}
\left\{\frac{H(x,p,\tau)}{\vert x-y\vert}
-\frac{H(x,y,\tau)}{\vert x-p\vert}
+2H(x,y,\tau)F(x,p,\tau)\right\}dS_x.
\end{array}
\label{3.14}
\end{equation}

\begin{Lemma}\label{Lemma 3.1}
$$\displaystyle
{}^tY_{22}(\tau)h(z)
=\frac{1}{2\pi}
\int_{\partial D}e^{-\tau\vert x-z\vert}H(x,z,\tau)h(x)dS_x,\,\,h\in C(\partial D), z\in\partial D.
$$

\end{Lemma}

{\it\noindent Proof.}
Let $f, h\in C(\partial D)$.
Since $Y_{22}(\tau)=S_{\partial D}(\tau)+\rho(x)V_D(\tau)$, we have
$$\begin{array}{c}
\displaystyle
\int_{\partial D}{}^tY_{22}(\tau)h(z)\cdot f(z)dS_z
=\int_{\partial D}h(x)\cdot Y_{22}(\tau)f(x)dS_x\\
\\
\displaystyle
=\int_{\partial D}dS_x\,h(x)\int_{\partial D}
\left\{\frac{\partial}{\partial\nu_x}E_{\tau}(x,z)
+\rho(x)E_{\tau}(x,z)\right\}f(z)dS_z\\
\\
\displaystyle
=\int_{\partial D}dS_z\,f(z)\int_{\partial D}
\left\{\frac{\partial}{\partial\nu_x}E_{\tau}(x,z)
+\rho(x)E_{\tau}(x,z)\right\}h(x)dS_x.
\end{array}
$$
This yields
$$\displaystyle
{}^tY_{22}(\tau)h(z)
=\int_{\partial D}\left\{\frac{\partial}{\partial\nu_x}E_{\tau}(x,z)
+\rho(x)E_{\tau}(x,z)\right\}h(x)dS_x,\,\,z\in\partial D.
$$
From this and (\ref{3.6}) we obtain the desired formula.

\noindent
$\Box$

Define
\begin{equation*}\displaystyle
M(\tau)={}^tY_{22}(\tau)(I-{}^tY_{22}(\tau))^{-1}.
\end{equation*}
One can write
\begin{equation}\begin{array}{c}
\displaystyle
(I-{}^tY_{22}(\tau))^{-1}
\displaystyle
=I+{}^tY_{22}(\tau)+({}^tY_{22}(\tau))^2(I-{}^tY_{22}(\tau))^{-1}\\
\\
\displaystyle
=I+{}^tY_{22}(\tau)+{}^tY_{22}(\tau)M(\tau).
\end{array}
\label{3.16}
\end{equation}

Define, for an arbitrary $z\not=x$ and $x\in\partial D$
\begin{equation*}
\begin{array}{c}
\displaystyle
H_0(x,z)=\frac{\nu_x\cdot (z-x)}{\vert x-z\vert^2},
\\
\\
\displaystyle
H_1(x,z)=\frac{1}{\vert x-z\vert}
\left(\frac{\nu_x\cdot(z-x)}{\vert x-z\vert^2}
+\rho(x)\right).
\end{array}
\end{equation*}
Since
\begin{equation}\displaystyle
H(x,z,\tau)=\tau H_0(x,z)+H_1(x,z),
\label{3.17}
\end{equation}
from lemma 3.1 we have
\begin{equation}
\displaystyle
{}^tY_{22}(\tau)=M_0(\tau)+\tilde{M}(\tau), 
\label{3.18}
\end{equation}
where
\begin{equation}
\begin{array}{c}
\displaystyle
M_0(\tau)h(z)
=\frac{\tau}{2\pi}\int_{\partial D}e^{-\tau\vert x-z\vert}H_0(x,z)h(x)dS_x,\\
\\
\displaystyle
\tilde{M}(\tau)h(z)
=\frac{1}{2\pi}\int_{\partial D}e^{-\tau\vert x-z\vert}H_1(x,z)h(x)dS_x.
\end{array}
\label{3.19}
\end{equation}

Now set
\begin{equation}\displaystyle
M_1(\tau)=\tilde{M}(\tau)+{}^tY_{22}(\tau)M(\tau)
\label{3.20}
\end{equation}
and
\begin{equation}\displaystyle
F_j(x,p,\tau)
=e^{\tau\vert x-p\vert}
\left(M_j(\tau)\left(
\frac{e^{-\tau\vert \,\cdot\,-p\vert}}{\vert\,\cdot\,-p\vert}\right)
\right)(x),\,\,j=0,1.
\label{3.21}
\end{equation}
From (\ref{3.16}), (\ref{3.18}) and (\ref{3.20}) we have 
$(I-{}^tY_{22}(\tau))^{-1}=I+M_0(\tau)+M_1(\tau)$
and thus (\ref{3.12}) can be rewritten as
\begin{equation*}
F(x,p,\tau)
=\frac{1}{\vert p-x\vert}+F_0(x,p,\tau)+F_1(x,p,\tau).
\end{equation*}

\noindent
Substituting this and (\ref{3.17}) into (\ref{3.14}), we obtain

\begin{Prop}\label{Proposition 3.1}
The decomposition
$$\displaystyle
I_0(\tau,p)=\tau I_{0\,0}(\tau,p)+I_{0\,1}(\tau,p),
$$
is valid, where
$$\begin{array}{c}
\displaystyle
G_0(x,y,p,\tau)
=H^{+}(x,y,p)+2H_0(x,y)(F_0(x,p,\tau)+F_1(x,p,\tau)),\\
\\
\displaystyle
G_1(x,y,p,\tau)=\frac{H_1(x,p)}{\vert x-y\vert}+\frac{H_1(x,y)}{\vert x-p\vert}
+2H_1(x,y)(F_0(x,p,\tau)+F_1(x,p,\tau))
\end{array}
$$
and
$$\displaystyle
I_{0\,j}(\tau,p)
=\left(\frac{1}{2\pi}\right)^2
\int_{\partial\Omega}dS_y\varphi(y,\tau)
\int_{\partial D}e^{-\tau\,l_p(x,y)}G_j(x,y,p,\tau)dS_x,\,\,j=0,1.
$$
\end{Prop}

\subsection{Basic estimates of integral kernels}
\label{Basic estimates of integral kernels}

We introduce basic estimates of 
the integral kernels of the operators $M_0(\tau)$ and $M_1(\tau)$
introduced in (\ref{3.19}) and (\ref{3.20}).
To obtain the asymptotic behaviour of $I_0(\tau,p)$, these estimates of 
the kernels are essentially needed in our proof.  
In this subsection we always assume that $D$ is a bounded domain with 
the boundary $\partial D$ of class $C^{2,\,\alpha_0}$ 
with $0<\alpha_0\le 1$.

It is well known that 
there exists a positive constant $C$ such that for all 
$x, z\in\partial D$
\begin{equation}\displaystyle
\vert \nu_x-\nu_z \vert\le C\vert x-z\vert,\,\,
\vert\nu_x \cdot(x-z)\vert
\le C\vert x-z\vert^2.
\label{3.22}
\end{equation}
From (\ref{3.19}) and (\ref{3.22}), we see that 
the integral kernel $M_0(x,z, \tau)$ of the operator 
$M_0(\tau)$ is given by
\begin{equation}\displaystyle
M_0(x,z, \tau)
=\frac{\tau}{2\pi}e^{-\tau\vert x-z\vert}\,
\frac{\nu_z\cdot (x-z)}{\vert x-z\vert^2}
\label{3.23}
\end{equation}
and has the estimate
\begin{equation}\displaystyle
\vert M_0(x,z, \tau)\vert
\le C_0\, \tau\,e^{-\tau\vert x-z\vert},\,\,x,z\in\,\partial D,
\,\,\tau > 0.
\label{3.24}
\end{equation}

For $M_1(\tau)$ we can obtain the following result:

\begin{Thm}\label{Theorem 3.1}
Assume that $\partial{D}$ is strictly convex. Then 
there exist positive constants $C$ and $\mu_0\ge 1$ such that:
for all $\tau \geq \mu_0$
the operator $M_1(\tau)$ has an integral kernel 
$M_1(x,z ,\tau)$ which
is measurable for $(x,z)\in\partial D\times\partial D$, 
continuous for $x\not=z$ and has the estimate
\begin{equation}\displaystyle
\vert M_1(x,z, \tau)\vert
\le C e^{-\tau\vert x-z\vert}
\left(1+\frac{1}{\vert x-z\vert}
+\min\,\left\{\tau(\tau\vert x-z\vert^3)^{1/2},
\,\frac{1}{\vert x-z\vert^3}\right\}\right).
\label{3.25}
\end{equation}

\end{Thm}

\begin{Remark}\label{Remark 3.1}
Since $\min\,\{\sqrt{a}, a^{-1}\}\le 1$ for all $a>0$, from
(\ref{3.25}) we get
\begin{equation}\displaystyle
\vert M_1(x,z, \tau)\vert\le 
C\left(\tau+\frac{1}{\vert x-z\vert}\right)
e^{-\tau\vert x-z\vert}.
\label{3.26}
\end{equation}
\end{Remark}

These estimates are essential to obtain theorem \ref{Theorem 2.1}.
As is described in section \ref{Proof of Theorem 1.1.}, 
for a proof of theorem \ref{Theorem 3.1} is given in \cite{IK2}.

\subsection{The structure of ${\cal M}(p)$}
\label{The structure of {cal M}(p)}

The last of the preliminaries, 
we study the structure of the set ${\cal M}(p)$.

\begin{Prop}\label{Proposition 4.1}
Assume that $\partial D$ is of class $C^2$.
Then it holds that:

\noindent
(1) if $(x_0,y_0)\in {\cal M}(p)$, then $\nu_{y_0}=(y_0-x_0)/\vert y_0-x_0\vert$;

\noindent
(2) if $(x_0,y_0)\in {\cal M}_1(p)$, then $\nu_{x_0}$ has to be on the plane
determined by the three points $p$, $x_0$, $y_0$ and the angle between
$p-x_0$ and $\nu_{x_0}$ coincides with the angle between $y_0-x_0$ and $\nu_{x_0}$;

\noindent
(3)  the set ${\cal M}(p)$ has the decomposition
$$\displaystyle
{\cal M}(p)={\cal M}_1(p)\cup {\cal M}_2^+(p)\cup {\cal M}_2^{-}(p)\cup {\cal M}_g(p);
$$

\noindent
(4)  if $(x_0,y_0)\in {\cal M}_2^{+}(p)\cup {\cal M}_2^{-}(p)\cup {\cal M}_g(p)$, then there exists $t\in]0,\,1[$ such that
$x_0=(1-t)p+ty_0$.

Further assume that $D$ is strictly convex.
Then it holds that:

\noindent
(5) if $(x_0,y_0)\in {\cal M}_2^-(p)$, then there exists a unique $x_0^*\in {\cal G}^+(p)$
such that $(x_0^*,y_0)\in {\cal M}_2^{+}(p)$;

\noindent
(6) if $(x_0,y_0)\in {\cal M}_2^{+}(p)$, then there exists a unique $x_0^*\in {\cal G}^-(p)$
such that $(x_0^*,y_0)\in {\cal M}_2^{-}(p)$.

\end{Prop}
{\it\noindent Proof.}
Let $(x_0,y_0)\in {\cal M}(p)$.  Choose a system of local coordinates
$x=x(\sigma)$, $\sigma=(\sigma_1,\sigma_2)$ with $x_0=x(0)$ in a
neighbourhood of $x_0\in\partial D$. Similarly choose a system of
local coordinates $y=y(\theta)$, $\theta=(\theta_1,\theta_2)$ with
$y_0=y(0)$ in a neighbourhood of $y_0\in\partial\Omega$. Then the
function $\tilde{l_p}(\sigma,\theta)=l_p(x(\sigma),y(\theta))$
takes the local minimum at $(\sigma,\theta)=(0,0)$. Thus we have,
for all $j=1,2$
$$\displaystyle
\frac{\partial}{\partial\sigma_j}\tilde{l_p}(0,0)=0,\,\,
\frac{\partial}{\partial\theta_j}\tilde{l_p}(0,0)=0.
$$
Since
\begin{equation}
\displaystyle
\frac{\partial}{\partial\sigma_j}\tilde{l_p}(\sigma,\theta)=
\left(\frac{x-p}{\vert x-p\vert}+\frac{x-y}{\vert x-y\vert}\right)\cdot\frac{\partial x}{\partial\sigma_j}
\label{4.1}
\end{equation}
and
\begin{equation}
\displaystyle
\frac{\partial}{\partial\theta_j}\tilde{l_p}(\sigma,\theta)
=-\frac{x-y}{\vert x-y\vert}\cdot\frac{\partial y}{\partial\theta_j},
\label{4.2}
\end{equation}
we get
\begin{equation}
\displaystyle
\left(\frac{x_0-p}{\vert x_0-p\vert}+\frac{x_0-y_0}{\vert x_0-y_0\vert}\right)\cdot\frac{\partial x}{\partial\sigma_j}(0,0)=0
\label{4.3}
\end{equation}
and
\begin{equation}
\displaystyle
\frac{x_0-y_0}{\vert x_0-y_0\vert}\cdot\frac{\partial y}{\partial\theta_j}(0,0)=0.
\label{4.4}
\end{equation}
This last equality yields that $\nu_{y_0}$ and $(y_0-x_0)/\vert y_0-x_0\vert$ 
have to be parallel.
Assume that $\nu_{y_0}=-(y_0-x_0)/\vert y_0-x_0\vert$.
Then one can find a point $y_0'$ outside $\Omega$ that
is located on the segment $x_0y_0$.  Since $x_0\in\Omega$, one can find a point 
$y_0''\in\partial\Omega$
on the segment $x_0y_0'$.  Then we have $l_p(x_0,y_0'')<l_p(x_0,y_0)$.  
This is against $(x_0,y_0)\in {\cal M}(p)$. 
Therefore (1) has to be true.

Write
$$\displaystyle
\frac{p-x_0}{\vert p-x_0\vert}=\alpha\nu_{x_0}
+\beta\frac{\partial x}{\partial\sigma_1}(0,0)
+\gamma\frac{\partial x}{\partial\sigma_2}(0,0)
$$
and
$$\displaystyle
\frac{y_0-x_0}{\vert y_0-x_0\vert}=\alpha'\nu_{x_0}
+\beta'\frac{\partial x}{\partial\sigma_1}(0,0)
+\gamma'\frac{\partial x}{\partial\sigma_2}(0,0).
$$
Since $\nu_{x_0}\cdot\partial x/\partial\sigma_j(0,0)=0$, we have
$$\displaystyle
\alpha=\frac{p-x_0}{\vert p-x_0\vert}\cdot\nu_{x_0},\,\,
\alpha'=\frac{y_0-x_0}{\vert y_0-x_0\vert}\cdot\nu_{x_0}.
$$
From (\ref{4.3}) we get
the system of the equations for $\beta+\beta'$ and $\gamma+\gamma'$:
\begin{equation}
\displaystyle
\left(\begin{array}{c}
\displaystyle
\frac{\partial x}{\partial\sigma_1}^T(0,0)\\
\\
\displaystyle
\frac{\partial x}{\partial\sigma_2}^T(0,0)
\end{array}
\right)
\left(\begin{array}{cc}
\displaystyle \frac{\partial x}{\partial\sigma_1}(0,0) &
\displaystyle \frac{\partial x}{\partial\sigma_2}(0,0)
\end{array}\right)
\left(\begin{array}{c}
\displaystyle \beta+\beta'\\
\\
\displaystyle \gamma+\gamma'
\end{array}\right)
=
\left(\begin{array}{c}
\displaystyle 0\\
\\
\displaystyle 0
\end{array}\right).
\label{4.5}
\end{equation}
Since the vectors $\partial x/\partial\sigma_j(0,0)$, $j=1,2$ are linearly independent,
the coefficients matrix of (\ref{4.5}) is invertible and one gets
$\beta+\beta'=0$ and $\gamma+\gamma'=0$.  This yields
\begin{equation}\displaystyle
\frac{p-x_0}{\vert p-x_0\vert}+
\frac{y_0-x_0}{\vert y_0-x_0\vert}=(\alpha+\alpha')\nu_{x_0}.
\label{4.6}
\end{equation}
Moreover since the vectors $(p-x_0)/\vert p-x_0\vert$ and $(y_0-x_0)/\vert y_0-x_0\vert$ have the unit length,
$\beta^2=\beta'^2$, $\gamma^2=\gamma'^2$ and $\beta\gamma=\beta'\gamma'$,
we get
\begin{equation}\displaystyle
\vert\alpha\vert=\vert\alpha'\vert.
\label{4.7}
\end{equation}

If $(x_0,y_0)\in {\cal M}_1(p)$, then both $\alpha$ and $\alpha'$ are positive and
from (\ref{4.6}) and (\ref{4.7}) we obtain $\alpha=\alpha'$
and 
\begin{equation}\displaystyle
\frac{p-x_0}{\vert p-x_0\vert}+
\frac{y_0-x_0}{\vert y_0-x_0\vert}
=2\alpha\nu_{x_0}.
\label{like the law of geometrical optics}
\end{equation}
This coincides with the law of reflection of the light and yields (2).

For the proof of (3) it suffices to prove
that the set ${\cal M}(p)$ is contained in 
${\cal M}_1(p)\cup {\cal M}_2^+(p)\cup
{\cal M}_2^{-}(p)\cup {\cal M}_g(p)$.  
We employ a contradiction argument.
Assume that there exists a 
$(x_0,y_0)\in {\cal M}(p)\setminus({\cal M}_1(p)\cup
{\cal M}_2^+(p)\cup {\cal M}_2^{-}(p)\cup {\cal M}_g(p))$. 
Since $(x_0,y_0)$ does not belong to ${\cal M}_g(p)$, 
we get $x_0\in {\cal G}^+(p)$ or $x_0\in {\cal G}^-(p)$.
Consider the case when $x_0\in {\cal G}^+(p)$.  Since the $(x_0,y_0)$
does not belong to ${\cal M}_1(p)\cup {\cal M}_2^{+}(p)$, we have
$(y_0-x_0)\cdot \nu_{x_0}=0$.  Then from (\ref{4.7}) we have
$(p-x_0)\cdot\nu_{x_0}=0$.  Contradiction. Next consider the case
when $x_0\in {\cal G}^-(p)$.  Since the $(x_0,y_0)$ does not belong to
${\cal M}_2^-(p)$, we have $(y_0-x_0)\cdot\nu_{x_0}\le 0$. 
(\ref{4.7}) yields
$(y_0-x_0)\cdot\nu_{x_0}<0$. Since $\partial D$ is $C^2$ at $x_0$,
one can find an open ball $B$ contained in $D$ such that $\partial
B\cap\partial D=\{x_0\}$.  Therefore the set of all numbers
$t\in\,]0,\,1[$ such that $(1-s)x_0+sy_0\in D$ for all $0<s<t$, is
not empty.  Denote by $t^*$ the least upper bound of the set. It
is easy to see that $0<t^*<1$ and the point
$x_0'=(1-t^*)x_0+t^*y_0\in\partial D$ and $x_0\not=x_0'$.
If the points $x_0'$, $p$ and $x_0$ form a triangle, then by
the triangle inequality we have $\vert p-x_0\vert+\vert x_0-y_0\vert>\vert
p-x_0'\vert+\vert x_0'-y_0\vert$. If they do not form a triangle, then $y_0$ has to be
on the segment $px_0$ since $(p - x_0)\cdot\nu_{x_0} < 0$ and $(y_0 - x_0)\cdot\nu_{x_0} < 0$.
Since $\vert p-x_0\vert=\vert p-x_0'\vert+\vert x_0-x_0'\vert$
and $\vert x_0-y_0\vert=\vert x_0-x_0'\vert+\vert x_0'-y_0\vert$,
we get $\vert p-x_0\vert+\vert x_0-y_0\vert=\vert p-x_0'\vert
+\vert x_0'-y_0\vert+2\vert x_0-x_0'\vert
>\vert p-x_0'\vert+\vert x_0'-y_0\vert$.
This against $(x_0,y_0)\in
{\cal M}(p)$.  Contradiction.  This completes the proof of (3).

The proof of (4) starts with the simple fact: if 
$(x_0,y_0)\in {\cal M}_2^{+}(p)\cup {\cal M}_2^{-}(p)$, then
the numbers $(p-x_0)\cdot\nu_{x_0}$ and $(y_0-x_0)\cdot\nu_{x_0}$ have different signature.
This together with (\ref{4.7}) yields $\alpha+\alpha'=0$ in (\ref{4.6}).   
If $(x_0,y_0)\in {\cal M}_g(p)$, then
$\alpha=0$ and (\ref{4.7}) gives again $\alpha+\alpha'=0$ in (\ref{4.6}).  
In any case we get
\begin{equation}\displaystyle
\frac{p-x_0}{\vert p-x_0\vert}+
\frac{y_0-x_0}{\vert y_0-x_0\vert}=0.
\label{4.8}
\end{equation}
Therefore $t=\vert p-x_0\vert/l_p(x_0,y_0)$ gives the desired conclusion.
(5) and (6) are trivial.

\noindent
$\Box$

%
\setcounter{equation}{0}
\section{Proof of theorem \ref{Theorem 2.1}}
\label{Proof of the main result}
%

Given $\delta>0$ define
$$\displaystyle
{\cal G}_{\delta}(p)=\{x\in\partial D\,\vert\,
\text{dist}(x,{\cal G}(p))\ge\delta\},\,\,
{\cal G}_{\delta}^{\pm}(p)={\cal G}_{\delta}(p)\cap {\cal G}^{\pm}(p).
$$
In this section first we state two crucial lemmas needed for establishing 
the asymptotic formula for $I_0(\tau,p)$.

The first lemma is concerned with the asymptotic behaviour of the
amplitudes of the integrals in proposition \ref{Proposition 3.1} and the proof is
given in section \ref{Asymptotic behaviour of F_j(x, p, lambda)}.

\begin{Lemma}\label{Lemma 5.1}
There exists a positive constant $\mu_0$ such that the following assertions are true.

\noindent
(1)  There exists a positive constant $C$ such that if 
$x\in\partial D$
and $\tau\ge\mu_0$, then
$$\displaystyle
\vert F_j(x,p,\tau)\vert\le C\tau,\,\, j=0,1.
$$

\noindent
(2)  Given $\delta>0$ there exists a positive constant $C_{\delta}$ such that if
$x\in {\cal G}_{\delta}^+(p)$ and $\tau\ge\mu_0$, then
$$\displaystyle
\vert F_j(x,p,\tau)\vert\le C_{\delta}\tau^{-1},\,\,j=0,1.
$$

\noindent
(3)  Given $\delta>0$ there exists a positive constant 
$C_{\delta}$ such that if $x\in {\cal G}_{\delta}^-(p)$, 
and $\tau\ge\mu_0$, then
$$\displaystyle
\vert F_1(x,p,\tau)\vert\le C_{\delta}\tau^{-1}.
$$

\noindent
(4)  Given $\delta>0$ there exists a positive constant 
$C_{\delta}$ such that if $x\in {\cal G}_{\delta}^{-}(p)$
and $\tau \ge \mu_0$, then
$$\displaystyle
\left\vert F_0(x,p,\tau)+\frac{1}{\vert x-p\vert}\right\vert
\le C_{\delta}\tau^{-\alpha_0/2}.
$$

\end{Lemma}

The following lemma gives the asymptotic behaviour of an integral 
with an exponential weight and the idea behind the derivation is 
called the Laplace method.

\begin{Lemma}\label{Lemma 5.2} 
Let $U$ be an arbitrary open set
of ${\Bbb R}^n$.  Let $f\in C^{2,\alpha_0}(\overline U)$ and satisfy
at a point $x_0\in U$, for all $x\in\overline U\setminus\{x_0\}$
$f(x)>f(x_0)$ and $\text{det}\,(\text{Hess}\,(f)(x_0))>0$. Then
given $\varphi\in C^{0,\alpha_0}(\overline U)$ it holds that
$$\displaystyle
\int_U e^{-\tau f(x)}\varphi(x)dx
=\frac{e^{-\tau f(x_0)}}{\sqrt{\text{det}\,(\text{Hess}\,(f)(x_0))}}
\left(\frac{2\pi}{\tau}\right)^{n/2}
\left(\varphi(x_0)+\Vert\varphi\Vert_{C^{0,\alpha_0}(\overline U)}
O(\tau^{-\alpha_0/2})\right).
$$
Moreover there exists a positive constant $C$ such that, 
for all $\tau \geq 1$
$$\displaystyle
\left\vert\int_Ue^{-\tau f(x)}\varphi(x)dx\right\vert\le
\frac{C e^{-\tau f(x_0)}}{\tau^{n/2}}
\Vert\varphi\Vert_{C(\overline U)}.
$$
\end{Lemma}

The proof of this lemma is given in Appendix A.
We now give a proof of theorem \ref{Theorem 2.1}.
Since ${\cal M}(p)$ is a finite set, one can write
$$\displaystyle
{\cal M}(p)=\{(x^{(j)},y^{(j)})\,\vert\,j=1,2,\cdots,N\}.
$$
However, by (5) and (6) of proposition \ref{Proposition 4.1} 
the counting number of the set ${\cal M}_2^{+}(p)$ coincides
with that of ${\cal M}_2^{-}(p)$.  Then (3) of proposition \ref{Proposition 4.1}
yields that 
the counting number of the set ${\cal M}(p)\setminus {\cal M}_1(p)$ has to 
be an even number.  Hence one can write
$$\begin{array}{c}
\displaystyle
{\cal M}_1(p)=\{(x^{(j)},y^{(j)})\,\vert\,j=1,2,\cdots,N_1\},\\
\\
\displaystyle
{\cal M}_2^{+}(p)
=\{(x^{(j)}, y^{(j)})\,\vert\,j=N_1+1,\cdots, N_1+N_2\},\\
\\
\displaystyle
{\cal M}_2^{-}(p)
=\{(x^{(j)}, y^{(j)})\,\vert\,j=N_1+N_2+1,\cdots,N_1+2N_2\},
\end{array}
$$
where $x^{(j)}=(x^{(j+N_2)})^*$, $j=N_1+1,\cdots,N_1+N_2$ and $N=N_1+2N_2$.

From the second equation in (\ref{3.2}) we have
$$\displaystyle
\psi(x,\tau)=(I-Y_{22}(\tau))^{-1}Y_{21}(\tau)\varphi(x,\tau).
$$
Then from the first equation of (\ref{3.2}) we obtain the equation of $\varphi$ only:
$$\displaystyle
\left\{I-Y_{11}(\tau)-Y_{12}(\tau)
(I-Y_{22}(\tau))^{-1}Y_{21}(\tau)\right\}\varphi(x,\tau)
=g(x,\tau).
$$
Since $\Vert Y_{ij}(\tau)\Vert=O(\tau^{-1})$ as 
$\tau\longrightarrow\infty$,
it follows from the equation that
\begin{equation}\displaystyle
\varphi(y,\tau)
=g(y,\tau)+O(\tau^{-1})\Vert g(\,\cdot\,,\tau)\Vert_{C(\partial\Omega)}
\label{5.1}
\end{equation}
as $\tau\longrightarrow\infty$ uniformly for 
$y\in\partial\Omega$.

Given $\delta>0$ set
$$\displaystyle U_{\delta}(x^{(j)})
=\{x\in\partial D\,\vert\,\vert x-x^{(j)}\vert<\delta\},\,\,
V_{\delta}(y^{(j)})
=\{y\in\partial\Omega\,\vert\,\vert y-y^{(j)}\vert<\delta\}.
$$
One can choose a sufficiently small $\delta>0$ such that, for $j=1,\cdots, N_1+2N_2$
$U_{2\delta}(x^{(j)})\cap {\cal G}(p)=\emptyset$
and $\left(U_{2\delta}(x^{(j)})\times V_{2\delta}(y^{(j)})\right)\cap {\cal M}(p)
=\{(x^{(j)},y^{(j)})\}$.
Moreover since $l_p(x,y)>l(p,D)$ for all $(x,y)$ in the compact set
$(\partial D\times\partial\Omega)\setminus
\left(\cup_{j=1}^{N_1+2N_2}U_{\delta/3}(x^{(j)})\times V_{\delta/3}(y^{(j)})\right)$,
one can find a positive constant $c_0$ such that
$$\displaystyle
l_p(x,y)\ge l(p,D)+c_0\,\,\text{for}\,(x,y)
\in (\partial D\times\partial\Omega)\setminus\left(\cup_{j=1}^{N_1+2N_2}
U_{\delta/3}(x^{(j)})\times V_{\delta/3}(y^{(j)})\right).
$$
From this, (1) of lemma \ref{Lemma 5.1}, (\ref{5.1}) 
and proposition \ref{Proposition 3.1} one gets, for $k=0,1$
\begin{equation}\displaystyle
I_{0k}(\tau,p)
=\left(\frac{1}{2\pi}\right)^2
\sum_{j=1}^{N_1+2N_2} I_{0kj}(\tau,p)
+e^{-\tau l(p,D)}O(e^{-c_0\tau/2})\Vert g\Vert_{C(\partial\Omega)}.
\label{5.2}
\end{equation}
Here, for $j=1,\cdots, N_1+2N_2$
$$\displaystyle
I_{0kj}(\tau,p)
=\int_{V_{\delta}(y^{(j)})}dS_y\varphi(y,\tau)
\int_{U_{\delta}(x^{(j)})}e^{-\tau l_p(x,y)}
\Psi_j(x,y)G_k(x,y,p,\tau)dS_x
$$
and $\Psi_j\in C_0^{2}(U_{\delta}(x^{(j)})\times V_{\delta}(y^{(j)}))$ 
is a cut-off function with
$\Psi_j(x,y)=1$ in $U_{\delta/2}(x^{(j)})\times V_{\delta/2}(y^{(j)})$ and
$\Psi_j(x,y)=0$ in $(U_{2\delta/3}(x^{(j)})\times V_{2\delta/3}(y^{(j)}))^c$.

We study the asymptotic behaviour of $I_{0kj}(\tau,p)$.  Choose local coordinate 
systems 
$x=s^{(j)}(\sigma)$ with $x^{(j)}=s^{(j)}(0)$ for $U_{\delta}(x^{(j)})$
and $y=\tilde{s}^{(j)}(\tilde{\sigma})$
with $y^{(j)}=\tilde{s}^{(j)}(0)$ for $V_{\delta}(y^{(j)})$.
Set $\tilde{\Psi}_j(\sigma,\tilde{\sigma})=\Psi_j(s^{(j)}(\sigma),
\tilde{s}^{(j)}(\tilde{\sigma}))$,
$$\displaystyle
J_j(\sigma,\tilde{\sigma})
=\sqrt{\text{det}\,
\left(\frac{\partial s^{(j)}}{\partial\sigma_p}(\sigma)\cdot
\frac{\partial s^{(j)}}{\partial\sigma_q}(\sigma)\right)
\,
\text{det}\,\left(\frac{\partial\tilde{s}^{(j)}}{\partial\tilde{\sigma}_p}
(\tilde{\sigma})\cdot
\frac{\partial\tilde{s}^{(j)}}{\partial\tilde{\sigma}_q}(\tilde{\sigma})\right)}
$$
and $\tilde{l_p}^{(j)}(\sigma,\tilde{\sigma})
=l_p(s^{(j)}(\sigma),\tilde{s}^{(j)}(\tilde{\sigma}))$.
A change of variables gives the expression
$$\displaystyle
I_{0kj}(\tau,p)
=\int_{{\Bbb R}^4} e^{-\tau\tilde{l_p}^{(j)}(\sigma,\tilde{\sigma})}
\varphi(\tilde{s}^{(j)}(\tilde{\sigma}),\tau)
\tilde{\Psi}_j(\sigma,\tilde{\sigma})
G_k(s^{(j)}(\sigma),\tilde{s}^{(j)}(\tilde{\sigma}),p,\tau)
J_j(\sigma,\tilde{\sigma})d\sigma d\tilde{\sigma}.
$$

Since the function $x\longmapsto (p-x)\cdot\nu_x$ is continuous and 
$\overline{U_{\delta}(x^{(j)})}\cap {\cal G}(p)=\emptyset$,
we have $\overline{U_{\delta}(x^{(j)})}\subset {\cal G}^{+}(p)$ for 
$j=1,\cdots, N_1+N_2$;
$\overline{U_{\delta}(x^{(j)})}\subset {\cal G}^{-}(p)$ for 
$j=N_1+N_2+1,\cdots, N_1+2N_2$.

Consider the case when $j=1,\cdots,N_1+N_2$.
It follows from (2) of lemma \ref{Lemma 5.1}
$$\displaystyle
G_0(x,y,p,\tau)=H^+(x,y,p)+O(\tau^{-1}),
\,\,
G_1(x,y,p,\tau)=O(1)
$$
as $\tau\longrightarrow\infty$ uniformly for $(x,y)\in\,
\overline{U_{\delta}(x^{(j)})}\times\partial\Omega$.
Since we have $\text{Hess}\,(\tilde{l_p}^{(j)})(0,0)>0$, 
from these estimate, (\ref{5.1}) and lemma \ref{Lemma 5.2},
we obtain
\begin{equation}\displaystyle
I_{01j}(\tau,p)
=e^{-\tau l(p, D)}\Vert g(\,\cdot, \,\tau)\Vert_{C(\partial\Omega)}O(\tau^{-2}),
\label{5.3}
\end{equation}
and
\begin{equation}
\begin{array}{c}
\displaystyle
I_{00j}(\tau,p)
=\frac{J_j(0,0)e^{-\tau l(p,D)}}
{\sqrt{\text{det}\,(\text{Hess}\,(\tilde{l_p}^{(j)})(0,0))}}
\left(\frac{2\pi}{\tau}\right)^2\\
\\
\displaystyle
\times
\left(g(y^{(j)},\tau)H^+(x^{(j)},y^{(j)},p)
+O(\tau^{-\alpha_0/2})\Vert g(\,\cdot\, ,\tau)
\Vert_{C^{0,\alpha_0}(\partial\Omega)}\right).
\end{array}
\label{5.4}
\end{equation}

Next consider the case when $j=N_1+N_2,\cdots, N_1+2N_2$.  
From (3) and (4) of lemma \ref{Lemma 5.1} we get
$$\displaystyle
G_0(x,y,p,\tau)=H^{-}(x,y,p)+O(\tau^{-\alpha_0/2}),\,\,
G_1(x,y,p,\tau)=O(1)
$$
as $\tau\longrightarrow\infty$ uniformly for $(x,y)\in\,
\overline{U_{\delta}(x^{(j)})}\times\partial\Omega$.
From these estimates, (\ref{5.1}) and lemma \ref{Lemma 5.2} 
we obtain
\begin{equation}\displaystyle
I_{01j}(\tau,p)
=e^{-\tau l(p,D)}\Vert g(\,\cdot\, ,\tau)\Vert_{C(\partial\Omega)}O(\tau^{-2})
\label{5.5}
\end{equation}
and
\begin{equation}\begin{array}{c}
\displaystyle
I_{00j}(\tau,p)
=\frac{J_j(0,0)e^{-\tau l(p,D)}}
{\sqrt{\text{det}\,(\text{Hess}\,(\tilde{l_p}^{(j)})(0,0))}}
\left(\frac{2\pi}{\tau}\right)^2\\
\\
\displaystyle
\times
\left(g(y^{(j)},\tau)H^-(x^{(j)},y^{(j)},p)
+O(\tau^{-\alpha_0/2})\Vert g(\,\cdot\, ,\tau)
\Vert_{C^{0,\alpha_0}(\partial\Omega)}\right).
\end{array}
\label{5.6}
\end{equation}

From proposition \ref{Proposition 3.1}, (\ref{5.2}) to 
(\ref{5.6}) and the fact that 
$H^+(x^{(j)}, y^{(j)},p)=0$ for $j=N_1+1,\cdots, N_1+N_2$ (see (\ref{4.8})),
we obtain the desired asymptotic formula (\ref{2.5}) for $I_0(\tau,p)$.
The coefficients $C(x_0,y_0)$ in (\ref{2.6}) for $(x_0,y_0)=(x^{(j)}, y^{(j)})$ 
are given by
$$\displaystyle
\frac{J_j(0,0)}
{\sqrt{\text{det}\,(\text{Hess}\,(\tilde{l_p}^{(j)})(0,0))}}
$$
and thus positive.  This completes the proof of theorem \ref{Theorem 2.1}.

%
\setcounter{equation}{0}
\section{Asymptotic behaviour of $F_j(x, p, \tau)$}
\label{Asymptotic behaviour of F_j(x, p, lambda)}
%

In this section, we prove lemma \ref{Lemma 5.1}. In the first two subsections, 
we prepare properties of the broken path and estimates of boundary integrals used to
show lemma \ref{Lemma 5.1}. The last subsection, we give a proof of 
lemma \ref{Lemma 5.1} 
using the estimates of the integral kernels of $M_0(\tau)$ and $M_1(\tau)$
given in (\ref{3.24}) and theorem \ref{Theorem 3.1}, respectively.    
\par
Throughout this section, we always assume that $\partial D$ is 
of class $C^{2,\alpha_0}$ with $0<\alpha_0\le 1$.  We denote by $B(x,r)$ 
the open ball centered at $x$ with radius $r$.

\subsection{Properties of the broken path }
\label{Properties of the broken path }

The aim of this subsection is
to study the behaviour of the function
$$\displaystyle
l_{(p,\,x)}(z)\equiv\vert p-z\vert+\vert z-x\vert
\,\,
$$
with the independent variable $z\in\partial D$, and given 
$p\in{\Bbb R}^3\setminus\overline\Omega$ and $x \in\partial D$.

We start with describing the following well known facts.

\begin{Lemma}\label{Lemma 6.1}
There exists $0<r_0$ such that, for all $x\in\partial D$, 
$\partial D\cap B(x,2r_0)$ can be represented as a graph of a function 
on the tangent plane of $\partial D$ at $x$, that is,
there exist an open neighbourhood $U_x$ of $(0,0)$ in ${\Bbb R}^2$
and a function $g=g_x\in C^{2,\alpha_0}({\Bbb R}^2)$ with $g(0,0)=0$ and 
$\nabla g(0,0)=0$ such that the map
$$
U_x\ni\,\sigma=(\sigma_1,\sigma_2)\mapsto
x+\sigma_1 e_1+\sigma_2e_2-g(\sigma_1,\sigma_2)\nu_x\in \partial D\cap B(x,2r_0)
$$
gives a system of local coordinates around $x$, where $\{e_1, e_2\}$ is 
an orthogonal basis for $T_{x}(\partial D)$.  Moreover
the norm $\Vert g\Vert_{C^{2,\alpha_0}({\Bbb R}^2)}$ has an
upper bound independent of $x\in\partial D$.
\end{Lemma}

\noindent
In this paper we call this system of coordinates the standard system 
of local coordinates around $x$.

The following lemma plays an important role in the proof of 
lemma \ref{Lemma 5.1}.

\begin{Lemma}\label{Lemma 6.2}
Assume that $\partial D$ is strictly convex.
If $x\in {\cal G}^+(p)\cup {\cal G}(p)$, then the function $l_{(p,\,x)}(z)$, 
$z\in\partial D$ attains the minimum only at $z=x$.  
If $x\in {\cal G}^{-}(p)$, then the points on $\partial D$
that attain the minimum are given by only two points $z=x, x^*$.
Moreover the following statements are true.

\noindent
(i)  Given $\delta>0$ there exists a positive constant $C_{\delta}$ such that
if $x\in {\cal G}_{\delta}^{+}(p)$, then for all $z\in\partial D$
we have
$$\displaystyle
l_{(p,\,x)}(z)\ge\vert p-x\vert+C_{\delta}\vert z-x\vert.
$$

\noindent
(ii)  Given $\delta>0$ there exists a constant $0<\delta'_0\le\delta$
such that if $x\in {\cal G}_{\delta}^{-}(p)$, then
$\vert x-x^*\vert\ge 2\delta_0'$.  Further, for any $0<\delta'\le\delta_0'$,
there exists a positive constant $C_{\delta'}$ such that, for all $x\in {\cal G}_{\delta}^-(p)$
and $z\in\partial D\setminus B(x^*,\delta')$ we have
$$\displaystyle
l_{(p,\,x)}(z)\ge\vert p-x\vert+C_{\delta'}\vert z-x\vert.
$$

\noindent
(iii)
Given $\delta>0$ there exist positive constants $C_{\delta}$ and $C'_{\delta}$
such that, if $0<\delta'\le C'_{\delta}$, then for all $x\in {\cal G}^-_{\delta}(p)$ 
and
$z\in\,\partial D\cap \overline{B(x^*,\delta')}$,
$$\displaystyle
l_{(p,\,x)}(z)\ge \vert p-x\vert+C_{\delta}\vert z-x^*\vert^2.
$$
\end{Lemma}
{\it\noindent Proof.}
It is clear that $\min_{z\in\,\partial D}\,l_{(p,x)}(z)=l_{(p,x)}(x)$.
Let $z\in\partial D$ be a point such that
$l_{(p,x)}(z)=l_{(p,x)}(x)$.

Consider the case $x\in {\cal G}^+(p)\cup {\cal G}(p)$.  Assume that $z\not=x$.
Since $\vert p-z\vert+\vert z-x\vert=\vert p-x\vert$,
$z$ has to be on the line segment determined by $p$ and $x$.
Since $(p-x)\cdot\nu_x\ge 0$, we have $(z-x)\cdot\nu_x\ge 0$.  
On the other hand, since $\partial D$ is strictly convex
and $z\not=x$, one gets $(z-x)\cdot\nu_x<0$. This is a contradiction. Thus $z=x$.

Next consider the case $x\in {\cal G}^-(p)$.  Assume that $z\not=x$.
Similarly to above one knows that $z$ is located on the line
segment determined by $p$ and $x$ and thus gets $z=x^*$. Therefore
the set of all points $z$ that attain the minimum of
$l_{(p,\,x)}(\,\cdot\,)$ is contained in the set $\{x, x^*\}$.
However since $l_{(p,\,x)}(x^*)=l_{(p,x)}(x)$, the function
$l_{(p,\,x)}(\,\cdot\,)$ really attains the minimum at $z=x, x^*$.

Now we give a proof of (i).  Let $z\not=x$.
We have
$$\begin{array}{c}
\displaystyle
\vert p-z\vert^2
=\vert p-x\vert^2+\vert z-x\vert^2
-2(p-x)\cdot (z-x)\\
\\
\displaystyle
=\left\{\vert p-x\vert-\vert z-x\vert
\left(\frac{z-x}{\vert z-x\vert}\cdot\frac{p-x}{\vert p-x\vert}\right)\right\}^2
+\vert z-x\vert^2\left\{1-
\left(\frac{z-x}{\vert z-x\vert}\cdot\frac{p-x}{\vert p-x\vert}\right)^2\right\}.
\end{array}
$$
This yields
$$\displaystyle
\vert p-z\vert\ge \vert p-x\vert-\vert z-x\vert
\left(\frac{z-x}{\vert z-x\vert}\cdot\frac{p-x}{\vert p-x\vert}\right).
$$
From this we obtain the estimate
\begin{equation}\displaystyle
l_{(p,\,x)}(z)\ge \vert p-x\vert+
\vert z-x\vert\left(1-\frac{z-x}{\vert z-x\vert}\cdot\frac{p-x}{\vert p-x\vert}\right).
\label{6.1}
\end{equation}
Let $z'$ be the orthogonal projection of $z$ onto $T_x(\partial D)$.  We see that
$(z-z')\cdot (p-x)\le 0$ since $(z-z')\cdot\nu_x\le 0$, $(p-x)\cdot\nu_x\ge 0$
and $z-z'$ is parallel to $\nu_x$.  It follows from this that
\begin{equation}
\displaystyle
\frac{z-x}{\vert z-x\vert}\cdot
\frac{p-x}{\vert p-x\vert}
=\frac{z-z'}{\vert z-x\vert}\cdot\frac{p-x}{\vert p-x\vert}
+\frac{z'-x}{\vert z-x\vert}\cdot\frac{p-x}{\vert p-x\vert}
\le
\frac{z'-x}{\vert z-x\vert}\cdot\frac{p-x}{\vert p-x\vert}.
\label{6.2}
\end{equation}
First consider the case $(z'-x)\cdot (p-x)\ge 0$.  
Since $\vert z-x\vert\ge\vert z'-x\vert$, from (\ref{6.2}) we have
\begin{equation}\displaystyle
\frac{z-x}{\vert z-x\vert}\cdot
\frac{p-x}{\vert p-x\vert}
\le
\frac{z'-x}{\vert z'-x\vert}\cdot\frac{p-x}{\vert p-x\vert}.
\label{6.3}
\end{equation}

Let $p'$ be the orthogonal projection of $p$ onto $T_x(\partial D)$.  Since
$z'-x$ and $p'-x$ are parallel to $T_x(\partial D)$,
we have
\begin{equation}\displaystyle
\left\vert\frac{p-x}{\vert p-x\vert}-
\frac{z'-x}{\vert z'-x\vert}\right\vert
\ge
\left\vert\frac{p-x}{\vert p-x\vert}
-\frac{p'-x}{\vert p-x\vert}\right\vert
\\
\\
\displaystyle
=\frac{(p-x)\cdot\nu_x}{\vert p-x\vert}.
\label{6.4}
\end{equation}

Set
$$\displaystyle
A_{\delta}\equiv\inf_{x\in {\cal G}^+_{\delta}(p)}
\frac{(p-x)\cdot\nu_x}{\vert p-x\vert}>0.
$$
From (\ref{6.4}) we have
\begin{equation}\displaystyle
\frac{z'-x}{\vert z'-x\vert}\cdot\frac{p-x}{\vert p-x\vert}
=1-\frac{1}{2}\left\vert\frac{p-x}{\vert p-x\vert}
-\frac{z'-x}{\vert z'-x\vert}\right\vert^2
\le 1-\frac{1}{2}A_{\delta}^2.
\label{6.5}
\end{equation}

Now from (\ref{6.1}), (\ref{6.3}) and (\ref{6.5}) we obtain
$$
\displaystyle l_{(p,\,x)}(z)\ge \vert p-x\vert
+\frac{A_{\delta}^2}{2}\vert z-x\vert
$$
provided $(z'-x)\cdot (p-x)\ge 0$.  If $(z'-x)\cdot (p-x)<0$,  (\ref{6.2}) gives
$$
\displaystyle \frac{z-x}{\vert z-x\vert}\cdot
\frac{p-x}{\vert p-x\vert}<0.
$$
Then from (\ref{6.1}) we have
$$
\displaystyle l_{(p,\,x)}(z)\ge \vert p-x\vert+\vert z-x\vert.
$$
Therefore (i) holds for $C_{\delta}=\min\,\{A_{\delta}^2,\,2\}/2$.

Next we give a proof of (ii).  It is clear that the map: 
${\cal G}^-(p)\ni x\mapsto x^*\in {\cal G}^+(p)$
is continuous.  Since the set ${\cal G}^-_{\delta}(p)$ is compact and 
$\vert x-x^*\vert>0$ for all $x\in {\cal G}^-_{\delta}(p)$, we have
$$\displaystyle
B_{\delta}\equiv\inf_{x\in {\cal G}^-_{\delta}(p)}\vert x-x^*\vert>0.
$$
Then  $\delta_0'=\min\,\{B_{\delta}/2, \delta\}$ satisfies the desired condition.
Next we prove that
\begin{equation}\displaystyle
D_{\delta'}\equiv\sup_{x\in {\cal G}_{\delta}^{-}(p)}\,
\sup_{z\in\,(\partial D\setminus\{x\})\setminus\overline{B(x^*,\delta')}}
\frac{z-x}{\vert z-x\vert}\cdot\frac{x^*-x}{\vert x^*-x\vert}<1.
\label{6.6}
\end{equation}
If this is not true, then the compactness of 
${\cal G}_{\delta}^-(p)$ and $\partial D$ yields 
the existence of points $x_0\in {\cal G}_{\delta}^{-}(p)$
and $z_0\in\partial D$ and sequences $\{x_n\}$ with 
$x_n\in {\cal G}_{\delta}^{-}(p)$ and $\{z_n\}$ with 
$z_n\in\,(\partial D\setminus\{x_n\})\setminus
\overline{B(x_n^*,\delta')}$ such that, 
as $n\longrightarrow\infty$ 
$x_n\longrightarrow x_0$, $z_n\longrightarrow z_0$ and
\begin{equation}\displaystyle
\frac{z_n-x_n}{\vert z_n-x_n\vert}\cdot
\frac{x_n^*-x_n}{\vert x_n^*-x_n\vert}\longrightarrow 1.
\label{6.7}
\end{equation}
Moreover, one may assume that the unit vectors $(z_n-x_n)/\vert
z_n-x_n\vert$ converges to a unit vector $\vartheta$.  Since
$\vert x_n-x_n^*\vert\ge 2\delta_0'$, from the continuity of the
map ${\cal G}_{\delta}^-(p)\ni x\longmapsto x^*\in\partial D$ we have
$x_0\not=x_0^*$.  Thus from (\ref{6.7}) we obtain
$$\displaystyle
\vartheta\cdot\frac{x_0^*-x_0}{\vert x_0^*-x_0\vert}=1.
$$
This gives $\vartheta=(x_0^*-x_0)/\vert x_0^*-x_0\vert$ and 
since $\partial D$ is strictly convex,
we obtain $\vartheta\cdot\nu_{x_0}<0$.

Consider the case when $z_0=x_0$.  From (\ref{3.22})
we obtain $\vartheta\cdot\nu_{x_0}=0$.  This is a contradiction.

Next consider the case when $z_0\not=x_0$.  In this case we obtain
$$
\displaystyle
\frac{z_0-x_0}{\vert z_0-x_0\vert}
=\frac{x_0^*-x_0}{\vert x_0^*-x_0\vert}.
$$
This yields that $z_0$ is located on the line determined by $x_0$ and $x_0^*$.
Since $\partial D$ is strictly convex, we have $z_0=x_0^*$.  However, we have also
$\vert z_0-x_0^*\vert\ge\delta'$.  Contradiction.

Therefore (\ref{6.6}) is valid.
Since $(x^*-x)/\vert x^*-x\vert=(p-x)/\vert p-x\vert$, from (\ref{6.1}) we have
$$\displaystyle
l_{(p,\,x)}(z)\ge \vert p-x\vert+
\vert z-x\vert\left(1-\frac{z-x}{\vert z-x\vert}\cdot
\frac{x^*-x}{\vert x^*-x\vert}\right).
$$
Now the final conclusion of (ii) is true for $C_{\delta'}=1-D_{\delta'}$.

Finally we give a proof of (iii).
Since $\vert p-x\vert=\vert p-x^*\vert+\vert x^*-x\vert$, we have
\begin{equation}
\displaystyle
l_{(p,x)}(z)-\vert p-x\vert=(\vert p-z\vert-\vert p-x^*\vert)
+(\vert z-x\vert-\vert x^*-x\vert).
\label{6.8}
\end{equation}
Set $\xi=z-x^*$.  We have
\begin{equation}\begin{array}{c}
\displaystyle
\vert p-z\vert-\vert p-x^*\vert
=\frac{\vert\xi\vert^2-2(p-x^*)\cdot\xi}
{\vert p-z\vert+\vert p-x^*\vert},\\
\\
\displaystyle
\vert z-x\vert-\vert x^*-x\vert
=\frac{\vert\xi\vert^2+2(x^*-x)\cdot\xi}
{\vert z-x\vert+\vert x^*-x\vert}.
\end{array}
\label{6.9}
\end{equation}
From (\ref{6.8}) and (\ref{6.9}) we have the expression
\begin{equation}\displaystyle
l_{(p,\,x)}(z)-\vert p-x\vert
=
\left(\frac{1}{\vert p-z\vert+\vert p-x^*\vert}
+\frac{1}{\vert z-x\vert+\vert x^*-x\vert}\right)\vert\xi\vert^2
+2R(z)\cdot\xi, 
\label{6.10}
\end{equation}
where
$$\displaystyle
R(z)=\frac{x^*-p}{\vert p-z\vert+\vert p-x^*\vert}
+\frac{x^*-x}{\vert z-x\vert+\vert x^*-x\vert}.
$$
Since
$$\begin{array}{c}
\displaystyle
\frac{x^*-p}{\vert p-z\vert+\vert p-x^*\vert}
-\frac{x^*-p}{2\vert p-x^*\vert}
=-\frac{(\vert p-z\vert-\vert p-x^*\vert)(x^*-p)}
{2(\vert p-z\vert+\vert p-x^*\vert)\vert p-x^*\vert},\\
\\
\displaystyle
\frac{x^*-x}
{\vert z-x\vert+\vert x^*-x\vert}-\frac{x^*-x}{2\vert x^*-x\vert}
=-\frac{(\vert z-x\vert-\vert x^*-x\vert)(x^*-x)}
{2(\vert z-x\vert+\vert x^*-x\vert)\vert x^*-x\vert},
\end{array}
$$
from (\ref{6.9}) one gets
$$\begin{array}{c}
\displaystyle
\frac{x^*-p}{\vert p-z\vert+\vert p-x^*\vert}
=\frac{1}{2}
\frac{x^*-p}{\vert p-x^*\vert}\\
\\
\displaystyle
+\frac{\{(p-x^*)\cdot\xi\}(x^*-p)}{(\vert p-z\vert+\vert p-x^*\vert)^2\vert p-x^*\vert}
-\frac{1}{2}\frac{\vert\xi\vert^2(x^*-p)}{(\vert p-z\vert+\vert p-x^*\vert)^2
\vert p-x^*\vert}
\end{array}
$$
and
$$\begin{array}{c}
\displaystyle
\frac{x^*-x}{\vert z-x\vert+\vert x^*-x\vert}
=\frac{1}{2}\frac{x^*-x}{\vert x^*-x\vert}\\
\\
\displaystyle
-\frac{\{(x^*-x)\cdot\xi\}(x^*-x)}
{(\vert z-x\vert+\vert x^*-x\vert)^2\vert x^*-x\vert}
-\frac{1}{2}\frac{\vert\xi\vert^2(x^*-x)}
{(\vert z-x\vert+\vert x^*-x\vert)^2\vert x^*-x\vert}.
\end{array}
$$
Since
$$
\frac{x^*-x}{\vert x^*-x\vert}+\frac{x^*-p}{\vert x^*-p\vert}=0,
$$
it follows that
$$\begin{array}{c}
\displaystyle
R(z)\cdot\xi
=-\frac{\{(x^*-p)\cdot\xi\}^2}
{(\vert p-z\vert+\vert p-x^*\vert)^2\vert p-x^*\vert}
-\frac{\{(x^*-x)\cdot\xi\}^2}
{(\vert z-x\vert+\vert x^*-x\vert)^2\vert x^*-x\vert}\\
\\
\displaystyle
-\frac{1}{2}\frac{\vert\xi\vert^2(x^*-p)\cdot\xi}{(\vert p-z\vert+\vert p-x^*\vert)^2
\vert p-x^*\vert}
-\frac{1}{2}\frac{\vert\xi\vert^2(x^*-x)\cdot\xi}
{(\vert z-x\vert+\vert x^*-x\vert)^2\vert x^*-x\vert}.
\end{array}
$$
Using the facts
\begin{equation}
\displaystyle
\inf_{(x,\,z)\in\, {\cal G}^{-}_{\delta}(p)\times\partial D}\,\vert p-z\vert
+\vert p-x^*\vert>0,
\,\,
\inf_{(x,\,z)\in\, {\cal G}^{-}_{\delta}(p)\times\partial D}\,\vert z-x\vert
+\vert x^*-x\vert>0,
\label{6.11}
\end{equation}
from (\ref{6.10}) we obtain
\begin{equation}\displaystyle
l_{(p,\,x)}(z)=\vert p-x\vert+K(z)\xi\cdot\xi+O(\vert\xi\vert^3)
\label{6.12}
\end{equation}
uniformly for $x\in {\cal G}^-_{\delta}(p)$,
where
$$\begin{array}{c}
\displaystyle
K(z)=\left(\frac{1}{\vert p-z\vert+\vert p-x^*\vert}
+\frac{1}{\vert z-x\vert+\vert x^*-x\vert}\right)I_3\\
\\
\displaystyle
-\frac{2(x^*-p)\otimes(x^*-p)}
{(\vert p-z\vert+\vert p-x^*\vert)^2\vert p-x^*\vert}
-\frac{2(x^*-x)\otimes (x^*-x)}
{(\vert z-x\vert+\vert x^*-x\vert)^2\vert x^*-x\vert}.
\end{array}
$$
Set
$$\displaystyle
\vartheta=\frac{x^*-p}{\vert x^*-p\vert}.
$$
Then
$$
\frac{x^*-x}{\vert x^*-x\vert}=-\vartheta
$$
and we have
$$\begin{array}{c}
\displaystyle
K(z)=\left(\frac{1}{\vert p-z\vert+\vert p-x^*\vert}
+\frac{1}{\vert z-x\vert+\vert x^*-x\vert}\right)I_3\\
\\
\displaystyle
-2\left(\frac{\vert x^*-p\vert}{(\vert p-z\vert+\vert p-x^*\vert)^2}
+\frac{\vert x^*-x\vert}{(\vert z-x\vert+\vert x^*-x\vert)^2}\right)
\vartheta\otimes\vartheta.
\end{array}
$$

Let $z=z(\sigma)=x^*+\sigma_1 e_1+\sigma_2 e_2-g(\sigma)\nu_{x^*}$ be the standard local coordinate system
around $x^*$ for $z\in\partial D\cap B(x^*,2r_0)$.  From lemma \ref{Lemma 6.1}
we know that, for a suitable constant $C$ depending only on $\partial D$
we have $\vert\sigma\vert\le\vert\xi\vert\le C\vert\sigma\vert$.
Using (\ref{6.11}) together with the following facts
$$
\displaystyle
\inf_{x\in\, {\cal G}^{-}_{\delta}(p)}\,\vert p-x^*\vert>0,\,\,
\inf_{x\in\, {\cal G}^{-}_{\delta}(p)}\,\vert x^*-x\vert>0,
$$
we have
$$
\displaystyle
\frac{1}{\vert p-z\vert+\vert p-x^*\vert}=\frac{1}{2\vert p-x^*\vert}
+O(\vert\xi\vert),\,\,
\frac{1}{\vert z-x\vert+\vert x^*-x\vert}=\frac{1}{2\vert x^*-x\vert}
+O(\vert\xi\vert)
$$
uniformly for $x\in\,{\cal G}^{-}_{\delta}(p)$.  These yield
$$\displaystyle
K(z)=\frac{1}{2}
\left(\frac{1}{\vert p-x\vert}+\frac{1}{\vert x^*-x\vert}\right)
\left(I_3-\vartheta\otimes\vartheta\right)
+O(\vert\xi\vert).
$$
Since $\xi=\sigma_1 e_1+\sigma_2 e_2+O(\vert\xi\vert^2)$, we obtain
$$\displaystyle
K(z)\xi\cdot\xi
=\frac{1}{2}\left(\frac{1}{\vert p-x\vert}+\frac{1}{\vert x^*-x\vert}\right)
\left(I_2-\vartheta'\otimes\vartheta'\right)\sigma\cdot\sigma
+O(\vert\xi\vert^3)
$$
where $\vartheta'=(\vartheta\cdot e_1,\,\vartheta\cdot e_2)^T$.
Here we note that the eigenvalues of the $2\times 2$-matrix 
$I_2-\vartheta'\otimes\vartheta'$
are given by $1$ and $1-\vert \vartheta\cdot e_1\vert^2
-\vert\vartheta\cdot e_2\vert^2=\vert\vartheta\cdot\nu_{x^*}\vert^2$.
Therefore we conclude that
$$\displaystyle
K(z)\xi\cdot\xi
\ge \frac{1}{2}\left(\frac{1}{\vert p-x\vert}+\frac{1}{\vert x^*-x\vert}\right)
\vert\vartheta\cdot\nu_{x^*}\vert^2\vert\sigma\vert^2
+O(\vert\xi\vert^3).
$$
Since
$$
\displaystyle
\inf_{x\in\,{\cal G}^-_{\delta}(p)}\frac{p-x^*}{\vert p-x^*\vert}\cdot\nu_{x^*}>0,
$$
from (\ref{6.12}) we obtain the desired conclusion.

\noindent
$\Box$

\begin{Remark}\label{Remark 6.1}
From the proof of (iii) we obtain the expression
$$\begin{array}{c}
\displaystyle
l_{(p,x)}(z)-\vert p-x\vert\\
\\
\displaystyle
=K(z)\xi\cdot\xi
-\frac{\vert\xi\vert^2(x^*-p)\cdot\xi}{(\vert p-z\vert+\vert p-x^*\vert)^2\vert p-x^*\vert}
-\frac{\vert\xi\vert^2(x^*-x)\cdot\xi}
{(\vert z-x\vert+\vert x^*-x\vert)^2\vert x^*-x\vert}
\end{array}
$$
with $\xi=z-x^*$.
To show theorem \ref{Theorem 2.1}, we need this equality.
\end{Remark}

\subsection{Estimates of integrals on the boundary $\partial{D}$}
\label{Estimates of integrals on the boundary partial{D}}

To show lemma \ref{Lemma 5.1}, we need the following estimates:
\begin{Lemma}\label{Lemma 6.3}
Let $r_0$ be the same as that of lemma \ref{Lemma 6.1}.
There exists a positive constant $C$ depending only on $\partial D$ such that

\noindent
(i)  for all $x\in\partial D$, $0<\rho'_0\le r_0$, 
$\tau>0$, $0\le k<2$
$$\displaystyle
\int_{B(x,\rho'_0)\,\cap\,\partial D}
\frac{e^{-\tau\vert x-z\vert}}{\vert x-z\vert^k}\,dS_z
\le\frac{C}{2-k}\,\min\,\{\tau^{-2+k},\,(\rho'_0)^{2-k}\};
$$

\noindent
(ii)  for all $x\in\partial D$, $\tau>0$, $0\le k<2$
$$\displaystyle
\int_{\partial D}\frac{e^{-\tau\vert x-z\vert}}{\vert x-z\vert^k}\,dS_z
\le\frac{C}{2-k}\tau^{-(2-k)}\left(1+\frac{\tau^{2-k}e^{-\tau r_0}}{r_0^k}\right).
$$
\end{Lemma}
{\it\noindent Proof.}
Let $z=s(\sigma)$ be the standard system of local coordinates around $x$ with
$\vert\sigma\vert^2+g(\sigma)^2<(2r_0)^2$.
We have
$$\begin{array}{c}
\displaystyle
\int_{B(x,\rho'_0)\cap\partial D}\,\frac{e^{-\tau\vert x-z\vert}}{\vert x-z\vert^k}\,dS_z
=\int_{\vert\sigma\vert^2+g(\sigma)^2<(\rho_0')^2}
\frac{e^{-\tau\sqrt{\vert\sigma\vert^2+g(\sigma)^2}}}
{(\vert\sigma\vert^2+g(\sigma)^2)^{k/2}}
\sqrt{1+\vert\nabla g(\sigma)\vert^2}d\sigma\\
\\
\displaystyle
\le C\int_0^{\rho'_0}\int_0^{2\pi}\frac{e^{-\tau r}}{r^k}\,rdrd\theta
\le 2\pi C\int_0^{\rho'_0}e^{-\tau r}r^{1-k}dr.
\end{array}
$$
Here note that
$$\displaystyle
\int_0^{\rho'_0}e^{-\tau r}r^{1-k}dr
\le\int_0^{\rho'_0}r^{1-k}dr
=\frac{(\rho'_0)^{2-k}}{2-k}
$$
and
$$\begin{array}{c}
\displaystyle
\int_0^{\rho'_0}e^{-\tau r}r^{1-k}dr
=\tau^{k-2}\int_0^{\tau\rho'_0}e^{-r}r^{1-k}dr\\
\\
\displaystyle
\le\tau^{k-2}\int_0^{\infty}e^{-r}r^{1-k}dr
=\tau^{k-2}\left(1+\frac{1}{2-k}\right)\le
\frac{3}{2-k}\,\tau^{k-2}.
\end{array}
$$
This proves (i).  To verify (ii) we compute
$$\displaystyle
\int_{\partial D\setminus B(x,r_0)} \frac{e^{-\tau\vert x-z\vert}}{\vert x-z\vert^k}dS_z
\le e^{-r_0\tau}\int_{\partial D}\frac{1}{r_0^k}dS_z
\le \frac{C}{r_0^k}e^{-r_0\tau}.
$$
From this and (i) for $\rho'_0=r_0$ we obtain (ii).
This completes the proof of lemma \ref{Lemma 6.3}.

\noindent
$\Box$

\subsection{Proof of lemma \ref{Lemma 5.1}}
\label{Proof of Lemma 5.1}

We start with the expression for $F_j(x,p,\tau)$ for $j=0,1$ 
(see (\ref{3.21})):
$$\displaystyle
F_j(x,p,\tau)
=e^{\tau\vert x-p\vert}\int_{\partial D}M_j(x,z, \tau)
\frac{e^{-\tau\vert z-p\vert}}{\vert z-p\vert}dS_z,\,\,
x\in\partial D.
$$
For the case $j=0$ $M_0(x,y, \tau)$ is given by (\ref{3.23}) and
the case $j=1$ is a consequence of theorem \ref{Theorem 3.1}.

First we prove (1) of lemma \ref{Lemma 5.1}.
From (\ref{3.24}) and (\ref{3.26}) we get
\begin{equation}
\begin{array}{c}
\displaystyle
\vert F_j(x,p,\tau)\vert
\le C e^{\tau\vert x-p\vert}
\int_{\partial D}\vert M_j(x,z, \tau)\vert\,
\frac{e^{-\tau\vert z-p\vert}}{\vert z-p\vert}\,dS_z\\
\\
\displaystyle
\le C \int_{\partial D}\left(\tau+\frac{1}{\vert x-z\vert}\right)
\frac{e^{-\tau (\vert x-z\vert+\vert z-p\vert-\vert x-p\vert)}}{\vert z-p\vert}\,dS_z.
\end{array}
\label{6.13}
\end{equation}
Since $\vert x-z\vert+\vert z-p\vert\ge\vert x-p\vert$, 
the right-hand side of (\ref{6.13}) has the bound
$$
\displaystyle
\frac{C}{\text{dist}\,(p,\,\partial D)}\int_{\partial D}
\left(\tau+\frac{1}{\vert x-z\vert}\right)\,dS_z.
$$
Applying the argument for the proof of (ii) in lemma \ref{Lemma 6.3} 
to the integral above,
we see that
$$\displaystyle
\sup_{x\in\partial D}\,\int_{\partial D}\frac{dS_z}{\vert x-z\vert}<\infty.
$$
Thus one concludes that (1) is true.

Second we prove (2) of lemma \ref{Lemma 5.1}.
Consider the case when $x\in {\cal G}_{\delta}^+(p)$.  One can apply 
(i) of lemma \ref{Lemma 6.2} to the integrand in
the right-hand side of (\ref{6.13}) and get
$$\begin{array}{c}
\displaystyle
\int_{\partial D}\left(\tau+\frac{1}{\vert x-z\vert}\right)
\frac{e^{-\tau (\vert x-z\vert+\vert z-p\vert-\vert x-p\vert)}}{\vert z-p\vert}\,dS_z
\le\int_{\partial D}\left(\tau+\frac{1}{\vert x-z\vert}\right)
\frac{e^{-\tau C_{\delta}\vert z-x\vert}}{\vert z-p\vert}\,dS_z\\
\\
\displaystyle
\le\frac{C}{\text{dist}\,(p,\,\partial D)}
\int_{\partial D}\left(\tau+\frac{1}{\vert x-z\vert}\right)
e^{-\tau C_{\delta}\vert z-x\vert}\,dS_z.
\end{array}
$$
Applying (ii) of lemma \ref{Lemma 6.3} to the integral of the right-hand side
above, one gets
$$\displaystyle
\int_{\partial D}\left(\tau+\frac{1}{\vert x-z\vert}\right)
e^{-\tau C_{\delta}\vert z-x\vert}\,dS_z
\le C(\tau\cdot\tau^{-2}+\tau^{-1}).
$$
Thus this together with (\ref{6.13}) yields that (2) is true.

Third we prove (3) of lemma \ref{Lemma 5.1}.
By (ii) and (iii) of lemma \ref{Lemma 6.2}, one can find $C_{\delta}>0$ 
and $\delta'>0$ such that, for all $x\in {\cal G}_{\delta}^-(p)$
\begin{equation}
\displaystyle
\vert x-x^*\vert\ge 2\delta',
\label{6.14}
\end{equation}
\begin{equation}\displaystyle
\vert x-z\vert+\vert z-p\vert
\ge \vert p-x\vert+C_{\delta}\vert z-x\vert,\,z\in\,\partial D\setminus B(x^*,\delta'),
\label{6.15}
\end{equation}
and
\begin{equation}\displaystyle
\vert x-z\vert+\vert z-p\vert
\ge \vert p-x\vert+C_{\delta}\vert z-x^*\vert^2,\,z\in\,\partial D\cap B(x^*,\delta').
\label{6.16}
\end{equation}
We decompose $\partial D$ into two parts 
$\partial{D}\cap B(x^*,\delta')$ and $\partial D\setminus B(x^*,\delta')$.
Then we have
$$
\displaystyle
\vert F_1(x,p,\tau)\vert\le I+II
$$
where
$$\begin{array}{c}
\displaystyle
I= e^{\tau\vert x-p\vert}\int_{\partial D\cap B(x^*,\delta')}
\vert M_1(x,z, \tau)
\vert\,\frac{e^{-\tau\vert z-p\vert}}{\vert z-p\vert}\,dS_z,\\
\\
\displaystyle
II=
e^{\tau\vert x-p\vert}\int_{\partial D\setminus B(x^*,\delta')}
\vert M_1(x,z, \tau)
\vert\,\frac{e^{-\tau\vert z-p\vert}}{\vert z-p\vert}\,dS_z.
\end{array}
$$

From (\ref{6.14}) we see that if $z\in\partial D\cap B(x^*,\delta')$, 
then $\vert z-x\vert\ge \delta'$.
This together with (\ref{3.25}) gives
$$\displaystyle
\vert M_1(x,z, \tau)\vert
\le C e^{-\tau\vert x-z\vert}\left(1+\frac{1}{\delta'}+\frac{1}{\delta'^3}\right).
$$
It follows from this and (\ref{6.16}) that
$$
\displaystyle
I\le \frac{C}{\text{dist}\,(p,\,\partial D)}
\int_{\partial D\cap B(x^*,\delta')} e^{-\tau C_{\delta}\vert z-x^*\vert^2}\,dS_z.
$$
Note that $\delta'$ can be arbitrary small and thus one may assume that
$\delta'<2r_0$, where $r_0$ is given in lemma \ref{Lemma 6.1}.  
Using the standard local coordinates around $x^*$, one obtains
$$\displaystyle
\int_{\partial D\cap B(x^*,\delta')}
e^{-\tau C_{\delta}\vert z-x^*\vert^2}\,dS_z
\le C\int_0^{\delta'}\int_0^{2\pi} e^{-\tau C_{\delta} r^2}rdrd\theta
\le C'\tau^{-1}
$$
and this thus yields $I\le C\tau^{-1}$.
For the estimation of $II$ we make use of (\ref{3.26}) and (\ref{6.15}).  This gives
$$\begin{array}{c}
\displaystyle
II\le \int_{\partial D\setminus B(x^*,\delta')}
\left(\tau+\frac{1}{\vert x-z\vert}\right)
\frac{e^{-\tau C_{\delta}\vert x-z\vert}}{\vert z-p\vert}dS_z\\
\\
\displaystyle
\le\frac{C}{\text{dist}\,(p,\,\partial D)}\int_{\partial D}
\left(\tau+\frac{1}{\vert x-z\vert}\right)
\,e^{-\tau C_{\delta}\vert x-z\vert}dS_z
\le C\tau^{-1}.
\end{array}
$$
Therefore $II\le C\tau^{-1}$ and this completes the proof of (3).

Finally we prove (4) of lemma \ref{Lemma 5.1}.
From (\ref{3.23}) one gets
\begin{equation}\displaystyle
F_0(x,p,\tau)
=\frac{\tau}{2\pi}\int_{\partial D\cap B(x^*,\delta')}
\frac{\nu_z\cdot(x-z)}{\vert x-z\vert^2}\,
\frac{e^{-\tau(\vert x-z\vert+\vert z-p\vert-\vert x-p\vert)}}{\vert z-p\vert}dS_z
+R(\tau),
\label{6.17}
\end{equation}
where
$$\displaystyle
R(\tau)=\frac{\tau}{2\pi}\int_{\partial D\setminus B(x^*,\delta')}
\frac{\nu_z\cdot(x-z)}{\vert x-z\vert^2}\,
\frac{e^{-\tau(\vert x-z\vert+\vert z-p\vert-\vert x-p\vert)}}{\vert z-p\vert}dS_z.
$$

For $r_0$ in lemma \ref{Lemma 6.1}, one can choose $\delta'$ above in such a way that
$\delta'\le 2r_0$ and (\ref{6.14})-(\ref{6.16}) are also satisfied.
A combination of the second inequality of (\ref{3.22}) and (\ref{6.15}) yields
\begin{equation}
\begin{array}{c}
\displaystyle
\vert R(\tau)\vert
\le\frac{C\tau}{\text{dist}(p,\,\partial D)}
\int_{\partial D} e^{-\tau C_{\delta}\vert x-z\vert}dS_z
=O(\tau^{-1}).
\end{array}
\label{6.18}
\end{equation}
Denote by $I(\tau)$ the first integral of the right-hand side of (\ref{6.17}).  
Let $z=s(\sigma)=x^*+\sigma_1\,e_1+\sigma_2\,e_2-g(\sigma)\,\nu_{x^*}$ be 
the standard local coordinates around $x^*$.
Choose a function $\chi\in C^{\infty}({\Bbb R}^2)$ such that
$\chi=1$ near $\vert s(\sigma)-x^*\vert\le \delta'/8$ and $\chi=0$ for 
$\vert s(\sigma)-x^*\vert\ge \delta'/4$.
Set
$$\displaystyle
L(\sigma)=\vert p-s(\sigma)\vert+\vert s(\sigma)-x\vert-\vert x-p\vert.
$$
It follows from (\ref{6.16}) that if 
$\delta'/4\ge\vert s(\sigma)-x^*\vert\ge\delta'/8$, then
$\displaystyle L(\sigma)\ge C_{\delta}(\delta'/8)^2$.
This yields
$$\displaystyle
I(\tau)
=\int_{{\Bbb R}^2}\chi(\sigma)\,\frac{\nu_{s(\sigma)}
\cdot(x-s(\sigma))}{\vert x-s(\sigma)\vert^2}\,
\frac{e^{-\tau L(\sigma)}}{\vert s(\sigma)-p\vert}
\sqrt{1+\nabla g(\sigma)^2}d\sigma
+O\left(e^{-C_{\delta}(\delta'/8)^2\,\tau}\right).
$$

Here we compute $\text{det}\,(\text{Hess}(L)(0))$.  
From remark \ref{Remark 6.1}, one can easily obtain
\begin{align*}
\frac{\partial^2}{\partial\sigma_i\partial\sigma_j}L(0)=2K(x^*)e_i\cdot e_j
&=\left(\frac{1}{\vert p-x^*\vert}+\frac{1}{\vert x-x^*\vert}\right)
\left(I_3-\vartheta\otimes\vartheta\right)e_i\cdot e_j
\nonumber
\\
&=\left(\frac{1}{\vert p-x^*\vert}+\frac{1}{\vert x-x^*\vert}\right)
\left(\delta_{ij}-(\vartheta\cdot e_i)(\vartheta\cdot e_j)\right)
\end{align*}
where $\vartheta=(p-x^*)/\vert p-x^*\vert$.  This together with the equation 
$\vert p-x^*\vert+\vert x-x^*\vert=\vert p-x\vert$
gives
\begin{align*}
\displaystyle
\text{det}\,(\text{Hess}\,(L)\,(0))
=\left(\frac{1}{\vert p-x^*\vert}+\frac{1}{\vert x-x^*\vert}\right)^2
\vert \vartheta\cdot\nu_{x^*}\vert^2
=\left(\frac{\vert p-x\vert}{\vert p-x^*\vert\,\vert x-x^*\vert}\right)^2
\vert\vartheta\cdot\nu_{x^*}\vert^2.
\end{align*}
Since $\vartheta\cdot\nu_{x^*}>0$, we obtain
$$\displaystyle
\sqrt{\text{det}\,(\text{Hess}\,(L)(0))}
=\frac{\vert p-x\vert}{\vert p-x^*\vert\,\vert x-x^*\vert}\vartheta\cdot\nu_{x^*}.
$$

Set
$$\displaystyle
\Phi(\sigma)
=\chi(\sigma)\,\frac{\nu_{s(\sigma)}\cdot(x-s(\sigma))}{\vert x-s(\sigma)\vert^2}\,
\frac{e^{-\tau L(\sigma)}}{\vert s(\sigma)-p\vert}\sqrt{1+\nabla g(\sigma)^2}.
$$
Since $(p-x^*)/\vert p-x^*\vert=-(x^*-x)/\vert x-x^*\vert$,
$\Phi(0)$ has the form
$$\displaystyle
\Phi(0)
=-\frac{\vartheta\cdot\nu_{x^*}}{\vert p-x^*\vert\vert x-x^*\vert}.
$$
This yields
$$\displaystyle
\frac{\Phi(0)}
{\sqrt{\text{det}\,(\text{Hess}(L)(0))}}
=-\frac{1}{\vert p-x\vert}.
$$
Now we are ready to apply lemma \ref{Lemma 5.2} 
to the integral $I(\tau)$.
The result is
$$\begin{array}{c}
\displaystyle
I(\tau)
=
\frac{e^{-\tau L(0)}}{\sqrt{\text{det}\,(\text{Hess}(L)(0))}}
\left(\frac{2\pi}{\tau}\right)^{2/2}
\left(\Phi(0)+O(\tau^{-\alpha_0/2})\right)+
O(\tau^{-\infty})\\
\\
\displaystyle
=-\frac{2\pi}{\tau}
\frac{1}{\vert p-x\vert}
+O(\tau^{-\alpha_0/2-1}).
\end{array}
$$

\noindent
From this together with (\ref{6.17}) and (\ref{6.18}) we obtain the desired conclusion.

\noindent
$\Box$

%
\setcounter{equation}{0}
\section{Sufficient conditions and examples}
\label{Sufficient conditions}
%

It is curious to know when assumptions of 
theorems \ref{Theorem 1.1} and \ref{Theorem 2.1}
are satisfied. We can give sufficient conditions 
to ensure that a point 
$(x_0,y_0)\in {\cal M}(p)\setminus {\cal M}_g(p)$
is a non-degenerate critical point of $l_p$ on $\partial D\times\partial\Omega$.
The conditions are given by using the Weingarten map of $C^2$ surfaces 
$S \subset \R^3$. Assume that $S$ is the $C^2$ boundary of a
bounded open set like $\partial\Omega$ and $\partial{D}$.
Let $\nu_x$ be the unit outer normal of $S$ at $x \in S$.
For a tangential vector field $v \in T_x(S)$ to $S$ at $x \in S$, 
the Weingarten map
${\cal A}_{S,x} $ is defined by 
${\cal A}_{S,x}(v) = D_v\nu_x$. By the original local coordinate 
of $\R^3$, $\nu_x$ is given as $\nu_x = (\nu_1(x), \nu_2(x), \nu_3(x))$,
and ${\cal A}_{S,x}$ is expressed as follows:
$$
{\cal A}_{S,x}(v) = \big(v_x(\nu_1), v_x(\nu_2), v_x(\nu_3)\big)
\qquad(v \in T_x(S)).
$$
From the definition, we can show that ${\cal A}_{S,x}$ is a linear
map on the tangent space $T_x(S)$, and ${\cal A}_{S,x} = 1/R$
when $S$ is a ball with radius $R > 0$.

\begin{Prop}\label{Proposition 2.1.}
Let $(x_0,y_0)\in {\cal M}(p)\setminus {\cal M}_g(p)$.
Assume that the Weingarten map ${\cal A}_{\partial\Omega,y_0}$ of $\partial\Omega$ at $y_0$ satisfies
\begin{equation}\displaystyle
{\cal A}_{\partial\Omega,y_0}<\frac{1}{l(p, D)}I.
\label{2.9}
\end{equation}
Then $(x_0,y_0)$ should be a non-degenerate critical point of $l_p$ on $\partial D\times\partial\Omega$.
\end{Prop}

\begin{Remark}\label{Remark 2.2.}
For every point $(x_0,y_0)\in{\cal M}_1(p)$, 
a sufficient condition (\ref{2.9}) for
non-degenerateness can be relaxed 
(cf. proposition \ref{Proposition 4.2}).
\end{Remark}

\vskip 1em

We need to check ${\cal M}_g(p)=\emptyset$ to apply 
theorems \ref{Theorem 1.1}, \ref{Theorem 2.1} and 
proposition \ref{Proposition 2.1.}
About this, we introduce a sufficient condition 
to satisfy  ${\cal M}_2^+(p)\cup {\cal M}_2^{-}(p)\cup{\cal M}_g(p)=\emptyset$.

\begin{Prop}\label{Proposition 2.2.}
If the set ${\cal L}(p)=\{y\in\partial\Omega\,\vert\,(y-p)/\vert y-p\vert\cdot\nu_y=1\}$ consists of a single point, 
then ${\cal M}_2^{+}(p)\cup {\cal M}_2^{-}(p)\cup {\cal M}_g(p)=\emptyset$.
\end{Prop}

\vskip 1em

\noindent
Note that ${\cal L}(p) \neq \emptyset$ since every point 
$ y_0 \in \partial\Omega$ attaining local maximum of 
the function $ \partial\Omega \ni y \to \vert y - p \vert $ belongs to 
${\cal L}(p)$.
If $\Omega$ is a ball, it is clear the assumption of 
proposition \ref{Proposition 2.2.} is satisfied.
However, even if $\Omega$ is convex, ${\cal L}(p)$ does not always consist of a single 
point.
For example, consider the case that $\partial\Omega$ contains a part of the 
sphere with 
the center $p$ and the radius $r = \max_{y \in \partial\Omega}\vert y - p \vert$.

From propositions \ref{Proposition 2.1.} and 
\ref{Proposition 2.2.} we can give examples for
corollary \ref{Corollary 2.1} deduced by theorems \ref{Theorem 1.1}
and \ref{Theorem 2.1} in sections \ref{Introduction} 
and \ref{Proof of Theorem 1.1.}.
We begin with introducing the following corollary:

\begin{Cor}\label{Corollary 2.2.}
Let $\Omega$ be the open ball with radius $R$ centered at the origin.  
Assume that $\partial D$ is strictly convex
and there exists a $\eta>0$ such that
$D$ contains the open ball with radius $R/2+\eta$ centered at the origin.
Let $p\in{\Bbb R}^3\setminus\overline\Omega$ satisfy 
${\rm dist}(p, \partial\Omega) < 2\eta$
Let $f(y,t)$ be the function of $(y,t)\in\partial\Omega\times\,]0,\,T[$ 
having the form $\tilde f(y)\varphi(t)$,
where $\tilde f\in C^{0,\alpha}(\partial\Omega)$ with 
$\tilde f(y)\not=0$ for all 
$y\in\partial\Omega$;
$\varphi\in L^2(0,\,T)$ satisfying the following condition:
there exists a $\tau>0$ such that
$$\displaystyle
0<\lim_{\tau\longrightarrow\infty}\tau^{\tau}\left\vert\int_0^Te^{-\tau^2t}\varphi(t)dt\right\vert<\infty.
$$
Then the formula (\ref{1.13}) is valid.
\end{Cor}

{\it\noindent\it Proof.}
Set $\epsilon={\rm dist}(p, \partial\Omega)$.
Since ${\rm dist}(p, \partial{D}) \le 
{\rm dist}(p, \partial{B_{R/2+\eta}})$ and
$$\displaystyle
{\rm dist}(p, \partial{B_{R/2+\eta}})
=\frac{R}{2}-\eta+\epsilon,
$$
we have
$$\displaystyle
2{\rm dist}(p, \partial{D})-{\rm dist}(p, \partial\Omega)
\le R-(2\eta-\epsilon).
$$
Since $\epsilon<2\eta$, this together with inequality
$
2{\rm dist}(p, \partial{D})-{\rm dist}(p, \partial\Omega)
\ge l(p, D)
$
which can be easily verified 
yields $l(p, D)<R$.
Now from this, propositions \ref{Proposition 2.1.}, 
\ref{Proposition 2.2.} and corollary \ref{Corollary 2.1} 
we obtain the desired conclusion.

\noindent
$\Box$

\par
Now, we give a simple example of a pair of $\Omega$ and $D$
in which the minimum length $l(p, D)$ can be obtained by 
the indicator function $I(\tau, p)$.


\begin{Example}\label{Example 2.1.}  
Let $R>r>0$.  Let $\Omega$ and $D$ be the open balls with radius $R$
and $r$, respectively centered at a common point $q$. 
Let $p$ be
an arbitrary point outside $\Omega$ with
$\text{dist}(p, \,\partial\Omega)=h>0$. Then by 
proposition \ref{Proposition 2.2.} 
one knows that
${\cal M}_g(p)={\cal M}_2^{+}(p)={\cal M}_2^{-}(p)=\emptyset$ 
and ${\cal M}(p)={\cal M}_1(p)$.
Let $(x_0,y_0)\in {\cal M}_1(p)$. Since $l(p, D)$ is the minimum of $l_p$, 
the function $ y \mapsto \vert y - x_0 \vert $ takes a local minimum at $y_0$.
This implies that $y_0-x_0$ and $\nu_{y_0}$ are parallel (precisely, we have
$\nu_{y_0} = \vert y_0-x_0 \vert^{-1}(y_0-x_0) $ as in (1) of 
proposition \ref{Proposition 4.1}).
Since
$D$ and $\Omega$ are spheres having the common center, the point
$x_0$ has to be on the line determined by $p$ and $q$.
Then one gets
${\cal M}_1(g)=\{(q+r(p-q)/\vert p-q\vert, q+R(p-q)/\vert p-q\vert)\}$ and 
$l(p,D)=h+2(R-r)$.
Assume
that we {\it know} $r_0\in\,]R/2,R[$ such that $r>r_0$.  This
$r_0$ can be considered as an a-priori information about unknown
$r$.  Choose $p$ in such a way that $h/2<r_0-R/2$.  Since
${\cal A}_{\partial\Omega,y_0}=(1/R)I$, the condition (\ref{2.9}) 
is satisfied.
In this case, from (\ref{2.6}) we see that 
the condition (\ref{2.7}) becomes
$\liminf_{\tau\longrightarrow\infty}\tau^{\tau}\vert
g(q+R(p-q)/\vert p-q\vert,\tau)\vert>0$.
Therefore if only this condition and corresponding one to (\ref{2.8}) 
are satisfied, one can 
extract the quantity $h+2(R-r)$ from (\ref{1.13}).
\end{Example}

As is mentioned in remark \ref{Remark 2.2.}, 
assumption (\ref{2.9}) in 
proposition \ref{Proposition 2.1.} 
can be relaxed.
Using this fact, we can cover other example 
containing example \ref{Example 2.1.},
which also justifies the fact that theorem \ref{Theorem 1.1} 
is considered as a three-dimensional 
analogue of (\ref{1.5}) 
(for this example, see subsection 
\ref{A sufficient condition of positive definiteness of l_p 
at (x_0,y_0) in {cal M}_1(p)}).

\subsection{Proof of proposition \ref{Proposition 2.2.}}

We give a proof of proposition \ref{Proposition 2.2.} in here.

\noindent
Step 1.  We prove: if $y'\in\partial\Omega$ attains
the maximum of the function $f(y)=\vert y-p\vert$, $y\in\partial\Omega$, then $y'\in {\cal L}(p)$.

\noindent Using a local coordinate at $y'$, we see that the
vectors $(y'-p)/\vert y'-p\vert$ and $\nu_{y'}$ are parallel each
other.  Assume that $(y'-p)/\vert y'-p\vert=-\nu_{y'}$. Since
$\partial\Omega$ is $C^2$, one can find a point $y''$ on
$\partial\Omega$ in such a way that $y''=y'-\tau\nu_{y'}$ with
a $\tau>0$. Then $f(y'')=\tau+f(y')>f(y')$.  This is a
contradiction.  Therefore it must hold that $(y'-p)/\vert
y'-p\vert=\nu_{y'}$.  This is nothing but $y'\in {\cal L}(p)$. Since the
existence of $y'$ is clear, this implies ${\cal L}(p)\not=\emptyset$.

\noindent
Step 2.  We prove: if
$(x_0,y_0)\in {\cal M}_2^+(p)\cup {\cal M}_2^{-}(p)\cup {\cal M}_g(p)$, 
then $y_0\in {\cal L}(p)$.

\noindent
By (1) of proposition \ref{Proposition 4.1},
$(y_0-x_0)/\vert y_0-x_0\vert=\nu_{y_0}$.  On the other hand from
(4) of proposition \ref{Proposition 4.1} we have 
$(y_0-x_0)/\vert y_0-x_0\vert=(y_0-p)/\vert y_0-p\vert$.
This gives $(y_0-p)/\vert y_0-p\vert\cdot\nu_{y_0}=1$, that is, 
$y_0\in {\cal L}(p)$.

Now assume that ${\cal L}(p)$ consists of a single point. Let
$(x_0,y_0)\in {\cal M}_2^+(p)\cup {\cal M}_2^{-}(p)\cup {\cal M}_g(p)$. 
From the second
step we have $y_0\in {\cal L}(p)$.  Then the first step implies that
$f(y)$ attains the maximum at $y_0$ only. Using the assumption
that $\partial D$ is $C^2$ at $x_0$ and choosing a suitable half
line that starts at $p$, one can conclude the existence of points
$x_1\in\partial D$, $y_1\in\partial\Omega$ with $y_1\not=y_0$ such
that $x_1$ is on the line determined by $p$ and $y_1$.  Hence
$f(y_1)=l_p(x_1,y_1)$.  By (4) of proposition \ref{Proposition 4.1}, 
we have
$f(y_0)=l_p(x_0,y_0)$.  Since $f(y_1)<f(y_0)$, we obtain
$l_p(x_0,y_0)>l_p(x_1,y_1)$. This is against 
$(x_0,y_0)\in {\cal M}(p)$.
Therefore one gets the desired conclusion.

\noindent
$\Box$

\subsection{Positive definiteness of the Hessian of $\tilde{l_p}(\sigma,\theta)$ at $(\sigma,\theta)=(0,0)$}

In this subsection, we show proposition \ref{Proposition 2.1.}.
Throughout this subsection, we always assume that $D$ is of class $C^2$ and
strictly convex.
 
As is in the proof of proposition \ref{Proposition 4.1}, 
we choose systems of local coordinates
$x=x(\sigma)$, $\sigma=(\sigma_1,\sigma_2)$ with $x_0=x(0)$ and
$y=y(\theta)$, $\theta=(\theta_1,\theta_2)$ with
$y_0=y(0)$ in a neighbourhood of $x_0\in\partial D$ and
$y_0\in\partial\Omega$ respectively. 
It suffices to prove that the Hessian of $\tilde{l}_p(\sigma,\theta)=l_p(x(\sigma), y(\theta))$ at $(\sigma,\theta)=(0,0)$ for
$(x_0,y_0)\in {\cal M}(p)\setminus {\cal M}_{g}(p)$ is positive definite under 
the constraint (\ref{2.9}) on the Weingarten map
for ${\cal A}_{\partial\Omega,y_0}$.  This is equivalent to the statement:
the quadratic form on ${\Bbb R}^2\times{\Bbb R}^2$
$$
\displaystyle
\sum_{j,k=1}^2(\tilde{l_p})_{\sigma_j\sigma_k}(0,0)\xi_j\xi_k
+2\sum_{j,k=1}^2(\tilde{l_p})_{\sigma_j\theta_k}(0,0)\xi_j\eta_k
+\sum_{j,k=1}^2(\tilde{l_p})_{\theta_j\theta_k}(0,0)\eta_j\eta_k,\,\,
(\xi,\eta)\in{\Bbb R}^2\times{\Bbb R}^2
$$
is positive definite.

First we give an expression for the form by using the Weingarten maps for surfaces.
Using (\ref{4.1})-(\ref{4.4}), one gets
\begin{align}
\displaystyle
(\tilde{l_p})_{\sigma_j\sigma_k}(0,0)
&=\nabla_x l_p(x_0,y_0)\cdot\frac{\partial^2 x}{\partial\sigma_j\partial\sigma_k}(0)
+\nabla_x\nabla_xl_p(x_0,y_0)\frac{\partial x}{\partial\sigma_k}(0)
\cdot\frac{\partial x}{\partial\sigma_j}(0)
\label{4.9}
\\
(\tilde{l_p})_{\sigma_j\theta_k}(0,0)
&=-\frac{1}{\vert x_0-y_0\vert}\frac{\partial x}{\partial\sigma_j}(0)\cdot\frac{\partial y}{\partial\theta_k}(0)
\label{4.10}
\\[-3mm]
\intertext{and \vskip-3mm}
(\tilde{l_p})_{\theta_j\theta_k}(0,0)
&=\nabla_yl_p(x_0,y_0)\cdot\frac{\partial^2 y}{\partial\theta_j\partial\theta_k}(0)
+\frac{1}{\vert x_0-y_0\vert}
\frac{\partial y}{\partial\theta_j}(0)\cdot\frac{\partial y}{\partial\theta_k}(0).
\label{4.11}
\end{align}

First we consider
\par\noindent
{\bf Case 1.  $(x_0,y_0)\in {\cal M}_1(p)$.}

\noindent
Let $S$ and $\tilde S$ be an spheroid and a sphere defined by
$$\displaystyle
S=\{x\in{\Bbb R}^3\,\vert\,l_p(x,y_0)=l_p(x_0,y_0)\}
\,\,\text{and}\,\,
\tilde{S}=\{y\in{\Bbb R}^3\,\vert\,\vert y-x_0\vert=\vert y_0-x_0\vert\}.
$$
We denote by $\hat{\nu}_x$, $\hat{\nu}_y$ the
unit outward normal vectors at $x\in S$, $y\in\tilde{S}$ on $S$, $\tilde{S}$, respectively.
Note that
$x_0\in\partial D\cap S$ and $y_0\in\partial\Omega\cap\tilde{S}$.

Here we claim:
\begin{align}
\displaystyle
\vert \nabla_xl_p(x_0,y_0)\vert&=\frac{2(p-x_0)\cdot\nu_{x_0}}{\vert x_0-p\vert}\not=0,
\quad
\displaystyle
\nu_{x_0}=-\frac{1}{\vert\nabla_x l_p(x_0,y_0)\vert}\nabla_x l_p(x_0,y_0)=-\hat{\nu}_{x_0},
\label{4.12}
\\
\displaystyle
\nu_{y_0}&=\hat{\nu}_{y_0}=\nabla_yl_p(x_0,y_0).
\nonumber
\end{align}
The first and second equations come from 
(\ref{like the law of geometrical optics}) in the proof of
(2) of proposition \ref{Proposition 4.1}.  The third equation
is nothing but (1) of proposition \ref{Proposition 4.1}.
(\ref{4.12}) yields that
$T_{x_0}(S)=T_{x_0}(\partial D)$ and 
$T_{y_0}(\tilde{S})=T_{y_0}(\partial\Omega)$.
Then one can choose a local coordinate system $x=x(\varphi)$ with 
$x_0=x(0)$ of $S$
in such a way that
$$\displaystyle
\frac{\partial x}{\partial\sigma_j}(0)=\frac{\partial x}{\partial\varphi_j}(0),\,\,j=1,2.
$$

Since $l_p(x(\varphi),y_0)=l_p(x_0,y_0)$, we have
$(\partial_{\varphi_j}\partial_{\varphi_k})l_p(x(\varphi),y_0)=0$.  That is,
\begin{equation}
\begin{array}{c}
\displaystyle
-\nabla_x l_p(x_0,y_0)\cdot\frac{\partial^2 x}{\partial\varphi_j\partial\varphi_k}(0)
=\nabla_x\nabla_x l_p(x_0,y_0)\frac{\partial x}{\partial\varphi_j}(0)
\cdot\frac{\partial x}{\partial\varphi_k}(0)\\
\\
\displaystyle
=\nabla_x\nabla_x l_p(x_0,y_0)\frac{\partial x}{\partial\sigma_j}(0)\cdot\frac{\partial x}{\partial\sigma_k}(0).
\end{array}
\label{4.13}
\end{equation}
Then from (\ref{4.9}), the second equation 
in (\ref{4.12}) and (\ref{4.13}) we obtain
\begin{equation}
(\tilde{l_p})_{\sigma_j\sigma_k}(0,0)
=-\vert\nabla_x l_p(x_0,y_0)\vert
\left\{\nu_{x_0}\cdot\frac{\partial^2 x}{\partial\sigma_j\partial\sigma_k}(0)
+\hat{\nu}_{x_0}\cdot\frac{\partial^2 x}{\partial\varphi_j\partial\varphi_k}(0)\right\}.
\label{4.14}
\end{equation}
Given $\xi\in{\Bbb R}^2$ set
$$\displaystyle
v(\xi)=\sum_{j=1}^2\xi_j\frac{\partial x}{\partial\sigma_j}(0).
$$
This vector in ${\Bbb R}^3$ belongs to $T_{x_0}(S)=T_{x_0}(\partial D)$.
Since we have
$$\begin{array}{c}
\displaystyle
{\cal A}_{\partial D, x_0}v(\xi)\cdot v(\xi)
=-\sum_{j,k=1}^2\nu_{x_0}\cdot\frac{\partial^2 x}{\partial\sigma_j\partial\sigma_k}(0)\xi_j\xi_k,\\
\\
\displaystyle
{\cal A}_{S, x_0}v(\xi)\cdot v(\xi)
=-\sum_{j,k=1}^2\hat{\nu}_{x_0}\cdot\frac{\partial^2 x}{\partial\varphi_j\partial\varphi_k}(0)\xi_j\xi_k,
\end{array}
$$
(\ref{4.14}) gives
\begin{equation}
\sum_{j,k=1}^2(\tilde{l_p})_{\sigma_j\sigma_k}(0,0)\xi_j\xi_k
=\vert\nabla_xl_p(x_0,y_0)\vert
\left({\cal A}_{\partial D,x_0}v(\xi)\cdot v(\xi)
+{\cal A}_{S,x_0}v(\xi)\cdot v(\xi)\right).
\label{4.15}
\end{equation}

Given $\eta\in{\Bbb R}^2$ set
$$\displaystyle
\tilde{v}(\eta)=\sum_{j=1}^2\eta_j\frac{\partial y}{\partial\theta_j}(0).
$$
This belongs to $T_{y_0}(\partial\Omega)=T_{y_0}(\tilde{S})$.
From (\ref{4.11}), the third equation 
in (\ref{4.12}) and a similar computation we obtain
\begin{equation}
\begin{array}{c}
\displaystyle
\sum_{j,k=1}^2(\tilde{l_p})_{\theta_j\theta_k}(0,0)\eta_j\eta_k
=-{\cal A}_{\partial\Omega, y_0}\tilde{v}(\eta)\cdot\tilde{v}(\eta)
+{\cal A}_{\tilde{S},y_0}\tilde{v}(\eta)\cdot\tilde{v}(\eta).
\end{array}
\label{4.16}
\end{equation}
And also (\ref{4.10}) gives
\begin{equation}
\sum_{j,k=1}^2(\tilde{l}_p)_{\sigma_j\theta_k}(0,0)\xi_j\eta_k
=-\frac{1}{\vert x_0-y_0\vert}v(\xi)\cdot\tilde{v}(\eta).
\label{4.17}
\end{equation}
Summing (\ref{4.15}), (\ref{4.16}) and (\ref{4.17}) 
up, we obtain the formula
\begin{equation}
\begin{array}{c}
\displaystyle
\sum_{j,k=1}^2(\tilde{l_p})_{\sigma_j\sigma_k}(0,0)\xi_j\xi_k
+2\sum_{j,k=1}^2(\tilde{l_p})_{\sigma_j\theta_k}(0,0)\xi_j\eta_k
+\sum_{j,k=1}^2(\tilde{l_p})_{\theta_j\theta_k}(0,0)\eta_j\eta_k\\
\\
\displaystyle
=\vert\nabla_x l_p(x_0,y_0)\vert
\left({\cal A}_{\partial D,x_0}v(\xi)\cdot v(\xi)
+{\cal A}_{S,x_0}v(\xi)\cdot v(\xi)\right)\\
\\
\displaystyle
-\frac{2}{\vert x_0-y_0\vert}v(\xi)\cdot\tilde{v}(\eta)
+{\cal A}_{\tilde{S},y_0}\tilde{v}(\eta)\cdot\tilde{v}(\eta)
-{\cal A}_{\partial\Omega, y_0}\tilde{v}(\eta)\cdot\tilde{v}(\eta).
\end{array}
\label{4.18}
\end{equation}

In order to prove the positive definiteness of the right hand side of (\ref{4.18})
first we consider the case when $\partial D$ is flat and $\partial\Omega$ is replaced with
a sphere in part.

\begin{Lemma}\label{Lemma 4.1}
Let $\tilde{S}'$ be the sphere centered at $\tilde{p}=y_0-l_p(x_0,y_0)\hat{\nu}_{y_0}$ with radius
$l_p(x_0,y_0)$.  Then, for all $(\xi,\eta)\in{\Bbb R}^2\times{\Bbb R}^2$ we have
$$\begin{array}{c}
\displaystyle
\vert\nabla_x l_p(x_0,y_0)\vert
{\cal A}_{S,x_0}v(\xi)\cdot v(\xi)
-\frac{2}{\vert x_0-y_0\vert}v(\xi)\cdot\tilde{v}(\eta)\\
\\
\displaystyle
+{\cal A}_{\tilde{S},y_0}\tilde{v}(\eta)\cdot\tilde{v}(\eta)
-{\cal A}_{\tilde{S}', y_0}\tilde{v}(\eta)\cdot\tilde{v}(\eta)\ge 0.
\end{array}
$$
\end{Lemma}

{\it\noindent Proof.}
Denote by $\Pi$ the set of all points $x$ such that $(x-x_0)\cdot\nu_{x_0}=0$.
Since $(x_0,y_0)\in {\cal M}_1(p)$, from (2) of 
proposition \ref{Proposition 4.1} one knows that
the points $p$ and $y_0$ are in the half space $(x-x_0)\cdot\nu_{x_0}>0$.
Choose a small neighbourhood $V$ of $y_0$.  Given $x\in\Pi$ and $y\in\tilde{S}'\cap V$
we have $\vert p-x\vert=\vert\tilde{p}-x\vert$ and $\vert\tilde{p}-y\vert=l_p(x_0,y_0)$.
The triangle inequality gives
$$\displaystyle
l_p(x,y)=\vert p-x\vert+\vert x-y\vert=\vert\tilde{p}-x\vert+\vert x-y\vert
\ge \vert\tilde{p}-y\vert=l_p(x_0,y_0).
$$
This yields that the function $l_p(x,y)$ on $\Pi\times(\tilde{S}'\cap V)$
attains the minimum value.  Therefore the Hessian of the local representation of the function
on $\Pi\times(\tilde{S}'\cap V)$ has to be non-negative at $(x_0,y_0)$.
This is nothing but the statement of lemma \ref{Lemma 4.1} 
since ${\cal A}_{\Pi,x_0}=0$.

\noindent
$\Box$

A combination of lemma \ref{Lemma 4.1} 
and (\ref{4.18}) gives
\begin{equation}
\begin{array}{c}
\displaystyle
\sum_{j,k=1}^2(\tilde{l_p})_{\sigma_j\sigma_k}(0,0)\xi_j\xi_k
+2\sum_{j,k=1}^2(\tilde{l_p})_{\sigma_j\theta_k}(0,0)\xi_j\eta_k
+\sum_{j,k=1}^2(\tilde{l_p})_{\theta_j\theta_k}(0,0)\eta_j\eta_k\\
\\
\displaystyle
\ge
\vert\nabla_x l_p(x_0,y_0)\vert
{\cal A}_{\partial D,x_0}v(\xi)\cdot v(\xi)
+\frac{1}{l_p(x_0,y_0)}\vert\tilde{v}(\eta)\vert^2
-{\cal A}_{\partial\Omega, y_0}\tilde{v}(\eta)\cdot\tilde{v}(\eta).
\end{array}
\label{4.19}
\end{equation}
Note that we made use of the fact
${\cal A}_{\tilde{S}',y_0}=(1/l_p(x_0,y_0))I$. Then assumption 
(\ref{2.9}) on
${\cal A}_{\partial\Omega,y_0}$ and strict convexity of $\partial D$ yield 
that the right hand side of (\ref{4.19}) is positive definite. 
This completes the proof
of proposition \ref{Proposition 2.1.} 
in the case when $(x_0,y_0)\in {\cal M}_1(p)$.

To complete the proof, from (3) of proposition \ref{Proposition 4.1},
it suffices to consider the following case:

{\bf Case 2.  $(x_0,y_0)\in {\cal M}^+_2(p)\cup {\cal M}^-_2(p)$.}

In this case (\ref{4.8}) in the proof of 
proposition \ref{Proposition 4.1} holds.  
This gives $\nabla_xl_p(x_0,y_0)=0$.
Then from (\ref{4.9}), (\ref{4.10}) and 
$\nu_{y_0}=(y_0-x_0)/\vert y_0-x_0\vert$ we get
$$\begin{array}{c}
\displaystyle
\sum_{j,k=1}^2(\tilde{l}_p)_{\sigma_j\sigma_k}(0,0)\xi_j\xi_k
=\frac{l_p(x_0,y_0)}{\vert p-x_0\vert\vert y_0-x_0\vert}\left\{v(\xi)\cdot v(\xi)-(v(\xi)\cdot\nu_{y_0})^2\right\},\\
\\
\displaystyle
\sum_{j,k=1}^2(\tilde{l}_p)_{\sigma_j\theta_k}(0,0)\xi_j\eta_k
=-\frac{1}{\vert x_0-y_0\vert}v(\xi)\cdot\tilde{v}(\eta)
\end{array}
$$
and we still have (\ref{4.16}).
Summing those up, we obtain
\begin{equation}
\begin{array}{c}
\displaystyle
\sum_{j,k=1}^2(\tilde{l_p})_{\sigma_j\sigma_k}(0,0)\xi_j\xi_k
+2\sum_{j,k=1}^2(\tilde{l_p})_{\sigma_j\theta_k}(0,0)\xi_j\eta_k
+\sum_{j,k=1}^2(\tilde{l_p})_{\theta_j\theta_k}(0,0)\eta_j\eta_k\\
\\
\displaystyle
=\frac{l_p(x_0,y_0)}
{\vert p-x_0\vert\vert x_0-y_0\vert}
\vert w(\xi)\vert^2
-\frac{2}{\vert x_0-y_0\vert}w(\xi)\cdot\tilde{v}(\eta)\\
\\
\displaystyle
+{\cal A}_{\tilde{S},y_0}\tilde{v}(\eta)\cdot\tilde{v}(\eta)-{\cal A}_{\partial\Omega, y_0}\tilde{v}(\eta)\cdot\tilde{v}(\eta),
\end{array}
\label{4.20}
\end{equation}
where $w(\xi)=v(\xi)-(v(\xi)\cdot\nu_{y_0})\nu_{y_0}$.
Here we note that
$$\begin{array}{l}
\displaystyle
\frac{l_p(x_0,y_0)}
{\vert p-x_0\vert\vert x_0-y_0\vert}
\vert w(\xi)\vert^2
-\frac{2}{\vert x_0-y_0\vert}w(\xi)\cdot\tilde{v}(\eta)
+{\cal A}_{\tilde{S},y_0}\tilde{v}(\eta)\cdot\tilde{v}(\eta)-\frac{1}{l_p(x_0,y_0)}\tilde{v}(\eta)\cdot\tilde{v}(\eta)\\
\\
\displaystyle
=
\frac{l_p(x_0,y_0)}
{\vert p-x_0\vert\vert x_0-y_0\vert}
\vert w(\xi)\vert^2
-\frac{2}{\vert x_0-y_0\vert}w(\xi)\cdot\tilde{v}(\eta)
+\left(\frac{1}{\vert x_0-y_0\vert}-\frac{1}{l_p(x_0,y_0)}\right)\tilde{v}(\eta)\cdot\tilde{v}(\eta)\\
\\
\displaystyle
=\frac{l_p(x_0,y_0)}
{\vert p-x_0\vert\vert x_0-y_0\vert}
\left\vert w(\xi)-\frac{\vert p-x_0\vert}{l_p(x_0,y_0)}\tilde{v}(\eta)\right\vert^2\ge 0.
\end{array}
$$
Thus the right hand side of (\ref{4.20}) becomes
$$\displaystyle
\frac{l_p(x_0,y_0)}
{\vert p-x_0\vert\vert x_0-y_0\vert}
\left\vert w(\xi)-\frac{\vert p-x_0\vert}{l_p(x_0,y_0)}\tilde{v}(\eta)\right\vert^2
+\frac{1}{l_p(x_0,y_0)}\tilde{v}(\eta)\cdot\tilde{v}(\eta)
-{\cal A}_{\partial\Omega,y_0}\tilde{v}(\eta)\cdot\tilde{v}(\eta).
$$
By the assumption on ${\cal A}_{\partial\Omega, y_0}$, for the positive
definiteness of (\ref{4.20}) it suffices to prove that if
$$\displaystyle
w(\xi)-\frac{\vert p-x_0\vert}{l_p(x_0,y_0)}\tilde{v}(\eta)=0,\,\,\tilde{v}(\eta)=0,
$$
then $v(\xi)=0$.  First $\tilde{v}(\eta)=0$ yields $w(\xi)=0$, that is, $v(\xi)=(v(\xi)\cdot\nu_{y_0})\nu_{y_0}$.
Since $v(\xi)\cdot\nu_{x_0}=0$, this yields that $v(\xi)\cdot\nu_{y_0}=0$ or $\nu_{x_0}\cdot\nu_{y_0}=0$.
However, if $\nu_{x_0}\cdot\nu_{y_0}=0$, then (\ref{4.8}) 
gives $(p-x_0)\cdot\nu_{x_0}=0$.  This is against
$(x_0,y_0)\not\in {\cal M}_g(p)$.  Therefore it must hold that $v(\xi)\cdot\nu_{y_0}=0$ and
thus this yields $v(\xi)=0$.
This completes the proof of 
proposition \ref{Proposition 2.1.}.

\noindent
$\Box$

\subsection{A sufficient condition of positive definiteness of 
$ l_p$ 
at $(x_0,y_0)\in {\cal M}_1(p)$}
\label{A sufficient condition of positive definiteness of l_p 
at (x_0,y_0) in {cal M}_1(p)}
For $(x_0,y_0)\in {\cal M}_1(p)$ one can relax 
condition (\ref{2.9}).

\begin{Prop}\label{Proposition 4.2}
Let $(x_0,y_0)\in{\cal M}_1(p)$.  Assume that:  there exists a constant $R>d_0\equiv\vert x_0-y_0\vert$
such that the one of the following holds.

\begin{equation}
\displaystyle
{\cal A}_{\partial\Omega, y_0}\le\frac{1}{R}I,\,\,
\vert\nabla_xl_p(x_0,y_0)\vert({\cal A}_{\partial D, x_0}+{\cal A}_{S,x_0})>\frac{R}{(R-d_0)d_0},
\label{4.21}
\end{equation}

\begin{equation}
\displaystyle
{\cal A}_{\partial\Omega, y_0}<\frac{1}{R}I,\,\,
\vert\nabla_xl_p(x_0,y_0)\vert({\cal A}_{\partial D, x_0}+{\cal A}_{S,x_0})\ge\frac{R}{(R-d_0)d_0}.
\label{4.22}
\end{equation}
Then we have the same conclusion as proposition 2.1.

\end{Prop}

{\it\noindent Proof.}
We start with rewriting (\ref{4.18}):

\begin{equation}
\begin{array}{c}
\displaystyle
\sum_{j,k=1}^2(\tilde{l_p})_{\sigma_j\sigma_k}(0,0)\xi_j\xi_k
+2\sum_{j,k=1}^2(\tilde{l_p})_{\sigma_j\theta_k}(0,0)\xi_j\eta_k
+\sum_{j,k=1}^2(\tilde{l_p})_{\theta_j\theta_k}(0,0)\eta_j\eta_k\\
\\
\displaystyle
=\vert\nabla_x l_p(x_0,y_0)\vert
\left({\cal A}_{\partial D,x_0}v(\xi)\cdot v(\xi)
+{\cal A}_{S,x_0}v(\xi)\cdot v(\xi)\right)\\
\\
\displaystyle
+\frac{1}{R}\vert\tilde{v}(\eta)\vert^2
-{\cal A}_{\partial\Omega, y_0}\tilde{v}(\eta)\cdot\tilde{v}(\eta)+I(\xi,\eta),
\end{array}
\label{4.23}
\end{equation}
where
$$\displaystyle
I(\xi,\eta)
=-\frac{1}{R}\vert\tilde{v}(\eta)\vert^2
-\frac{2}{d_0}v(\xi)\cdot\tilde{v}(\eta)
+{\cal A}_{\tilde{S},y_0}\tilde{v}(\eta)\cdot\tilde{v}(\eta).
$$

From the equation 
$$\displaystyle
{\cal A}_{\tilde{S},y_0}\tilde{v}(\eta)\cdot\tilde{v}(\eta)=d_0^{-1}\vert\tilde{v}(\eta)\vert^2
$$
it follows that
$$\begin{array}{c}
\displaystyle
I(\xi,\eta)
=\frac{R-d_0}{Rd_0}\vert\tilde{v}(\eta)\vert^2-\frac{2}{d_0}v(\xi)\cdot\tilde{v}(\eta)\\
\\
\displaystyle
=\frac{R-d_0}{Rd_0}
\left\vert\tilde{v}(\eta)-\frac{R}{R-d_0}v(\xi)\right\vert^2-\frac{R}{(R-d_0)d_0}\vert v(\xi)\vert^2.
\end{array}
$$
Thus the right-hand side of (\ref{4.23}) becomes
\begin{equation}
\begin{array}{c}
\displaystyle
\vert\nabla_x l_p(x_0,y_0)\vert
\left({\cal A}_{\partial D,x_0}v(\xi)\cdot v(\xi)
+{\cal A}_{S,x_0}v(\xi)\cdot v(\xi)\right)-\frac{R}{(R-d_0)d_0}\vert v(\xi)\vert^2\\
\\
\displaystyle
+\frac{1}{R}\vert\tilde{v}(\eta)\vert^2
-{\cal A}_{\partial\Omega, y_0}\tilde{v}(\eta)\cdot\tilde{v}(\eta)
+\frac{R-d_0}{Rd_0}
\left\vert\tilde{v}(\eta)-\frac{R}{R-d_0}v(\xi)\right\vert^2.
\end{array}
\label{4.24}
\end{equation}
Now it is easy to see that (\ref{4.21}) or 
(\ref{4.22}) ensure the positive definiteness 
of (\ref{4.24}).

\noindent
$\Box$

\subsection{An example covered by proposition 
\ref{Proposition 4.2}}

\indent
Using propositions \ref{Proposition 2.2.} 
and \ref{Proposition 4.2}, we can give 
another example including example \ref{Example 2.1.}. 
%
%

Let $\Omega$ be an open ball centered at the origin $O$ with radius 
$R$ and $p$ be a point outside $\Omega$.
Let $D$ be an open ball centered at $q$ with radius $r$.  
Assume that: $\overline D\subset\Omega$,
that is, $\vert q\vert+r<R$;  the line determined by two points 
$p$ and $q$ passes the origin and $\vert p-q\vert\le\vert p\vert$.
By proposition \ref{Proposition 2.2.} one knows that
${\cal M}_g(p)={\cal M}_2^{+}(p)={\cal M}_2^{-}(p)=\emptyset$ 
and ${\cal M}(p)={\cal M}_1(p)$.
Let $(x_0,y_0)\in {\cal M}_1(p)$.  By (1) of proposition 4.1 one knows that
$y_0-x_0$ and $\nu_{y_0}$ are parallel.  
This yields that $x_0$ has to be on the line
determined by $y_0$ and the origin $O$.  By (2) of 
proposition \ref{Proposition 4.1} one knows
that the angle between $p-x_0$ and $\nu_{x_0}$ coincides with the one between
$y_0-x_0$ and $\nu_{x_0}$.  This yields that $x_0$ has to be on the line
determined by $p$ and $q$.  Then one gets
${\cal M}_1(p)=\{(q+r(p-q)/\vert p-q\vert, 
q+(R-\vert q \vert)(p-q)/\vert p-q\vert)\}$ and 
$l(p,D)=\vert p-x_0\vert+\vert x_0-y_0\vert$.
We point out that the condition
\begin{equation}
\displaystyle
\vert\nabla_xl_p(x_0,y_0)
\vert({\cal A}_{\partial D, x_0}+{\cal A}_{S,x_0})>\frac{R}{(R-d_0)d_0},
\label{4.25}
\end{equation}
is satisfied.  Since ${\cal A}_{\partial\Omega,y_0}=(1/R)I$, we conclude 
that (\ref{4.21}) is satisfied.

The condition (\ref{4.25}) for this example is checked as follows.
For the ellipsoid $S$ with the focal points $p$ and $y_0$, it follows that
\begin{equation}
\displaystyle
{\cal A}_{S,x_0}=\frac{l_0}{\displaystyle 2(l_0-d_0)d_0}I,\,l_0=l_p(x_0,y_0).
\label{4.26}
\end{equation}
The proof of (\ref{4.26}) is given in Appendix B.  Using
(\ref{4.26}), the equations 
${\cal A}_{\partial D, x_0}=(1/r)I$ and 
$\vert\nabla_{x_0}l_p(x_0,y_0)\vert=2$
we know that (\ref{4.25}) is equivalent to 
the condition
\begin{equation}
\displaystyle
\frac{1}{r}+\frac{l_0}{\displaystyle 2(l_0-d_0)d_0}>\frac{R}{2(R-d_0)d_0}.
\label{4.27}
\end{equation}
This condition itself is checked by 
a direct computation, however, we present here the detail
for the convenience of the reader.  
Set $h=\vert p-y_0\vert$.  We have
$l_0=h+2d_0$.  Noting that $d_0=R-(r+\vert q\vert)$, one gets
$$\begin{array}{c}
\displaystyle
\frac{1}{r}
+\frac{l_0}{2(l_0-d_0)d_0}
-\frac{R}{2(R-d_0)d_0}
=\frac{2(l_0-d_0)(R-d_0)d_0+rl_0(R-d_0)-r(l_0-d_0)R}
{2rd_0(l_0-d_0)(R-d_0)}.
\end{array}
$$
The numerator of the right-hand side is written as
$$\begin{array}{c}
\displaystyle
2(l_0-d_0)(R-d_0)d_0+rl_0(R-d_0)-r(l_0-d_0)R
\\
\\
\displaystyle
= (l_0-d_0)(2\vert q\vert +r)d_0+rd_0(R-d_0) > 0
\end{array}
$$
Therefore (\ref{4.27}) is valid.

\subsection{Upper bound of the location of $D$}
\label{Upper bound of the location of D}

It should be pointed out that knowing $l(p, D)$, one 
can obtain an upper bound of the location of $D$. 
For $p \notin \overline{\Omega}$, and 
$y \in \partial\Omega$, 
we put ${\mathcal E}_p(y) = \left\{ x \in{\Bbb R}^3 \,\vert\, 
\vert p - x \vert +\vert x - y \vert\ge l(p, D) \,\right\}$,
${\mathcal E}_p=\cap_{y\in\partial\Omega}\,{\mathcal E}_p(y)$,
${\mathcal R}_p = \left\{ x \in{\Bbb R}^3\,\vert\, 
\vert p - x \vert\ge 
2^{-1}(d_{\partial\Omega}(p)+l(p, D)) \,\right\}$,
$d_{\partial\Omega}(p)=\inf\{\vert y-p\vert\,\vert\,y\in\partial\Omega\}$
and
$d_{{\mathcal E}_p\cap\Omega}(p)
=\inf\{\vert x-p\vert\,\vert\,x\in{\mathcal E}_p\cap\Omega\}$.

\begin{Prop}\label{Proposition 4.3}  
It holds that:
(i)  $\overline D\subset{\mathcal E}_p\cap\Omega$;
(ii) ${\mathcal E}_p\cap\Omega\subset{\mathcal R}_p$
and $2d_{{\mathcal E}_p\cap\Omega}(p)\ge l(p, D)+d_{\partial\Omega}(p)$.
\end{Prop}
{\it\noindent Proof.}
For $ x \in D$, one can find $x_0 \in \partial{D}$ such that $x_0$ is on the segment
connecting $x$ with $p$. The definition of $l(p, D)$ implies that
$l(p, D) \leq \vert p - x_0 \vert + \vert x_0 - y \vert$ 
for any $y \in \partial\Omega$.
If $x_0$, $x$ and $y$ are not on a line, from triangle inequality we have 
$l(p, D) < \vert p - x_0 \vert + \vert x_0 - x \vert + \vert x - y \vert
= \vert p - x \vert+ \vert x - y \vert$. If 
$x_0$, $x$ and $y$ are on a line, then this line should be the line passing 
points $x$ and $p$.
If $y$ is located on the segment $xp$, then we have 
$\vert p - x_0 \vert < \vert p - x \vert$ and
$\vert x_0 - y \vert < \vert x - y \vert$, which also implies that
$l(p, D) < \vert p - x \vert+ \vert x - y \vert$.  
If $y$ is outside of segment $px$, then
we have $l(p, D)\le\vert p-x_0\vert+\vert x_0-y\vert=\vert p-x\vert+\vert x-y\vert$.
Hence $D\subset{\mathcal E}_p(y)$.  
Since ${\mathcal E}_p(y)$ is closed and $\overline D\subset\Omega$, 
we obtain (i).

To show (ii), assume that $x \in {\mathcal E}_p\cap\Omega$. 
From $x \in \Omega$ and 
$p \notin \overline{\Omega}$, there exists $t > 0$ such that
$y_t = p+t(x-p)/\vert{x-p}\vert \in \partial\Omega$. 
Since $x \in {\mathcal E}_p(y_t)$ we have 
$l(p, D)\le \vert p - x \vert + \vert x - y_t \vert
= 2\vert p - x \vert - t = 2\vert p - x \vert - \vert p - y_t \vert$,
which implies that $2\vert p - x \vert\ge l(p, D)+\vert p - y_t \vert
\geq l(p, D)+d_{\partial\Omega}(p)$. 
This yields (ii).

\noindent
$\Box$

\noindent
Hence if we have a set $\Lambda \subset {\Bbb R}^3 \setminus\overline{\Omega}$
such that for every $p \in \Lambda$, $l(p, D)$ can be calculated by formula
(\ref{1.13}), 
then for an arbitrary set $\Gamma\subset\partial\Omega$
we have $\overline D \subset 
\cap_{(p,y) \in \Lambda\times\Gamma}\,{\mathcal E}_p(y)\cap\Omega$.

In example \ref{Example 2.1.}, we also note that 
if we put $\Lambda = \{ p \in {\Bbb R}^3\setminus\overline{\Omega}
\,\vert\, h/2 < r_0 - R/2 \,\}$ we have 
$l(p, D) = h+2(R-r)$ for $p \in \Lambda$.
Hence proposition \ref{Proposition 4.3} implies that 
$$
\overline D \subset \cap_{p \in \Lambda}
\left\{ x \in \Omega \,\vert\, \vert p - x \vert\ge h+(R-r) \,\right\} 
= \left\{ x \in \Omega \,\vert\, \vert q - x \vert\le r \,\right\} =\overline D.
$$
Thus the estimate given in proposition \ref{Proposition 4.3} 
is optimal.
This can be extended as follows:

Assume that $\Omega$ is convex and consider the case that 
there exists a point $(x_0,y_0)\in {\cal M}_1(p)$ 
corresponding to the one-dimensional case 
(i.e. $y_0$ is on the line segment $px_0$). 
In this case, the argument showing (ii) of 
proposition \ref{Proposition 4.3} implies that
$l(p, D) = \vert p - x_0 \vert + \vert x_0 - y_0 \vert
= 2\vert p - x_0 \vert - \vert p - y_0 \vert$.
Note that from convexity of $\Omega$, we can characterize $y_0$
as the unique point $y_{min}(p) \in \partial\Omega$ as the point 
attaining the minimum 
$\min_{y \in \partial\Omega}\vert y - p \vert = d_{\partial\Omega}(p)$, and thus
we have $\vert p - y_0 \vert = d_{\partial\Omega}(p)$.
Indeed, from (1) of proposition \ref{Proposition 4.1}, 
one can know that 
$\nu_{y_0} = \vert y_0 - x_0 \vert^{-1}(y_0 - x_0)$.
 Hence we have $\nu_{y_0} = \vert p - y_0 \vert^{-1}(p - y_0)$.
The convexity of $\Omega$ implies that 
the point $y_1 \in \partial\Omega $ satisfying 
$\nu_{y_1} = \vert p - y_1 \vert^{-1}(p - y_1)$ should be 
coincide with $y_{min}(p)$.
Thus the point $x_0$ is determined by 
$x_0 = p+2^{-1}(d_{\partial\Omega}(p)+l(p, D))
\omega(p)$, 
where $\omega(p) = (y_{min}(p) - p)/\vert y_{min}(p) - p\vert$.
Note also that we have
$$
h_{D}(-\omega(p)) = -p\cdot\omega(p)-2^{-1}(d_{\partial\Omega}(p)+l(p, D)).
$$
Thus in the case corresponding to the one-dimensional case
we can find the value of the support function in the direction 
$-\omega(p)$ like as is in the original enclosure method.

Unfortunately, this equality for $h_{D}(-\omega(p))$ does not always holds. 
Even the estimate 
$h_{D}(-\omega(p)) \geq -p\cdot\omega(p)-2^{-1}(d_{\partial\Omega}(p)+l(p, D))$
may not always be true. 
Note also that even the set 
$\cap_{(p,y)\in \Lambda\times\partial\Omega}\,{\mathcal E}_p(y)$ does 
not always coincide with $D$.
However, instead of lines $\omega(p)\cdot{x} = t$, if we use
${\mathcal E}_p(y)$ for $p \notin \overline{\Omega}$ 
and $y \in \partial\Omega$, from proposition \ref{Proposition 4.3}, 
we can give estimates of $D$.

$$\quad$$

\centerline{{\bf Acknowledgements}}

MI was partially supported by Grant-in-Aid for
Scientific Research (C)(nos 21540162 and 25400155) of Japan  Society for
the Promotion of Science.
MK was partially supported by Grant-in-Aid for
Scientific Research (C)(No. 22540194) of Japan  Society for
the Promotion of Science.

$$\quad$$

%

%
\setcounter{equation}{0}
\appendix
\renewcommand{\theequation}{A.\arabic{equation}}

\section{Proof of lemma \ref{Lemma 5.2}}
The Taylor theorem gives
$$
\displaystyle
f(x)= f(x_0) + \left(A(x)(x - x_0), x - x_0 \right)_{{\Bbb R}^n},
$$
where
$$\displaystyle
\\
A(x) = (A_{ij}(x)), \qquad
A_{ij}(x) = \int_0^1(1-\theta)f_{x_ix_j}(x_0+\theta(x - x_0))d\theta.
$$
Since $f \in C^{2, \alpha_0}(\overline{U})$, we have, for a 
positive constant $C > 0$
$$
\vert A_{ij}(x) - A_{ij}(x_0)\vert\le C\vert x - x_0\vert^{\alpha_0}\qquad
(x \in \overline{U}, i,j = 1, 2, \dots ,n).
$$
By the assumption, $A(x_0) = \displaystyle\frac{1}{2}({\rm
Hess}(f)(x_0)) > 0$. Let $\mu_1 \geq \mu_2 \geq \cdots \geq \mu_n
> 0$ be the eigenvalues of $A(x_0)$. Then there exists an
orthogonal matrix $P$ such that ${}^tPA(x_0)P = {\rm diag }(\mu_1,
\mu_2, \cdots ,\mu_n)$. Define $y = {}^tP(x - x_0)$ and set $B(y) =
A(x) - A(x_0)$, $\tilde\varphi(y) = \varphi(x)$. Then there exist
constants $\delta_1 > 0$, $C > 0$ such that
$$
\vert B_{ij}(y)\vert \leq C\vert y\vert^{\alpha_0}
\qquad (\vert y\vert \leq \delta_1, {i, j = 1,2, \ldots, n}).
$$
Thus we have $\vert\Phi_1(y)\vert< \displaystyle\frac{1}{2}\Phi_0(y)$
$(\vert y\vert \leq \delta_1)$, where
$\Phi_0(y) = \displaystyle\frac{1}{2}\sum_{j = 1}^n\mu_jy_j^2$ and 
$\Phi_1(y) = (B(y)y, y)_{{\Bbb R}^n}$.
\par
Let $\delta>0$.  Then, one can choose a constant $c_1 > 0$ such that
$f(x) \geq f(x_0)+ c_1\,\,(x \in \overline{U}, \vert x  - x_0\vert \geq \delta)$, which yields
$$
\left\vert\int_{\{x \in U \,\vert\, \vert x - x_0\vert 
\geq \delta\}}
e^{-\tau{f(x)}}\varphi(x)dx \right\vert
\leq {C{e^{-\tau{f(x_0)}}}}e^{-c_1\tau}\max_{x \in U}\vert \varphi(x)\vert.
$$
Therefore choosing a suitable function $\psi \in C^\infty_0({\Bbb R}^n)$ with
$\psi(y) = 1$ $(\vert y\vert \leq
\delta_1/3)$, $\psi(y) = 0$ $(\vert y\vert \geq 2\delta_1/3)$, we get
$$
\displaystyle
\int_{U}e^{-\tau{f(x)}}\varphi(x)dx
= e^{-\tau{f(x_0)}}\Big\{
\int_{{\Bbb R}^n}\,e^{-\tau{(\Phi_0(y)+\Phi_1(y))}}
\tilde\varphi(y)\psi(y)dy
+ O(e^{-c_1\tau})\Vert \varphi\Vert_{C(\overline{U})}\Big\}
$$
as  $\tau\longrightarrow\infty$.

Here since $\vert e^X - 1\vert = \vert X\int_0^1e^{\theta{X}}d\theta\vert
\leq \vert X\vert e^{\vert \text{Re}\, X\vert}$,  we have
$$\displaystyle
\vert e^{-\tau\Phi_0(y)}\big(e^{-\tau\Phi_1(y)} - 1\big)\vert
\leq C\vert \tau\vert \vert y\vert^{2+\alpha_0}
e^{-\tau\Phi_0(y)}e^{\tau\vert \Phi_1(y)\vert}
\le C\vert \tau\vert \vert y\vert^{2+\alpha_0}e^{-\tau\Phi_0(y)/2}.
$$
This gives
$$\begin{array}{c}
\displaystyle
\Big\vert
\int_{{\Bbb R}^n}\big\{\,e^{-\tau(\Phi_0(y)+\Phi_1(y))}
- e^{-\tau\Phi_0(y)}\big\}
\tilde\varphi(y)\psi(y)dy\Big\vert
\\
\\
\displaystyle
\leq
C\vert \tau\vert\int_{{\Bbb R}^n}\vert y\vert^{2+\alpha_0}
e^{-\tau\Phi_0(y)/2}
\tilde\varphi(y)\psi(y)dy
\\
\\
\displaystyle
\leq C\vert \tau\vert\Vert \varphi\Vert_{C(\overline{U})}
\int_{{\Bbb R}^n}\vert y\vert^{2+\alpha_0}e^{-\tau\Phi_0(y)/2}dy
\\
\\
\displaystyle
= C\vert \tau\vert\Vert \varphi\Vert_{C(\overline{U})}
\int_{{\Bbb R}^n}\vert y\vert^{\alpha_0+2}e^{-\Phi_0(y)/2}dy
\tau^{-(n+2+\alpha_0)/2}
\leq C'_{\delta_0}\Vert \varphi\Vert_{C(\overline{U})}
\tau^{-(n+\alpha_0)/2}.
\end{array}
$$
Summing up, we obtain
$$\displaystyle
\int_{U}e^{-\tau{f(x)}}\varphi(x)dx
= e^{-\tau{f(x_0)}}\Big\{
\int_{{\Bbb R}^n}\,e^{-\tau{\Phi_0(y)}}
\tilde\varphi(y)\psi(y)dy
+ O(\tau^{-(n+\alpha_0)/2})\Vert \varphi\Vert_{C(\overline{U})}\Big\}
$$
as $\tau\longrightarrow\infty$.

Since $\varphi \in C^{0,\alpha_0}(\overline U)$, we have
$\vert \tilde\varphi(y) - \tilde\varphi(0)\vert
\leq \Vert \varphi\Vert_{C^{0, \alpha_0}(\overline{U})}\vert y\vert^{\alpha_0}$
$(\vert y\vert \leq \delta_1)$.
Using a similar argument, one can replace $\tilde{\varphi}(y)$ in
the integrand above with $\tilde{\varphi}(0)=\varphi(x_0)$ and obtain
$$
\displaystyle
\int_{U}e^{-\tau{f(x)}}\varphi(x)dx
= e^{-\tau{f(x_0)}}\Big\{\varphi(x_0)
\int_{{\Bbb R}^n}\,e^{-\tau{\Phi_0(y)}}\psi(y)dy
+ O(\tau^{-(n+\alpha_0)/2})
\Vert \varphi\Vert_{C^{0, \alpha_0}(\overline{U})}\Big\}
$$
as $\tau\longrightarrow\infty$.
Using the asymptotics
$$\begin{array}{c}
\displaystyle
\int_{{\Bbb R}^n}e^{-\tau{\Phi_0(y)}}\psi(y)dy=
\int_{{\Bbb R}^n}e^{-\tau{\Phi_0(y)}}dy
+ \int_{{\Bbb R}^n}e^{-\tau{\Phi_0(y)}}(\psi(y) - 1)dy
\\
\\
\displaystyle
= \frac{1}{\sqrt{\mu_1\mu_2\cdots\mu_n}}
\Big(\frac{2\pi}{\tau}\Big)^{n/2}\big\{1+
O(e^{-\tau(\mu_n\delta_1^2/36)})\big\}
\end{array}
$$
and the equality $\det({\rm Hess}(f)(x_0)) = \mu_1\mu_2\cdots\mu_n$,
we obtain the desired asymptotic formula.

For the estimation choose a cut-off function $\psi \in C^\infty_0(U)$ such that
$\psi = 1$ near $\{x_0\}$, $ \psi \geq 0$ in $U$.
Then, there exists a constant $c_2>0$ such that
$f(x)\ge f(x_0)+c_2\,(x\in\text{supp}\,(1-\psi))$.
This gives
\begin{equation}
\begin{array}{c}
\displaystyle
\left\vert\int_{U}e^{-\tau{f(x)}}\varphi(x)dx \right\vert
\le\left\vert\int_U e^{-\tau f(x)}\varphi(x)\psi(x)dx\right\vert
+\left\vert\int_U e^{-\tau f(x)}(1-\psi(x))\varphi(x)dx\right\vert
\\
\\
\displaystyle
\leq \Vert \varphi\Vert_{C(\overline{U})}
\left\{
\int_{U}e^{-\tau{f(x)}}\psi(x)dx
+Ce^{-\tau(f(x_0)+c_2)}\right\}
\quad (\tau\in\,{\rm \bf C}_{\delta_0}).
\end{array}
\label{A.1}
\end{equation}
Applying the asymptotic formula established above, we have
\begin{equation}
\displaystyle
\int_{U}e^{-\tau{f(x)}}\psi(x)dx =
\frac{e^{-\tau{f(x_0)}}}{\sqrt{\det({\rm Hess}(f)(x_0))}}
\Big(\,\frac{2\pi}{\tau}\Big)^{n/2}\Big(\psi(x_0) +
\Vert \psi\Vert_{C^{0, 1}(\overline{U})}O(\tau^{-1/2})
\Big)
\label{A.2}
\end{equation}
as $\tau\longrightarrow\infty$ uniformly in $\tau \in\,{\rm \bf C}_{\delta_0}$.
Then a combination of (\ref{A.1}) and (\ref{A.2}) 
yields the desired estimate.

\noindent
$\Box$

\section{Proof of (\ref{4.26})}

Set 
$$\begin{array}{c}
\displaystyle
p_0=2^{-1}(p+y_0),\\
\\
\displaystyle
\mbox{\boldmath $e$}_1=(y_0-p)/\vert y_0-p\vert,\,\,
\mbox{\boldmath $e$}'=(x_0-p_0)-((x_0-p_0)\cdot\mbox{\boldmath $e$}_1)\mbox{\boldmath $e$}_1.
\end{array}
$$
If $\mbox{\boldmath $e$}'\not=0$, then set 
$\mbox{\boldmath $e$}_2=\vert\mbox{\boldmath $e$}'\vert^{-1}\mbox{\boldmath $e$}'$ and choose a unit vector 
$\mbox{\boldmath $e$}_3$ in such a way that
$\mbox{\boldmath $e$}_1$, $\mbox{\boldmath $e$}_2$ and $\mbox{\boldmath $e$}_3$ form orthogonal bases of ${\Bbb R}^3$; if $\mbox{\boldmath $e$}'=0$, then choose unit vectors
$\mbox{\boldmath $e$}_2$ and $\mbox{\boldmath $e$}_3$ in such a way that $\mbox{\boldmath $e$}_1$, $\mbox{\boldmath $e$}_2$ and 
$\mbox{\boldmath $e$}_3$ form orthogonal bases of ${\Bbb R}^3$.
Then one can write the equation
for the ellipsoid $S$ as
$$\displaystyle
x=s(\sigma_1,\sigma_2)
=p_0+f(\sigma_1)\mbox{\boldmath $e$}_1+g(\sigma_1)(\cos\,\sigma_2)\mbox{\boldmath $e$}_2+g(\sigma_1)(\sin\,\sigma_2)
\mbox{\boldmath $e$}_3
$$
where 
$$\begin{array}{c}
\displaystyle
f(\sigma_1)=a\cos\,\sigma_1,\,\,g(\sigma_1)=b\sin\,\sigma_1,\\
\\
\displaystyle
a=l_0/2, b=\sqrt{(l_0/2)^2-(\vert p-y_0\vert/2)^2}.
\end{array}
$$
We have
$$\begin{array}{c}
\displaystyle
\frac{\partial s}{\partial\sigma_1}
=f'(\sigma_1)\mbox{\boldmath $e$}_1+g'(\sigma_1)\cos\,(\sigma_2)\mbox{\boldmath $e$}_2+g'(\sigma_1)\sin\,(\sigma_2)
\mbox{\boldmath $e$}_3,\\
\\
\displaystyle
\frac{\partial s}{\partial\sigma_2}
=-g(\sigma_1)\sin\,(\sigma_2)\mbox{\boldmath $e$}_2+g(\sigma_1)\cos\,(\sigma_2)\mbox{\boldmath $e$}_3,\\
\\
\displaystyle
\frac{\partial^2 s}{\partial\sigma_1^2}
=f''(\sigma_1)\mbox{\boldmath $e$}_1+g''(\sigma_1)\cos\,(\sigma_2)\mbox{\boldmath $e$}_2
+g''(\sigma_1)\sin\,(\sigma_2)\mbox{\boldmath $e$}_3,\\
\\
\displaystyle
\frac{\partial^2 s}{\partial\sigma_2^2}
=-g(\sigma_1)\cos\,(\sigma_2)\mbox{\boldmath $e$}_2-g(\sigma_1)\sin\,(\sigma_2)\mbox{\boldmath $e$}_3,\\
\\
\displaystyle
\frac{\partial^2 s}{\partial\sigma_1\partial\sigma_2}
=-g'(\sigma_1)\sin\,(\sigma_2)\mbox{\boldmath $e$}_2+g'(\sigma_1)\cos\,(\sigma_2)\mbox{\boldmath $e$}_3.
\end{array}
$$
Denoting by $\hat{\nu}(\sigma)$ the unit normal vector at $s(\sigma)$ outward to $S$, one gets
$$\displaystyle
\hat{\nu}(\sigma)
=\frac{1}{F(\sigma)}
(g'(\sigma_1)\mbox{\boldmath $e$}_1-f'(\sigma_1)\cos\,(\sigma_2)\mbox{\boldmath $e$}_2
-f'(\sigma_1)\sin\,(\sigma_2)\mbox{\boldmath $e$}_3)
$$
where
$$\displaystyle
F(\sigma)=\sqrt{f'(\sigma_1)^2+g'(\sigma_1)^2}.
$$
A direct computation gives
$$\begin{array}{c}
\displaystyle
\hat{\nu}(\sigma)\cdot\frac{\partial^2 s}{\partial\sigma_1^2}
=\frac{f''(\sigma_1)g'(\sigma_1)-f'(\sigma_1)g''(\sigma_1)}{F(\sigma)},
\,\,
\hat{\nu}(\sigma)\cdot\frac{\partial^2 s}{\partial\sigma_1\partial\sigma_2}=0,\\
\\
\displaystyle
\hat{\nu}(\sigma)\cdot\frac{\partial^2 s}{\partial\sigma_2^2}
=\frac{f'(\sigma_1)g(\sigma_1)}{F(\sigma)}.
\end{array}
$$
Given $\xi=(\xi_1,\xi_2)$ set
$$\displaystyle
v(\xi)=\xi_1\frac{\partial s}{\partial\sigma_1}
+\xi_2\frac{\partial s}{\partial\sigma_2}\in T_{x}S.
$$
We have
\begin{equation}
\begin{array}{c}
\displaystyle
{\cal A}_{S,x_0}(v(\xi))\cdot v(\xi)
=-\sum_{j,k=1}\hat{\nu}(\sigma)\cdot\frac{\partial^2 s}{\partial\sigma_j\partial\sigma_k}(\sigma)\xi_j\xi_k\\
\\
\displaystyle
=(A\xi,\xi)_{{\Bbb R}^2}
\end{array}
\label{A.3}
\end{equation}
where
$$\begin{array}{c}
\displaystyle
A=\frac{1}{F(\sigma)}
\left(\begin{array}{cc}
\displaystyle f'(\sigma_1)g''(\sigma_1)-f''(\sigma_1)g'(\sigma_1) &
\displaystyle 0\\
\\
\displaystyle
0 &
\displaystyle -f'(\sigma_1)g(\sigma_1)
\end{array}\right)\\
\\
\displaystyle
=\frac{ab}{\sqrt{a^2\sin^2\sigma_1+b^2\cos^2\sigma_1}}
\left(\begin{array}{cc}
\displaystyle
1 & 0\\
\\
\displaystyle
0 &\displaystyle \sin^2\sigma_1
\end{array}
\right).
\end{array}
$$
It is easy to see that
\begin{equation}\displaystyle
v(\xi)\cdot v(\xi')
=(G\xi,\xi')_{{\Bbb R}^2},\,\,\xi,\xi'\in{\Bbb R}^2
\label{A.4}
\end{equation}
where
$$\displaystyle
G=\left(\begin{array}{cc}
\displaystyle a^2\sin^2\sigma_1+b^2\cos^2\sigma_1 &
\displaystyle 0\\
\\
\displaystyle
0 & \displaystyle b^2\sin^2\sigma_1
\end{array}
\right).
$$
A combination of (\ref{A.3}) and (\ref{A.4}) gives
$$\displaystyle
{\cal A}_{S,x}v(\xi)=v(G^{-1}A\xi).
$$
This means that the representation matrix of ${\cal A}_{S,x_0}$ with respect to the basis
$\partial s/\partial\sigma_1$ and $\partial s/\partial\sigma_2$ is given by $G^{-1}A$.
Now set $\sigma=(0,0)$.  Since $x_0=p_0+a\mbox{\boldmath $e$}_1$,
$l_0=2d_0+\vert p-y_0\vert$ and $G^{-1}A=(a/b^2)I$, 
we obtain (\ref{4.26}).

\noindent
$\Box$

\section{Proof of (\ref{2.3})}

We put $\epsilon_f(x, \tau) = w(x, \tau)-w_0(x, \tau)$.
From (\ref{2.1-1}) and (\ref{2.2}), it follows that 
$\epsilon_f(\cdot, \tau) \in H^1(\Omega\setminus\overline{D})$
satisfies 
$$
\left\{
\begin{array}{l}
\displaystyle
(\triangle-\tau^2)\epsilon_f = u(x,T)e^{-\tau^2T}\,\,
\text{in}\,\Omega\setminus\overline D,
\\[2mm]
\displaystyle
\frac{\partial \epsilon_f}{\partial\nu}+\rho(x)\epsilon_f=0\,\,\text{on}\,\partial D,
\quad
\displaystyle
\frac{\partial \epsilon_f}{\partial\nu}=0\,\,\text{on}\,\partial\Omega
\end{array}
\right.
$$
in the weak sense. 
Integration by parts implies that
$$
\begin{array}{l}
\displaystyle
\int_{\Omega\setminus\overline{D}}\{\vert\nabla_x\epsilon_f\vert^2
+\tau^2\vert\epsilon_f\vert^2\}dx 
= -e^{-\tau^2T}\int_{\Omega\setminus\overline{D}}u(x,T)\overline{\epsilon_f(x)}dx
+\int_{\partial{D}}\rho\vert\epsilon_f\vert^2dS(x).
\end{array}
$$
This equality yields
\begin{equation}
\begin{array}{l}
\Vert \nabla_x\epsilon_f \Vert^2+\tau^2\Vert \epsilon_f \Vert^2
\leq  
e^{-\tau^2T}\Vert u(\cdot, T) \Vert \Vert \epsilon_f \Vert
+\Vert\rho\Vert_{C^0(\partial{D})}
\Vert \epsilon_f \Vert_{L^2(\partial{D})}^2
(= Q).
\end{array}
\label{A.5}
\end{equation}
where $\Vert \cdot \Vert = \Vert \cdot \Vert_{L^2(\Omega\setminus\overline{D})}$.

Note that there exists a constant 
$C > 0$ depends only on $\partial{D}$ and $\partial\Omega$ 
such that
$$
\Vert \epsilon_f \Vert_{L^2(\partial\Omega)}^2 
+
\Vert \epsilon_f \Vert_{L^2(\partial{D})}^2 
\leq C\left\{\varepsilon\Vert\nabla_x\epsilon_f\Vert^2
+ \frac{1}{\varepsilon}\Vert\epsilon_f\Vert^2\right\}
 \qquad(0 < \varepsilon < 1).
$$
For $\tau > 1$, 
taking $\varepsilon = \tau^{-1}$, we obtain
\begin{equation}
\Vert \epsilon_f \Vert_{L^2(\partial\Omega)}^2 
+
\Vert \epsilon_f \Vert_{L^2(\partial{D})}^2 
\leq C\tau^{-1}\left\{\Vert \nabla_x\epsilon_f \Vert^2
+\tau^2\Vert \epsilon_f \Vert^2
\right\}
\quad(\tau > 1),
\label{A.6}
\end{equation}
which yields
$$
\begin{array}{l}
Q \leq \tau^{-2}
e^{-2\tau^2T}\Vert u(\cdot, T) \Vert^2
+ 4^{-1}\tau^{2}\Vert \epsilon_f \Vert^2
+C\Vert\rho\Vert_{C^0(\partial{D})}
\tau^{-1}
\left\{\Vert \nabla_x\epsilon_f \Vert^2
+\tau^2\Vert \epsilon_f \Vert^2
\right\}.
\end{array}
$$
From the above estimate and (\ref{A.5}), it follows that
there exist constants $C > 0$ and $\mu_0 > 1$ depending
only on $\partial{D}$ and $\rho$ such that
$$
\Vert \nabla_x\epsilon_f \Vert^2
+\tau^2\Vert \epsilon_f \Vert^2
\leq C\tau^{-2}e^{-2\tau^2T}
\Vert u(\cdot, T) \Vert^2
\qquad(\tau \geq \mu_0).
$$ 
Combining the above estimate with (\ref{A.6}), we obtain
$$
\Vert \epsilon_f(\cdot, \tau) \Vert_{L^2(\partial\Omega)} 
\leq C\tau^{-3/2}e^{-\tau^2T}
\quad(\tau \geq \mu_0).
$$
Hence, we get (\ref{2.3}) since 
$$
\begin{array}{c}\displaystyle
I(\tau, p) - I_0(\tau, p)
= \int_{\partial\Omega}\left(\frac{\partial E_{\tau}(y,p)}
{\partial\nu}\epsilon_f(y,\tau)
-\frac{\partial \epsilon_f(y,\tau)}{\partial\nu}E_{\tau}(y,p)\right)dS(y)
\\[4mm]
\displaystyle
= \int_{\partial\Omega}\frac{\partial E_{\tau}(y,p)}
{\partial\nu}\epsilon_f(y,\tau)dS(y)
\qquad\,\,
\end{array}
$$
and 
$$
\left\Vert \frac{\partial E_{\tau}(\cdot,p)}{\partial\nu}
\right\Vert_{L^2(\partial\Omega)}
\leq C\tau 
\qquad(\tau > 0).
$$
\noindent
$\Box$

\end{document}